\pdfoutput=1
\documentclass[12pt,reqno,a4paper]{amsart}

\usepackage[a4paper,left=26mm,right=26mm,top=25mm,bottom=25mm]{geometry}
\usepackage{amsmath,amssymb,amsthm,mathtools,mathrsfs,bm,etoolbox}
\usepackage{xcolor}
\definecolor{LinkBlue}{RGB}{0,51,153}
\definecolor{ReferencePurple}{RGB}{76,29,111}
\usepackage[colorlinks=true,linkcolor=LinkBlue,citecolor=ReferencePurple,urlcolor=ReferencePurple]{hyperref}

\allowdisplaybreaks[2]
\numberwithin{equation}{section}

\makeatletter
\AtBeginDocument{\setlength{\linespacing}{\baselineskip}}
\g@addto@macro\@adjustvertspacing{\setlength{\jot}{5pt}}
\patchcmd{\abstract}{\Small}{\small}{}{}
\renewcommand\section{\@startsection{section}{1}\z@
  {1.1\linespacing\@plus.4\linespacing}{.55\linespacing}%
  {\normalfont\scshape\centering}}
\renewcommand\subsection{\@startsection{subsection}{2}\z@
  {.95\linespacing\@plus.3\linespacing}{.45\linespacing}%
  {\normalfont\normalsize\bfseries}}
\DeclareRobustCommand{\SkipTocEntry}[5]{}
\makeatother

\newtheorem{theorem}{Theorem}[section]
\newtheorem{lemma}[theorem]{Lemma}
\newtheorem{proposition}[theorem]{Proposition}
\theoremstyle{remark}
\newtheorem{remark}[theorem]{Remark}

\title[Smooth forward self-similar Navier--Stokes solutions]{Smooth forward self-similar solutions of the\\ incompressible Navier--Stokes equations\\ in high dimensions}

\author[B. Lai]{Baishun Lai}
\address{LCSM (MOE) and School of Mathematics and Statistics, Hunan Normal University, Changsha 410081, Hunan, P. R. China}
\email{laibaishun@hunnu.edu.cn}

\author[R. Zhang]{Ruiyi Zhang}
\address{School of Mathematical Sciences, Beihang University, Beijing 100191, P. R. China}
\email{ruiyizhang@buaa.edu.cn}

\author[X. Zheng]{Xiaoxin Zheng}
\address{School of Mathematical Sciences, Beihang University, Beijing 100191, P. R. China}
\email{xiaoxinzheng@buaa.edu.cn}


\begin{document}

\begin{abstract}
We construct globally smooth forward self-similar solutions to the incompressible Navier--Stokes equations in \(\mathbb R^n\) (\(n \ge 5\)) for arbitrary \((-1)\)-homogeneous, divergence-free initial data \(\mathbf u_0\) smooth on \(\mathbb S^{n-1}\), attaining \(\mathbf u_0\) strongly in \(H^1_{\mathrm{loc}}(\mathbb R^n)\) as \(t \downarrow 0\).
The profile is obtained as a perturbation of a short-time Navier--Stokes evolution of truncated data whose parabolic scaling defect decays like \((1 + |\boldsymbol y|)^{-n-1}\).
Global smoothness and \((1 + |\boldsymbol y|)^{-n-1}\) remainder decay follow from a modified Bernoulli bound via weighted rotational artificial compressibility, a lift to the stationary Navier--Stokes system in \(\mathbb R^{n+1}\), and \(A_2\)-weighted energy estimates.
\end{abstract}

\maketitle

\tableofcontents

\section{Introduction}\label{sec:intro}

We consider the incompressible Navier--Stokes equations in \(\mathbb R^n\) (\(n \ge 3\)),
\begin{equation}\label{eq:NS}
  \partial_t \mathbf u - \Delta \mathbf u
  + (\mathbf u \cdot \nabla) \mathbf u + \nabla p = \boldsymbol 0,
  \qquad \operatorname{div} \mathbf u = 0.
\end{equation}
System \eqref{eq:NS} is invariant under the parabolic scaling
\begin{equation}\label{eq:NSscaling}
  \mathbf u_\lambda(t, \boldsymbol x)
  = \lambda \mathbf u(\lambda^2 t, \lambda \boldsymbol x),
  \qquad
  p_\lambda(t, \boldsymbol x)
  = \lambda^2 p(\lambda^2 t, \lambda \boldsymbol x),
  \qquad \lambda > 0.
\end{equation}
A solution invariant under \eqref{eq:NSscaling} for all \(\lambda > 0\) is called forward self-similar and takes the form
\begin{equation}\label{eq:similarity}
  \mathbf u(t, \boldsymbol x)
  = \frac{1}{\sqrt t} \mathbf U\left(\frac{\boldsymbol x}{\sqrt t}\right),
  \qquad
  p(t, \boldsymbol x)
  = \frac{1}{t} \Pi\left(\frac{\boldsymbol x}{\sqrt t}\right),
  \qquad t > 0;
\end{equation}
its initial velocity \(\mathbf u_0\) is necessarily \((-1)\)-homogeneous on \(\mathbb R^n \setminus \{\boldsymbol 0\}\), meaning that \(\mathbf u_0(\lambda \boldsymbol x) = \lambda^{-1} \mathbf u_0(\boldsymbol x)\) for all \(\lambda > 0\).
In terms of the similarity variable \(\boldsymbol y = \boldsymbol x / \sqrt t\), substituting \eqref{eq:similarity} into \eqref{eq:NS} shows that the profile pair \((\mathbf U, \Pi)\) satisfies the stationary system on \(\mathbb R^n\)
\begin{equation}\label{eq:profile}
  -\Delta \mathbf U - \frac{1}{2} \boldsymbol y \cdot \nabla \mathbf U - \frac{1}{2} \mathbf U
  + (\mathbf U \cdot \nabla) \mathbf U + \nabla \Pi = \boldsymbol 0,
  \qquad \operatorname{div} \mathbf U = 0,
\end{equation}
while sending \(t \downarrow 0\) at fixed \(\boldsymbol x \ne \boldsymbol 0\) corresponds formally to the far-field asymptotics \(\mathbf U(\boldsymbol y) \sim \mathbf u_0(\boldsymbol y)\) as \(|\boldsymbol y| \to \infty\).

For sufficiently small \((-1)\)-homogeneous initial data, smooth forward self-similar solutions can be obtained by fixed-point arguments in scale-invariant spaces.
Although Kato's \(L^n(\mathbb R^n)\) framework \cite{Kato1984} is scale-critical, nontrivial \((-1)\)-homogeneous fields lie only in weak \(L^n(\mathbb R^n) = L^{n,\infty}(\mathbb R^n)\).
To accommodate homogeneous data, Cannone, Meyer, and Planchon \cite{CannoneMeyerPlanchon1994} and Cannone and Planchon \cite{CannonePlanchon1996} introduced critical Besov and weak-type spaces in \(\mathbb R^3\), Barraza \cite{Barraza1996} constructed small self-similar solutions in \(L^{n,\infty}(\mathbb R^n)\), Koch and Tataru \cite{KochTataru2001} proved small-data well-posedness in \(\mathrm{BMO}^{-1}(\mathbb R^n)\), and Brandolese \cite{Brandolese2009} established fine spatial asymptotics for small self-similar profiles.

For large \((-1)\)-homogeneous data, perturbative contraction arguments no longer apply.
In dimension three, Jia and \v{S}ver\'ak \cite{JiaSverak2014} combined local-in-space estimates near the initial time with Leray--Schauder degree theory to construct smooth forward self-similar solutions for arbitrary \((-1)\)-homogeneous data smooth on \(\mathbb S^2\).
They further showed in \cite{JiaSverak2015} that spectral instability and non-uniqueness of large self-similar profiles give a mechanism for non-uniqueness of Leray--Hopf weak solutions \cite{Leray1934}, which was established for the forced Navier--Stokes equations by Albritton, Bru\'e, and Colombo \cite{AlbrittonBrueColombo2022}.
See also Gruji\'c \cite{Grujic2006} on interior regularity of self-similar suitable weak solutions, Lai, Miao, and Zheng \cite{LaiMiaoZheng2019,LaiMiaoZheng2021}, Lai \cite{Lai2024}, and Miao and Zheng \cite{MiaoZheng2024} for fractional and critical dissipation, and Jia, \v{S}ver\'ak, and Tsai \cite{JiaSverakTsai2018} for a survey of self-similar Navier--Stokes solutions.
By energy methods, Tsai \cite{Tsai2014}, Korobkov and Tsai \cite{KorobkovTsai2016}, Bradshaw and Tsai \cite{BradshawTsai2017,BradshawTsai2019}, and Chae and Wolf \cite{ChaeWolf2018} constructed self-similar and discretely self-similar local Leray solutions \cite{LemarieRieusset2002} in \(\mathbb R^3\) and \(\mathbb R^3_+\).

In dimensions \(3 \le n \le 5\), Lai, Miao, and Zheng \cite[Theorems~1.2 and~1.5]{LaiMiaoZheng2026} recently constructed forward self-similar weak solutions by perturbing around the linear heat profile \(e^\Delta \mathbf u_0\) in \(H^1(\mathbb R^n)\), and further established global smoothness and spatial decay in the energy-critical dimension \(n = 4\).
For \(n \ge 5\), a nontrivial \((-1)\)-homogeneous datum \(\mathbf u_0\) smooth on \(\mathbb S^{n-1}\) belongs to \(H^1_{\mathrm{loc}}(\mathbb R^n)\), whereas the stationary profile system \eqref{eq:profile} is \(H^1\)-supercritical.
Furthermore, the nonlinear residual \(-\mathbb P \operatorname{div}(e^\Delta \mathbf u_0 \otimes e^\Delta \mathbf u_0)\) of the heat profile (where \(\mathbb P\) denotes the Leray projection on \(\mathbb R^n\)) is in general only \(O(|\boldsymbol y|^{-3})\) as \(|\boldsymbol y| \to \infty\), which does not ensure square-integrability when \(n \ge 6\) \cite{BradshawPhelps2023,Brandolese2009,BrandoleseVigneron2007}.
Consequently, global regularity of large-data profiles for \(n = 5\) and both their existence and regularity for \(n \ge 6\) have remained open.

Writing \(B_\ell := \{\boldsymbol y \in \mathbb R^n \mid |\boldsymbol y| < \ell\}\) (\(\ell > 0\)) and \(\langle \boldsymbol y \rangle := (1 + |\boldsymbol y|^2)^{1/2}\), our main result is the following.

\begin{theorem}\label{thm:main}
Let \(n \ge 5\), and let \(\mathbf u_0\) be a \((-1)\)-homogeneous vector field on \(\mathbb R^n \setminus \{\boldsymbol 0\}\) satisfying \(\mathbf u_0|_{\mathbb S^{n-1}} \in C^\infty(\mathbb S^{n-1}; \mathbb R^n)\) and \(\operatorname{div} \mathbf u_0 = 0\) in \(\mathscr{D}'(\mathbb R^n)\).
Then the profile system \eqref{eq:profile} admits a smooth solution \((\mathbf U, \Pi) \in C^\infty(\mathbb R^n; \mathbb R^n) \times C^\infty(\mathbb R^n)\) with the following properties:
\begin{enumerate}
\item[\textup{(i)}]
There exists a constant \(C > 0\) such that for all \(\boldsymbol y \in \mathbb R^n\) and every \(q \in (n / 2, \infty)\),
\begin{equation}\label{eq:maindecay}
  |\mathbf U(\boldsymbol y)| \le \frac{C}{\langle \boldsymbol y \rangle},
  \qquad
  \Pi \in L^q(\mathbb R^n),
\end{equation}
and moreover \(|\mathbf U(\boldsymbol y) - \mathbf u_0(\boldsymbol y)| \le C |\boldsymbol y|^{-3}\) for all \(|\boldsymbol y| \ge 1\).

\item[\textup{(ii)}]
The pair \((\mathbf u, p)\) defined by \eqref{eq:similarity} belongs to \(C^\infty((0, \infty) \times \mathbb R^n)\) and solves \eqref{eq:NS} classically for \(t > 0\).
As \(t \downarrow 0\), for every \(\ell > 0\),
\begin{equation}\label{eq:initialtrace}
  \mathbf u(t) \longrightarrow \mathbf u_0
  \quad \text{strongly in } H^1(B_\ell; \mathbb R^n).
\end{equation}
Moreover, \((\mathbf u, p)\) solves \eqref{eq:NS} distributionally on \([0, \infty) \times \mathbb R^n\) and satisfies the local energy equality up to \(t = 0\) (in particular, it is a local Leray solution attaining \(\mathbf u_0\)).
\end{enumerate}
\end{theorem}

The profile in Theorem~\ref{thm:main} is constructed as a perturbation of a short-time Navier--Stokes evolution of regularized data.
For \(R \ge 1\), we regularize \(\mathbf u_0\) near the origin into a smooth divergence-free datum \(\mathbf u_{0,R}\) vanishing on \(B_R\) and coinciding with \(\mathbf u_0\) outside \(B_{2R}\).
For sufficiently large \(R\), \(\mathbf u_{0,R}\) generates a bounded mild solution \(\mathbf u_R\) of \eqref{eq:NS} on \([0, 1] \times \mathbb R^n\), whose parabolic scaling defect (the infinitesimal scaling derivative at \(\lambda = 1\)) is given for \(t > 0\) by
\begin{equation}\label{eq:scalingdefect}
  \mathbf S_R(t, \boldsymbol y)
  := \left.\frac{\mathrm d}{\mathrm d \lambda}\right|_{\lambda=1}
     \lambda \mathbf u_R(\lambda^2 t, \lambda \boldsymbol y)
  = 2t \partial_t \mathbf u_R(t, \boldsymbol y)
    + \boldsymbol y \cdot \nabla \mathbf u_R(t, \boldsymbol y)
    + \mathbf u_R(t, \boldsymbol y).
\end{equation}
We write \(\|\cdot\|_2 := \|\cdot\|_{L^2(\mathbb R^n)}\), \(H^1_\sigma(\mathbb R^n) := \{\mathbf w \in H^1(\mathbb R^n; \mathbb R^n) \mid \operatorname{div} \mathbf w = 0\}\), and \(\mathcal R_i = \partial_i (-\Delta)^{-1/2}\) (\(1 \le i \le n\)) for the Riesz transforms on \(\mathbb R^n\).

\begin{theorem}\label{thm:structure}
Under the hypotheses of Theorem~\ref{thm:main}, there exists \(R_0 \ge 1\) such that for every \(R \ge R_0\), there exist a constant \(C_R > 0\) (depending on \(n\), \(\mathbf u_0\), and \(R\)) and a smooth divergence-free cutoff datum \(\mathbf u_{0,R} \in C^\infty(\mathbb R^n; \mathbb R^n)\) with \(\mathbf u_{0,R}|_{B_R} = \boldsymbol 0\) and \(\mathbf u_{0,R}|_{\mathbb R^n \setminus B_{2R}} = \mathbf u_0\) generating a unique bounded mild solution \(\mathbf u_R\) to \eqref{eq:NS} on \([0, 1] \times \mathbb R^n\), smooth on \((0, 1] \times \mathbb R^n\), whose scaling defect \eqref{eq:scalingdefect} satisfies
\[
  |\mathbf S_R(1, \boldsymbol y)| \le \frac{C_R}{\langle \boldsymbol y \rangle^{n+1}},
  \qquad \boldsymbol y \in \mathbb R^n,
\]
and, with \(\mathbf B := \mathbf u_R(1)\), the smooth solution \((\mathbf U, \Pi)\) of Theorem~\ref{thm:main} can be chosen of the form \(\mathbf U = \mathbf B + \mathbf W\) with \(\mathbf W \in H^1_\sigma(\mathbb R^n)\) such that:
\begin{enumerate}
\item[\textup{(i)}]
The remainder \(\mathbf W\) satisfies the energy inequality
\begin{equation}\label{eq:mainenergy}
  \|\nabla \mathbf W\|_2^2 + \frac{n - 2}{4} \|\mathbf W\|_2^2
  + \int_{\mathbb R^n} ((\mathbf W \cdot \nabla) \mathbf B) \cdot \mathbf W \, \mathrm d \boldsymbol y
  \le \frac{1}{2} \int_{\mathbb R^n} \mathbf S_R(1) \cdot \mathbf W \, \mathrm d \boldsymbol y,
\end{equation}
the pressure is given by
\begin{equation}\label{eq:mainpressure}
  \Pi = \sum_{i=1}^n \sum_{j=1}^n \mathcal R_i \mathcal R_j (U_i U_j),
  \qquad
  \Pi_B := \sum_{i=1}^n \sum_{j=1}^n \mathcal R_i \mathcal R_j (B_i B_j),
\end{equation}
and, with \(\Phi_B := \Pi_B + \frac{1}{2} |\mathbf B|^2 - \frac{1}{2} \boldsymbol y \cdot \mathbf B\), the modified Bernoulli function satisfies
\begin{equation}\label{eq:mainmax}
  \Phi(\boldsymbol y)
  := \Pi(\boldsymbol y) + \frac{1}{2} |\mathbf U(\boldsymbol y)|^2
     - \frac{1}{2} \boldsymbol y \cdot \mathbf U(\boldsymbol y)
  \le \sup_{\mathbb R^n} \Phi_B < \infty,
  \qquad \boldsymbol y \in \mathbb R^n.
\end{equation}

\item[\textup{(ii)}]
The remainder satisfies the sharp kernel-rate decay
\begin{equation}\label{eq:remdecay}
  |\mathbf W(\boldsymbol y)| \le \frac{C_R}{\langle \boldsymbol y \rangle^{n+1}},
  \qquad \boldsymbol y \in \mathbb R^n.
\end{equation}
\end{enumerate}
\end{theorem}

\begin{remark}\label{rem:main}
\textup{(1)} The smoothness and spatial decay statements in Theorems~\ref{thm:main} and~\ref{thm:structure} apply to the profile constructed below via weighted artificial compressibility, which satisfies the modified Bernoulli bound \eqref{eq:mainmax}; whether every weak solution of \eqref{eq:profile} with \(\mathbf U - \mathbf B \in H^1_\sigma(\mathbb R^n)\) obeying \eqref{eq:mainenergy} is smooth remains open.

\textup{(2)} We assume \(\mathbf u_0|_{\mathbb S^{n-1}} \in C^\infty(\mathbb S^{n-1}; \mathbb R^n)\) only for simplicity of presentation.
Constructing \(\mathbf B\) via the \(k \in \{0, 1\}\) version of Lemma~\ref{lem:evolution} and passing to the incompressible limit require only \(\mathbf u_0|_{\mathbb S^{n-1}} \in C^1(\mathbb S^{n-1}; \mathbb R^n)\).
Lemma~\ref{lem:evolution} estimates derivatives up to order \(k = 2\) and \eqref{eq:bgdiff} uses three derivatives of \(\mathbf u_0\), so with a \(C^3\) cutoff datum \(\mathbf u_{0,R} \in C^3(\mathbb R^n; \mathbb R^n)\), all remaining conclusions of Theorems~\ref{thm:main} and~\ref{thm:structure} for \(\mathbf u_R\), \((\mathbf U, \Pi)\), and \((\mathbf u, p)\) hold for \(\mathbf u_0|_{\mathbb S^{n-1}} \in C^3(\mathbb S^{n-1}; \mathbb R^n)\).

\textup{(3)} For general large data \(\mathbf u_0\), non-uniqueness of forward self-similar profiles is conjectured \cite{AlbrittonBrueColombo2022,JiaSverak2015}; in view of Theorem~\ref{thm:main}\textup{(ii)}, any multiplicity of distinct smooth velocity profiles \(\mathbf U\) satisfying the conclusions of Theorem~\ref{thm:main}\textup{(ii)} for a given \(\mathbf u_0\) in \(n \ge 5\) yields distinct local Leray solutions \cite{LemarieRieusset2002} smooth for \(t > 0\) attaining the same \(H^1_{\mathrm{loc}}(\mathbb R^n)\) initial datum strongly in \(H^1_{\mathrm{loc}}(\mathbb R^n)\).
By contrast, when \(\|\mathbf u_0\|_{L^{n,\infty}(\mathbb R^n)}\) (or \(\|\mathbf u_0\|_{\mathrm{BMO}^{-1}(\mathbb R^n)}\)) is sufficiently small, perturbative fixed-point arguments yield a unique forward self-similar solution in a small ball of the corresponding critical space \cite{Barraza1996,KochTataru2001}.

\textup{(4)} For nontrivial \(\mathbf u_0\) and \(\mathbf W\), dimension \(n \ge 5\) is the sharp threshold both for \(\mathbf u_0 \in H^1_{\mathrm{loc}}(\mathbb R^n)\) and for the vanishing gradient trace \(t^{(n-4)/4} \|\nabla \mathbf W\|_2 \to 0\) in \eqref{eq:initialtrace}.
At \(n = 4\), Proposition~\ref{thm:stationaryreg} still applies to the lifted flow in \(\mathbb R^{n+1} = \mathbb R^5\), while the pressure-potential and exterior estimates in Sections~\ref{sec:profile}--\ref{sec:decay} require four-dimensional modifications (such as \(L^{3/2}\) potential bounds, logarithmic cutoffs on bounded balls, and \(H^1(\mathbb R^4) \hookrightarrow L^4(\mathbb R^4)\); cf.\ \cite{LaiMiaoZheng2026}).
\end{remark}

We now outline the main ideas of the proof.
First, because the linear similarity operator \(-\Delta - \frac{1}{2} \boldsymbol y \cdot \nabla - \frac{1}{2}\) annihilates \(e^\Delta \mathbf u_0\), perturbing around \(e^\Delta \mathbf u_0\) leaves \(-\mathbb P \operatorname{div}(e^\Delta \mathbf u_0 \otimes e^\Delta \mathbf u_0)\) as the forcing term.
However, this term is in general only \(O(|\boldsymbol y|^{-3})\) as \(|\boldsymbol y| \to \infty\) \cite{Brandolese2009}, which does not ensure membership in \(L^2(\mathbb R^n)\) when \(n \ge 6\) because the radial \(L^2\) integral of \(|\boldsymbol y|^{-3}\) diverges as \(\ell \to \infty\):
\begin{equation}\label{eq:dimensionobstacle}
  \int_1^\ell r^{n-7} \, \mathrm d r
  =
  \begin{cases}
    \log \ell, & n = 6, \\
    \dfrac{\ell^{n-6} - 1}{n - 6}, & n > 6,
  \end{cases}
  \qquad \ell > 1.
\end{equation}

We take \(\mathbf B := \mathbf u_R(1)\) and \(\Pi_B := p_R(1)\) from the time-one Navier--Stokes evolution of \(\mathbf u_{0,R}\).
Substituting \(\mathbf U = \mathbf B + \mathbf W\) into \eqref{eq:profile} yields the remainder system
\begin{align}
  &-\Delta \mathbf W - \frac{1}{2} \boldsymbol y \cdot \nabla \mathbf W - \frac{1}{2} \mathbf W
   + ((\mathbf B + \mathbf W) \cdot \nabla) \mathbf W
   + (\mathbf W \cdot \nabla) \mathbf B + \nabla(\Pi - \Pi_B) \nonumber \\
  ={}& -\left[-\Delta \mathbf B - \frac{1}{2} \boldsymbol y \cdot \nabla \mathbf B - \frac{1}{2} \mathbf B
                  + (\mathbf B \cdot \nabla) \mathbf B + \nabla \Pi_B\right]. \label{eq:introresidual}
\end{align}
Termwise, \(\Delta \mathbf B\) and \((\mathbf B \cdot \nabla) \mathbf B + \nabla \Pi_B\) are in general still only \(O(|\boldsymbol y|^{-3})\) at infinity, which does not ensure membership in \(L^2(\mathbb R^n)\) when \(n \ge 6\) by \eqref{eq:dimensionobstacle}.
Since \((\mathbf u_R, p_R)\) solves \eqref{eq:NS} at \(t = 1\), however, the sum \(-\Delta \mathbf B + (\mathbf B \cdot \nabla) \mathbf B + \nabla \Pi_B\) equals \(-\partial_t \mathbf u_R(1)\).
Combining this with the similarity drift reduces the profile residual to the scaling defect \eqref{eq:scalingdefect}:
\begin{equation}\label{eq:introidentity}
  -\Delta \mathbf B - \frac{1}{2} \boldsymbol y \cdot \nabla \mathbf B - \frac{1}{2} \mathbf B
  + (\mathbf B \cdot \nabla) \mathbf B + \nabla \Pi_B
  = -\frac{1}{2} \mathbf S_R(1).
\end{equation}
Differentiating \eqref{eq:NS} with respect to the scaling parameter \(\lambda\) at \(\lambda = 1\) shows that \(\mathbf S_R\) solves the homogeneous linearized system
\begin{equation}\label{eq:introlinearized}
  \partial_t \mathbf S_R - \Delta \mathbf S_R
  + \mathbb P \operatorname{div}
    (\mathbf u_R \otimes \mathbf S_R + \mathbf S_R \otimes \mathbf u_R) = \boldsymbol 0,
  \qquad
  \mathbf S_R(0) = (\boldsymbol y \cdot \nabla + 1) \mathbf u_{0,R}.
\end{equation}

This identity motivates constructing \(\mathbf u_{0,R}\) by cutting off a \(0\)-homogeneous antisymmetric tensor potential for \(\mathbf u_0\), so that \(\mathbf u_{0,R}\) vanishes in \(B_R\) and agrees with \(\mathbf u_0\) outside \(B_{2R}\).
By Euler's homogeneous function theorem \((\boldsymbol y \cdot \nabla + 1) \mathbf u_0 = \boldsymbol 0\), the initial datum \(\mathbf S_R(0)\) is smooth and supported in the annulus \(\{R \le |\boldsymbol y| \le 2R\}\).
Weighted Duhamel estimates for \eqref{eq:introlinearized} then propagate the \(\langle \boldsymbol y \rangle^{-n-1}\) spatial decay of \(e^{\tau\Delta} \mathbb P \operatorname{div}\) to \(\mathbf S_R(1)\) in every dimension \(n \ge 5\), while parabolic rescaling from unit scale ensures \(\|\nabla \mathbf B\|_\infty \le C R^{-2}\), providing coercivity for the linearized operator in \eqref{eq:introresidual} once \(R\) is chosen large enough.

Second, once a weak remainder \(\mathbf W \in H^1_\sigma(\mathbb R^n)\) is obtained, standard elliptic bootstrap arguments break down for \(n \ge 5\).
Indeed, Sobolev embedding gives only \(\mathbf U \in L^{2n/(n-2)}_{\mathrm{loc}}(\mathbb R^n)\), so \((\mathbf U \cdot \nabla) \mathbf U\) belongs only to \(L^{n/(n-1)}_{\mathrm{loc}}(\mathbb R^n)\).
Linear Stokes regularity then returns \(\nabla \mathbf U \in L^{n/(n-2)}_{\mathrm{loc}}(\mathbb R^n)\) with \(n / (n - 2) < 2\), and for \(n \ge 7\) the velocity \(\mathbf U\) need not even belong to \(L^3_{\mathrm{loc}}(\mathbb R^n)\).

Our regularity argument builds on the head-pressure theory for the stationary Navier--Stokes equations \cite{FR94maximum,FR94bounded,FR98criterion,GW78,LY22,Struwe1988} (see also \cite{BangGuiLiuWangXie2025,GLLWX26transverse,GLLWX26homogeneous,Ma26,WY26} for recent high-dimensional stationary constructions).
Specifically, for smooth solutions of \eqref{eq:profile}, the modified Bernoulli function \(\Phi\) in \eqref{eq:mainmax} satisfies
\begin{equation}\label{eq:introPhi}
  -\Delta \Phi + \left(\mathbf U - \frac{1}{2} \boldsymbol y\right) \cdot \nabla \Phi
  = -\frac{1}{2} |\Omega_{\mathbf U}|^2 \le 0,
\end{equation}
where \((\Omega_{\mathbf v})_{ij} = \partial_j v_i - \partial_i v_j\) is the antisymmetric vorticity tensor in the rotational (Lamb) formulation \cite{Lamb1932}.
This differs essentially from the backward Leray equations \cite{Leray1934} studied by Ne\v{c}as, R\r{u}\v{z}i\v{c}ka, and \v{S}ver\'ak \cite{NRS96}, where the outward drift \(+\frac{1}{2} \boldsymbol y\) combined with decay at infinity gives \(\Phi_{\mathrm{back}} \le 0\) and leads to triviality.
Here the forward drift \(-\frac{1}{2} \boldsymbol y\) points inward and \(\boldsymbol y \cdot \mathbf U\) remains bounded at infinity with \(\sup_{\mathbb R^n} \Phi_B \ge 0\); in whole-space barrier arguments, the inward drift contributes a strictly positive zero-order adjoint term \(+\frac{n}{2}\) on \(\mathbb R^n\), allowing us to establish \(\Phi \le \sup_{\mathbb R^n} \Phi_B\) directly across all dimensions \(n \ge 5\).

Constructing large-data weak solutions obeying this bound requires an approximation on \(\mathbb R^n\) that preserves the second-order whole-space maximum principle for \(\Phi\) (cf.\ \cite{FR94maximum,LY22}).
Even for second-order artificial compressibility \cite{Chorin67,Temam68}, a constant penalty parameter on \(\mathbb R^n\) fails to control the far-field divergence cross term produced by the similarity drift \(-\frac{1}{2} \boldsymbol y\).
We therefore approximate \eqref{eq:profile} on \(\mathbb R^n\) by a rotational second-order system with a quadratically weighted divergence penalty:
\[
  -\operatorname{div} \Omega_{\mathbf U_\varepsilon}
  - \frac{1}{2} \Omega_{\mathbf U_\varepsilon} \boldsymbol y
  + \Omega_{\mathbf U_\varepsilon} \mathbf U_\varepsilon
  + \nabla \Phi_\varepsilon = \boldsymbol 0,
  \qquad
  \Phi_\varepsilon := \Phi_B - \left(1 + \frac{\langle \boldsymbol y \rangle^2}{\varepsilon^2}\right) \operatorname{div} \mathbf W_\varepsilon.
\]
Antisymmetry of \(\Omega_{\mathbf U_\varepsilon}\) preserves both convective orthogonality and \(\operatorname{div}(\operatorname{div} \Omega_{\mathbf U_\varepsilon}) = 0\) without \(\operatorname{div} \mathbf W_\varepsilon = 0\), while the weight \(\varepsilon^{-2} \langle \boldsymbol y \rangle^2\) absorbs the drift cross term in energy and gives \(\Phi_\varepsilon \in L^2_{\mathrm{loc}}(\mathbb R^n)\).
Renormalizing via DiPerna--Lions commutators \cite{DL89} and testing against algebraic radial barriers for \(n \ge 5\) then gives \(\Phi_\varepsilon \le \sup_{\mathbb R^n} \Phi_B\) uniformly in \(\varepsilon\), and elliptic recovery of the canonical pressure passes this bound to the limit \(\varepsilon \downarrow 0\).

To deduce smoothness from \(\Phi \le \sup_{\mathbb R^n} \Phi_B\), note that classical stationary regularity theorems apply to unforced divergence-free flows, whereas \(\mathbf U - \frac{1}{2} \boldsymbol y\) has divergence \(-\frac{n}{2}\) on \(\mathbb R^n\).
We restore incompressibility by superimposing \((\mathbf U(\boldsymbol y), 0)\) onto the \((n + 1)\)-dimensional irrotational uniaxial straining flow \((-\frac{1}{2} \boldsymbol y, \frac{n}{2} y_{n+1})\), setting
\[
  \mathbf V(\boldsymbol y, y_{n+1})
  := \left(\mathbf U(\boldsymbol y) - \frac{1}{2} \boldsymbol y, \, \frac{n}{2} y_{n+1}\right),
  \qquad
  \widetilde \Pi(\boldsymbol y, y_{n+1})
  := \Pi(\boldsymbol y) - \frac{1}{8} |\boldsymbol y|^2 - \frac{n^2}{8} y_{n+1}^2.
\]
The stretching component \(+\frac{n}{2} y_{n+1}\) cancels the divergence \(-\frac{n}{2}\) of the horizontal drift, and the quadratic straining terms cancel identically in the total head pressure:
\[
  \widetilde \Pi(\boldsymbol y, y_{n+1}) + \frac{1}{2} |\mathbf V(\boldsymbol y, y_{n+1})|^2
  = \Phi(\boldsymbol y) \le \sup_{\mathbb R^n} \Phi_B.
\]
Thus \((\mathbf V, \widetilde \Pi)\) is an unforced weak solution of the stationary Navier--Stokes equations in \(\mathbb R^{n+1}\).
Quantitative stationary head-pressure regularity \cite{FR94maximum,FR98criterion,LY22} (Proposition~\ref{thm:stationaryreg} in Appendix~\ref{app:stationary}) then yields global smoothness of \((\mathbf U, \Pi)\) for all \(n \ge 5\), and the same lift later provides drift-balanced \(C^n\) bounds at scale \(O(|\boldsymbol y_0|^{-1})\) in the far-field decay analysis.

Finally, inverting the linear similarity operator \(\mathcal L := -\Delta - \frac{1}{2} \boldsymbol y \cdot \nabla - \frac{1}{2}\) via the logarithmic-time heat semigroup \cite{Brandolese2009,CannonePlanchon1996} propagates spatial decay up to the kernel rate \(O(\langle \boldsymbol y \rangle^{-n-1})\) of \(e^{\tau\Delta} \mathbb P \operatorname{div}\) \cite{Brandolese2009,BrandoleseVigneron2007}, provided one first establishes a preliminary algebraic decay bound \(|\mathbf W(\boldsymbol y)| \le C \langle \boldsymbol y \rangle^{-q_0}\) for some \(q_0 > 0\).
Passing from \(\mathbf W \in H^1_\sigma(\mathbb R^n)\) to such an algebraic decay bound is governed by two exponent constraints on spatial weights \(\langle \boldsymbol y \rangle^\beta\).
At distance \(r = |\boldsymbol y_0| \gg 1\), the similarity drift has amplitude \(O(r)\) and balances viscosity at the local scale \(O(r^{-1})\).
Rescaling the lifted stationary flow from balls of radius \(O(r^{-1})\) to the unit ball provides uniform local \(C^k\) bounds for every \(k \ge 1\), and overcoming the local \(L^2\) scaling factor in Gagliardo--Nirenberg interpolation requires \(\langle \boldsymbol y \rangle^\beta \mathbf W \in L^2(\mathbb R^n)\) with \(\beta > n/2\) (specifically \(\beta > (n + 1)/2\) when \(k = n\)).
Conversely, testing the remainder equation against \(\langle \boldsymbol y \rangle^{2\beta} \mathbf W\) transfers a derivative onto the weight in the pressure pairing and requires bounding \(\langle \boldsymbol y \rangle^{\beta - 1} (\Pi - \Pi_B)\) in \(L^2(\mathbb R^n)\).
For the double Riesz transforms in \eqref{eq:mainpressure} and \(\beta \ge 0\), the weight \(\langle \boldsymbol y \rangle^{2\beta - 2}\) belongs to Muckenhoupt's class \(A_2(\mathbb R^n)\) if and only if \(\beta < (n + 2)/2\) \cite{Muckenhoupt1972,Stein1993}.
Both requirements hold whenever
\[
  \frac{n}{2} < \beta < \frac{n + 2}{2}.
\]
To control the nonlinear terms for such \(\beta\), we bound local scale-critical Morrey norms of \(\mathbf U\) on unit balls \(B_1(\boldsymbol y_0)\).
Combining the modified Bernoulli bound \eqref{eq:mainmax} with the uniform local \(L^2\) energy of \(\mathbf U\) improves the direct \(O(\langle \boldsymbol y_0 \rangle^2)\) drift bound to the subquadratic Morrey estimate
\[
  \|\mathbf U\|_{L^{2,n-2}(B_1(\boldsymbol y_0))}^2
  \le C \langle \boldsymbol y_0 \rangle^{2 - 4/n},
\]
which gives the weighted \(H^1\) bound and starts the semigroup iteration to \eqref{eq:maindecay} and \eqref{eq:remdecay}.

The rest of the paper is organized as follows.
In Section~\ref{sec:evolution}, we construct the divergence-free cutoff datum \(\mathbf u_{0,R}\) and establish the \(\langle \boldsymbol y \rangle^{-n-1}\) decay of the scaling defect \(\mathbf S_R(1)\) together with the boundedness of \(\Phi_B\).
In Section~\ref{sec:profile}, we construct the self-similar profile via weighted rotational artificial compressibility, prove the modified Bernoulli upper bound, and deduce global smoothness together with \eqref{eq:initialtrace} by lifting to the stationary Navier--Stokes system in \(\mathbb R^{n+1}\).
Section~\ref{sec:decay} derives the subquadratic local Morrey and weighted \(H^1\) estimates for the remainder, deduces pointwise decay at the \(\langle \boldsymbol y \rangle^{-n-1}\) kernel rate, and completes the proofs of Theorems~\ref{thm:main} and~\ref{thm:structure}.
Appendix~\ref{app:evolution} collects the classical kernel and short-time mild-solution estimates, and Appendix~\ref{app:stationary} gives the quantitative stationary head-pressure interior regularity theorem used in Sections~\ref{sec:profile} and~\ref{sec:decay}.

\addtocontents{toc}{\SkipTocEntry}\subsection*{Notation}

Throughout Sections~\ref{sec:evolution}--\ref{sec:decay} we fix an integer \(n \ge 5\).
For \(\ell > 0\) and \(\boldsymbol y_0 \in \mathbb R^n\), we write \(B_\ell(\boldsymbol y_0) = \{\boldsymbol y \in \mathbb R^n \mid |\boldsymbol y - \boldsymbol y_0| < \ell\}\) (with \(B_\ell = B_\ell(\boldsymbol 0)\)), \(\mathbb S^{n-1} = \partial B_1\), and \(B_\ell^{n+1}\) for the open ball of radius \(\ell\) centered at the origin in \(\mathbb R^{n+1}\).
For a measurable set \(E\), \(\mathbf 1_E\) denotes its indicator function, and \(D_1 \Subset D_2\) means that \(\overline{D_1}\) is compact and contained in \(D_2\).
Latin vector and second-order tensor fields are typeset in bold upright font (\(\mathbf u, \mathbf U, \mathbf B, \mathbf W, \mathbf A, \mathbf F\)), whereas spatial coordinates and Greek vector test functions are typeset in bold italic font (\(\boldsymbol x, \boldsymbol y, \boldsymbol z, \boldsymbol \xi, \boldsymbol \varphi\)).
For vector fields \(\mathbf f, \mathbf g\) and second-order tensor fields \(\mathbf F = (F_{ij}), \mathbf G = (G_{ij})\), we use the componentwise conventions (\(1 \le i, j \le n\))
\[
  (\mathbf f \otimes \mathbf g)_{ij} = f_i g_j, \qquad
  (\nabla \mathbf f)_{ij} = \partial_j f_i, \qquad
  (\operatorname{div} \mathbf F)_i = \sum_{j=1}^n \partial_j F_{ij}, \qquad
  \mathbf F : \mathbf G = \sum_{i=1}^n \sum_{j=1}^n F_{ij} G_{ij},
\]
together with the matrix-vector product \((\mathbf F \mathbf f)_i = \sum_{j=1}^n F_{ij} f_j\) and Frobenius norm \(|\mathbf F| := (\mathbf F : \mathbf F)^{1/2}\); likewise, \(|\nabla^k \mathbf f|\) denotes the Euclidean norm of the \(k\)-th derivative tensor.
For a locally differentiable vector field \(\mathbf v\), its antisymmetric vorticity tensor \(\Omega_{\mathbf v} := \nabla \mathbf v - (\nabla \mathbf v)^\top\) has components
\[
  (\Omega_{\mathbf v})_{ij} = \partial_j v_i - \partial_i v_j,
  \qquad
  (\Omega_{\mathbf v} \mathbf f)_i = \sum_{j=1}^n (\partial_j v_i - \partial_i v_j) f_j,
  \qquad
  |\Omega_{\mathbf v}|^2 = \sum_{i=1}^n \sum_{j=1}^n (\partial_j v_i - \partial_i v_j)^2,
\]
and we write \(s_+ := \max\{s, 0\}\) for \(s \in \mathbb R\) and \(\langle \cdot \rangle := (1 + |\cdot|^2)^{1/2}\) for scalars as well as vectors in \(\mathbb R^n\) or \(\mathbb R^{n+1}\) (so that \(\langle (\boldsymbol y, \sqrt \tau) \rangle = (1 + |\boldsymbol y|^2 + \tau)^{1/2}\) for \(\tau > 0\)).
Lebesgue and Sobolev spaces \(L^p, W^{k,p}, H^k\) with omitted domain are taken over \(\mathbb R^n\), with \(\|\cdot\|_p := \|\cdot\|_{L^p(\mathbb R^n)}\), \(C_0(\mathbb R^n)\) denotes the space of continuous functions on \(\mathbb R^n\) vanishing at infinity, \(\mathscr{D}'(D)\) denotes the space of distributions on an open set \(D\), and vector or tensor target spaces (\(; \mathbb R^n\), \(; \mathbb R^{n\times n}\)) are omitted when clear from context.
For a bounded open set \(D \subset \mathbb R^m\), \(1 \le q < \infty\), and \(0 \le \mu < m\), we recall the Morrey norm
\begin{equation}\label{eq:morreynorm}
  \|\mathbf f\|_{L^{q,\mu}(D)}
  := \sup_{\boldsymbol z \in D,\, 0 < \rho \le \operatorname{diam} D}
     \left(\frac{1}{\rho^\mu} \int_{D \cap B_\rho(\boldsymbol z)} |\mathbf f(\boldsymbol \xi)|^q \, \mathrm d \boldsymbol \xi\right)^{1/q},
\end{equation}
and the same expression with \(\rho > 0\) unrestricted when \(D = \mathbb R^m\).
Generic positive constants \(C\) may vary from line to line and depend on \(n\), the fixed radial cutoff function, and finitely many derivatives of \(\mathbf u_0\) on \(\mathbb S^{n-1}\); once the cutoff radius \(R\) is fixed in Theorem~\ref{thm:structure} and Section~\ref{sec:profile}, constants may also depend on \(R\), while dependence on additional parameters is indicated by subscripts.

\section{Parabolic background and scaling defect}\label{sec:evolution}

We construct the divergence-free regularization \(\mathbf u_{0,R}\) by cutting off an explicit \(0\)-homogeneous antisymmetric tensor potential (cf.\ \cite{BradshawTsai2017,Tsai2014} for annular Bogovskii cutoffs).

Fix a smooth radial function \(\chi \in C^\infty(\mathbb R^n)\) satisfying \(\chi|_{B_1} = 0\) and \(\chi|_{\mathbb R^n \setminus B_2} = 1\), and set \(\chi_R(\boldsymbol y) := \chi(\boldsymbol y / R)\) for \(R \ge 1\).
For \(\boldsymbol y \ne \boldsymbol 0\), define the antisymmetric tensor potential \(\mathbf A = (A_{ij})_{1 \le i,j \le n}\) by
\begin{equation}\label{eq:potential}
  \mathbf A(\boldsymbol y)
  := \frac{\mathbf u_0(\boldsymbol y) \otimes \boldsymbol y
           - \boldsymbol y \otimes \mathbf u_0(\boldsymbol y)}{n - 2},
  \qquad
  A_{ij}(\boldsymbol y)
  = \frac{y_j (\mathbf u_0)_i(\boldsymbol y) - y_i (\mathbf u_0)_j(\boldsymbol y)}{n - 2}.
\end{equation}
Since \(\chi_R\) vanishes on \(B_R\), extending \(\chi_R \mathbf A\) by zero inside \(B_R\) defines a smooth tensor field on \(\mathbb R^n\).

\begin{lemma}\label{lem:cutoff}
Let \(n \ge 5\), and let \(\mathbf u_0\) be a \((-1)\)-homogeneous vector field on \(\mathbb R^n \setminus \{\boldsymbol 0\}\) satisfying \(\mathbf u_0|_{\mathbb S^{n-1}} \in C^\infty(\mathbb S^{n-1}; \mathbb R^n)\) and \(\operatorname{div} \mathbf u_0 = 0\) in \(\mathscr{D}'(\mathbb R^n)\).
For \(R \ge 1\), let \(\chi_R(\boldsymbol y) = \chi(\boldsymbol y / R)\) and let \(\mathbf A\) be given by \eqref{eq:potential}.
Then the vector field \(\mathbf u_{0,R} := \operatorname{div}(\chi_R \mathbf A)\) belongs to \(C^\infty(\mathbb R^n; \mathbb R^n)\), satisfies \(\operatorname{div} \mathbf u_{0,R} = 0\) on \(\mathbb R^n\), and obeys
\begin{equation}\label{eq:cutoffbounds}
  \begin{gathered}
  \mathbf u_{0,R} = \boldsymbol 0 \quad \text{in } B_R, \qquad
  \mathbf u_{0,R} = \mathbf u_0 \quad \text{in } \mathbb R^n \setminus B_{2R}, \\
  \mathbf u_{0,R}(\boldsymbol y)
    = \frac{1}{R} \mathbf u_{0,1}\left(\frac{\boldsymbol y}{R}\right),
  \qquad
  |\nabla^k \mathbf u_{0,R}(\boldsymbol y)|
    \le \frac{C_k}{(R + |\boldsymbol y|)^{k+1}}
    \quad (k \ge 0).
  \end{gathered}
\end{equation}
Moreover, the initial scaling defect is supported in \(\{R \le |\boldsymbol y| \le 2R\}\) and satisfies
\begin{equation}\label{eq:initialdefect}
  (\boldsymbol y \cdot \nabla + 1) \mathbf u_{0,R}
  = \operatorname{div}\bigl((\boldsymbol y \cdot \nabla \chi_R) \mathbf A\bigr)
  \in C_c^\infty(\mathbb R^n; \mathbb R^n).
\end{equation}
\end{lemma}

\begin{proof}
On \(\mathbb R^n \setminus \{\boldsymbol 0\}\), Euler's homogeneous function theorem \(\boldsymbol y \cdot \nabla \mathbf u_0 = -\mathbf u_0\) and \(\operatorname{div} \mathbf u_0 = 0\) yield
\[
  \sum_{j=1}^n \partial_j \bigl[y_j (\mathbf u_0)_i - y_i (\mathbf u_0)_j\bigr]
  = (n - 1) (\mathbf u_0)_i + \boldsymbol y \cdot \nabla (\mathbf u_0)_i
    - y_i \operatorname{div} \mathbf u_0
  = (n - 2) (\mathbf u_0)_i,
\]
so \(\operatorname{div} \mathbf A = \mathbf u_0\) away from the origin.
By the product rule, \(\mathbf u_{0,R} = \chi_R \mathbf u_0 + \mathbf A \nabla \chi_R\), whence \(\mathbf u_{0,R}\) vanishes in \(B_R\) and coincides with \(\mathbf u_0\) on \(\mathbb R^n \setminus B_{2R}\).
Antisymmetry \(A_{ij} = -A_{ji}\) further gives
\[
  \operatorname{div} \mathbf u_{0,R}
  = \sum_{i=1}^n \sum_{j=1}^n \partial_i \partial_j (\chi_R A_{ij}) = 0.
\]

As \(\mathbf A\) is \(0\)-homogeneous, we have
\[
  \chi_R(\boldsymbol y) \mathbf A(\boldsymbol y) = (\chi_1 \mathbf A)\left(\frac{\boldsymbol y}{R}\right),
  \qquad
  \mathbf u_{0,R}(\boldsymbol y) = \frac{1}{R} \mathbf u_{0,1}\left(\frac{\boldsymbol y}{R}\right).
\]
Since \(\mathbf u_{0,1} \in C^\infty(\mathbb R^n; \mathbb R^n)\) agrees with \(\mathbf u_0\) outside \(B_2\), every derivative obeys
\[
  |\nabla^k \mathbf u_{0,1}(\boldsymbol z)| \le \frac{C_k}{(1 + |\boldsymbol z|)^{k+1}},
\]
and rescaling gives
\[
  |\nabla^k \mathbf u_{0,R}(\boldsymbol y)|
  = \frac{1}{R^{k+1}} \left|(\nabla^k \mathbf u_{0,1})\left(\frac{\boldsymbol y}{R}\right)\right|
  \le \frac{C_k}{R^{k+1} \bigl(1 + \frac{|\boldsymbol y|}{R}\bigr)^{k+1}}
  = \frac{C_k}{(R + |\boldsymbol y|)^{k+1}},
\]
proving \eqref{eq:cutoffbounds}.
Finally, using the commutator identity \((\boldsymbol y \cdot \nabla + 1) \operatorname{div} = \operatorname{div}(\boldsymbol y \cdot \nabla)\) together with Euler's homogeneous function theorem \(\boldsymbol y \cdot \nabla \mathbf A = \boldsymbol 0\), we obtain
\[
  (\boldsymbol y \cdot \nabla + 1) \mathbf u_{0,R}
  = (\boldsymbol y \cdot \nabla + 1) \operatorname{div}(\chi_R \mathbf A)
  = \operatorname{div}\bigl[\boldsymbol y \cdot \nabla (\chi_R \mathbf A)\bigr]
  = \operatorname{div}\bigl((\boldsymbol y \cdot \nabla \chi_R) \mathbf A\bigr),
\]
which is supported in \(\operatorname{supp} \nabla \chi_R \subset \{R \le |\boldsymbol y| \le 2R\}\).
\end{proof}

Appendix~\ref{app:evolution} collects the weighted \(L^\infty\) estimates for the heat kernel \(G_\tau\) and the differentiated Oseen kernel \(\mathcal K_\tau\) representing \(e^{\tau\Delta} \mathbb P \operatorname{div}\) (Lemma~\ref{lem:kernel}), together with the unit-scale mild solution \(\mathbf u_1\) solving \eqref{eq:mild} and its parabolic rescaling \(\mathbf u_R\) (Lemma~\ref{lem:evolution}).
We prove the \(\langle \boldsymbol y \rangle^{-n-1}\) decay of the scaling defect \(\mathbf S_R\) first at unit scale by forming parabolic scaling difference quotients of \eqref{eq:mild} and then rescaling to \(R\).

\begin{proposition}\label{prop:defect}
Let \(n \ge 5\), and let \(\mathbf u_0\), \(\mathbf u_{0,R}\), \(T_0 \in (0, 1/4]\), and \(\mathbf u_R\) be as in Lemmas~\ref{lem:cutoff} and~\ref{lem:evolution}.
If \(R \ge 1\) satisfies \(R^{-2} \le T_0\), then the scaling defect \(\mathbf S_R\) defined in \eqref{eq:scalingdefect} satisfies
\begin{equation}\label{eq:defectdecay}
  \sup_{0 < t \le 1} \sup_{\boldsymbol y \in \mathbb R^n}
    \langle \boldsymbol y \rangle^{n+1} |\mathbf S_R(t, \boldsymbol y)|
  \le C_R,
\end{equation}
and hence \(\mathbf S_R(1) \in L^1(\mathbb R^n) \cap L^2(\mathbb R^n) \cap L^\infty(\mathbb R^n)\).
Moreover, for \(\mathbf B = \mathbf u_R(1)\), the functions \(\Pi_B\) and \(\Phi_B\) defined in \eqref{eq:mainpressure}--\eqref{eq:mainmax} belong to \(C^\infty(\mathbb R^n) \cap L^\infty(\mathbb R^n)\) with \(\sup_{\mathbb R^n} \Phi_B \ge 0\), and the background residual identity holds classically and in \(L^2(\mathbb R^n)\):
\begin{equation}\label{eq:bgresidual}
  -\Delta \mathbf B - \frac{1}{2} \boldsymbol y \cdot \nabla \mathbf B - \frac{1}{2} \mathbf B
  + (\mathbf B \cdot \nabla) \mathbf B + \nabla \Pi_B
  = -\frac{1}{2} \mathbf S_R(1).
\end{equation}
\end{proposition}

\begin{proof}
Work first at unit scale \(R = 1\).
For \(t \in [0, T_0]\) and \(\lambda \in (1/2, 2) \setminus \{1\}\), we have \(\lambda^2 t < 4T_0\), and we define the difference quotient
\[
  \mathbf D_\lambda(t, \boldsymbol y)
  := \frac{\lambda \mathbf u_1(\lambda^2 t, \lambda \boldsymbol y) - \mathbf u_1(t, \boldsymbol y)}{\lambda - 1}.
\]
Setting \(\mathbf S_1(0) := (\boldsymbol y \cdot \nabla + 1) \mathbf u_{0,1} \in C_c^\infty(\mathbb R^n; \mathbb R^n)\) and differentiating \(\theta \mapsto \theta \mathbf u_{0,1}(\theta \boldsymbol y)\) gives
\[
  \mathbf D_\lambda(0, \boldsymbol y)
  = \frac{1}{\lambda - 1} \int_1^\lambda \mathbf S_1(0, \theta \boldsymbol y) \, \mathrm d \theta.
\]
By Lemma~\ref{lem:cutoff}, \(\mathbf S_1(0)\) is supported in the annulus \(\{1 \le |\boldsymbol z| \le 2\}\).
Since \(\theta\) lies between \(1\) and \(\lambda \in (1/2, 2)\), it follows that \(\operatorname{supp} \mathbf D_\lambda(0) \subset \overline{B_4}\) and
\[
  \sup_{\boldsymbol y \in \mathbb R^n} \langle \boldsymbol y \rangle^{n+1} |\mathbf D_\lambda(0, \boldsymbol y)|
  \le C
\]
uniformly in \(\lambda \in (1/2, 2) \setminus \{1\}\).

Evaluating the mild integral equation \eqref{eq:mild} at \((\lambda^2 t, \lambda \boldsymbol y)\) with the differentiated Oseen kernel \(\mathcal K_\tau\) from \eqref{eq:kerneldef} and multiplying by \(\lambda\), we write
\[
  \begin{aligned}
  \lambda \mathbf u_1(\lambda^2 t, \lambda \boldsymbol y)
  &= \lambda \int_{\mathbb R^n} G_{\lambda^2 t}(\lambda \boldsymbol y - \boldsymbol z) \, \mathbf u_{0,1}(\boldsymbol z) \, \mathrm d \boldsymbol z \\
  &\quad - \lambda \int_0^{\lambda^2 t} \!\! \int_{\mathbb R^n}
       \mathcal K_{\lambda^2 t - s}(\lambda \boldsymbol y - \boldsymbol z) \,
       (\mathbf u_1 \otimes \mathbf u_1)(s, \boldsymbol z) \, \mathrm d \boldsymbol z \, \mathrm d s.
  \end{aligned}
\]
Changing variables \((s, \boldsymbol z) \mapsto (\lambda^2 s, \lambda \boldsymbol z)\) and using the parabolic scaling identities
\[
  \lambda^n G_{\lambda^2 t}(\lambda \boldsymbol \xi) = G_t(\boldsymbol \xi),
  \qquad
  \lambda^{n+1} \mathcal K_{\lambda^2(t-s)}(\lambda \boldsymbol \xi) = \mathcal K_{t-s}(\boldsymbol \xi),
\]
we obtain the rescaled integral representation
\[
  \begin{aligned}
  \lambda \mathbf u_1(\lambda^2 t, \lambda \boldsymbol y)
  &= \int_{\mathbb R^n} G_t(\boldsymbol y - \boldsymbol z) \, \lambda \mathbf u_{0,1}(\lambda \boldsymbol z) \, \mathrm d \boldsymbol z \\
  &\quad - \int_0^t \!\! \int_{\mathbb R^n}
       \mathcal K_{t-s}(\boldsymbol y - \boldsymbol z) \,
       \bigl[\lambda \mathbf u_1(\lambda^2 s, \lambda \boldsymbol z)
             \otimes \lambda \mathbf u_1(\lambda^2 s, \lambda \boldsymbol z)\bigr]
       \, \mathrm d \boldsymbol z \, \mathrm d s.
  \end{aligned}
\]
Subtracting \eqref{eq:mild} for \(\mathbf u_1\) and dividing by \(\lambda - 1\) yields, on \([0, T_0]\),
\begin{equation}\label{eq:Volterra}
  \mathbf D_\lambda(t)
  = e^{t\Delta} \mathbf D_\lambda(0)
    - \int_0^t \mathcal K_{t-s} * \mathbf F_\lambda(s) \, \mathrm d s,
\end{equation}
where we set
\[
  \mathbf F_\lambda(s, \boldsymbol z)
  := \lambda \mathbf u_1(\lambda^2 s, \lambda \boldsymbol z) \otimes \mathbf D_\lambda(s, \boldsymbol z)
     + \mathbf D_\lambda(s, \boldsymbol z) \otimes \mathbf u_1(s, \boldsymbol z).
\]
By \eqref{eq:unitweighted}, \(\mathbf u_1\) is uniformly bounded on \([0, 4T_0] \times \mathbb R^n\), so
\[
  |\mathbf F_\lambda(s, \boldsymbol z)| \le C |\mathbf D_\lambda(s, \boldsymbol z)|
\]
uniformly in \(\lambda \in (1/2, 2) \setminus \{1\}\).

For each \(\ell \ge 2\), since \(\mathbf D_\lambda \in C([0, T_0]; C_0(\mathbb R^n; \mathbb R^n))\), the truncated weighted supremum
\[
  \mu_{\lambda,\ell}(t)
  := \sup_{\boldsymbol y \in \mathbb R^n} \min\{\langle \boldsymbol y \rangle, \ell\}^{n+1} |\mathbf D_\lambda(t, \boldsymbol y)|
\]
is finite and continuous on \([0, T_0]\).
We split \(\mathbf F_\lambda(s) = \mathbf 1_{B_\ell} \mathbf F_\lambda(s) + \mathbf 1_{\mathbb R^n \setminus B_\ell} \mathbf F_\lambda(s)\).
Because \(\min\{\langle \boldsymbol z \rangle, \ell\}\) is comparable to \(\langle \boldsymbol z \rangle\) on \(B_\ell\) and equals \(\ell\) outside \(B_\ell\), we have
\[
  \sup_{\boldsymbol z \in \mathbb R^n} \langle \boldsymbol z \rangle^{n+1} |(\mathbf 1_{B_\ell} \mathbf F_\lambda)(s, \boldsymbol z)|
  \le C \mu_{\lambda,\ell}(s),
  \qquad
  \|\mathbf 1_{\mathbb R^n \setminus B_\ell} \mathbf F_\lambda(s)\|_\infty
  \le C \ell^{-n-1} \mu_{\lambda,\ell}(s).
\]
In \eqref{eq:Volterra}, we bound \(\min\{\langle \boldsymbol y \rangle, \ell\}^{n+1}\) by \(\langle \boldsymbol y \rangle^{n+1}\) on the linear term and the inner-ball convolution, and by \(\ell^{n+1}\) on the exterior convolution:
\[
  \begin{aligned}
  \min\{\langle \boldsymbol y \rangle, \ell\}^{n+1} |\mathbf D_\lambda(t, \boldsymbol y)|
  &\le \langle \boldsymbol y \rangle^{n+1} \bigl|(e^{t\Delta} \mathbf D_\lambda(0))(\boldsymbol y)\bigr| \\
  &\quad + \int_0^t \langle \boldsymbol y \rangle^{n+1}
           \bigl|\bigl(\mathcal K_{t-s} * (\mathbf 1_{B_\ell} \mathbf F_\lambda(s))\bigr)(\boldsymbol y)\bigr| \, \mathrm d s \\
  &\quad + \int_0^t \ell^{n+1}
           \bigl|\bigl(\mathcal K_{t-s} * (\mathbf 1_{\mathbb R^n \setminus B_\ell} \mathbf F_\lambda(s))\bigr)(\boldsymbol y)\bigr| \, \mathrm d s.
  \end{aligned}
\]
Applying \eqref{eq:kernelweighted} in Lemma~\ref{lem:kernel} with \(q = n + 1\) to the first two terms and \eqref{eq:kernel} to the third gives
\[
  \begin{aligned}
  \langle \boldsymbol y \rangle^{n+1} \bigl|(e^{t\Delta} \mathbf D_\lambda(0))(\boldsymbol y)\bigr|
  &\le C \sup_{\boldsymbol z \in \mathbb R^n} \langle \boldsymbol z \rangle^{n+1} |\mathbf D_\lambda(0, \boldsymbol z)|
   \le C, \\
  \langle \boldsymbol y \rangle^{n+1}
  \bigl|\bigl(\mathcal K_{t-s} * (\mathbf 1_{B_\ell} \mathbf F_\lambda(s))\bigr)(\boldsymbol y)\bigr|
  &\le \frac{C}{\sqrt{t - s}}
       \sup_{\boldsymbol z \in \mathbb R^n} \langle \boldsymbol z \rangle^{n+1} |(\mathbf 1_{B_\ell} \mathbf F_\lambda)(s, \boldsymbol z)|
   \le \frac{C}{\sqrt{t - s}} \mu_{\lambda,\ell}(s), \\
  \ell^{n+1}
  \bigl|\bigl(\mathcal K_{t-s} * (\mathbf 1_{\mathbb R^n \setminus B_\ell} \mathbf F_\lambda(s))\bigr)(\boldsymbol y)\bigr|
  &\le \ell^{n+1} \|\mathcal K_{t-s}\|_1 \|\mathbf 1_{\mathbb R^n \setminus B_\ell} \mathbf F_\lambda(s)\|_\infty
   \le \frac{C}{\sqrt{t - s}} \mu_{\lambda,\ell}(s).
  \end{aligned}
\]
Substituting these bounds and taking the supremum over \(\boldsymbol y \in \mathbb R^n\) yields, on \([0, T_0]\),
\[
  \mu_{\lambda,\ell}(t)
  \le C + C \int_0^t \frac{\mu_{\lambda,\ell}(s)}{\sqrt{t - s}} \, \mathrm d s,
\]
with a constant \(C\) independent of \(\lambda\) and \(\ell\).
Iterating once and applying Fubini's theorem gives
\[
  \begin{aligned}
  \mu_{\lambda,\ell}(t)
  &\le C + C \int_0^t \frac{1}{\sqrt{t - s}}
       \left(C + C \int_0^s \frac{\mu_{\lambda,\ell}(\sigma)}{\sqrt{s - \sigma}} \, \mathrm d \sigma\right) \mathrm d s \\
  &= C + 2 C^2 \sqrt t
     + C^2 \int_0^t \left(\int_\sigma^t \frac{\mathrm d s}{\sqrt{(t - s)(s - \sigma)}}\right)
       \mu_{\lambda,\ell}(\sigma) \, \mathrm d \sigma \\
  &\le C + C \int_0^t \mu_{\lambda,\ell}(\sigma) \, \mathrm d \sigma,
  \end{aligned}
\]
where we used the Beta-integral identity \(\int_\sigma^t (t - s)^{-1/2} (s - \sigma)^{-1/2} \, \mathrm d s = \pi\).
By Gr\"onwall's inequality,
\[
  \sup_{0 \le t \le T_0} \mu_{\lambda,\ell}(t) \le C.
\]
Sending \(\ell \to \infty\) and then \(\lambda \to 1\) pointwise on \((0, T_0] \times \mathbb R^n\), where \(\mathbf u_1\) is smooth and \(\mathbf D_\lambda\) converges to \(\mathbf S_1 := (2t \partial_t + \boldsymbol y \cdot \nabla + 1) \mathbf u_1\), we obtain
\[
  |\mathbf S_1(\tau, \boldsymbol z)| \le \frac{C}{\langle \boldsymbol z \rangle^{n+1}},
  \qquad (\tau, \boldsymbol z) \in (0, T_0] \times \mathbb R^n.
\]
For \(R^{-2} \le T_0\) and \(0 < t \le 1\), rescaling via \eqref{eq:rescale} gives
\[
  |\mathbf S_R(t, \boldsymbol y)|
  = \frac{1}{R} \left|\mathbf S_1\left(\frac{t}{R^2}, \, \frac{\boldsymbol y}{R}\right)\right|
  \le \frac{C}{R \langle \boldsymbol y / R \rangle^{n+1}}
  \le \frac{C_R}{\langle \boldsymbol y \rangle^{n+1}},
\]
establishing \eqref{eq:defectdecay}.

Next, since \(R \ge 1\), \eqref{eq:bgbounds} yields
\[
  |\mathbf B \otimes \mathbf B| + |\nabla(\mathbf B \otimes \mathbf B)|
  \le \frac{C}{\langle \boldsymbol y \rangle^2},
\]
so \(\mathbf B \otimes \mathbf B \in W^{1,n+1}(\mathbb R^n)\).
Boundedness of \(\mathcal R_i\) on \(W^{1,n+1}(\mathbb R^n)\) and the Sobolev embedding \(W^{1,n+1}(\mathbb R^n) \hookrightarrow C_0(\mathbb R^n)\) (by Morrey's inequality and density of \(C_c^\infty(\mathbb R^n)\); see \cite[Corollary~9.14]{Brezis2011}) give \(\Pi_B \in C_0(\mathbb R^n)\).
Because \(\mathbf B\) is smooth, interior elliptic regularity applied to
\[
  -\Delta \Pi_B = \sum_{i=1}^n \sum_{j=1}^n \partial_i \partial_j (B_i B_j) \in C^\infty(\mathbb R^n)
\]
shows that \(\Pi_B \in C^\infty(\mathbb R^n)\), and boundedness of \(\boldsymbol y \cdot \mathbf B\) then gives \(\Phi_B \in C^\infty(\mathbb R^n) \cap L^\infty(\mathbb R^n)\).
Furthermore, since \(\operatorname{div} \mathbf B = 0\), the divergence theorem yields
\[
  \int_{\partial B_\ell} \boldsymbol y \cdot \mathbf B \, \mathrm d \sigma
  = \ell \int_{B_\ell} \operatorname{div} \mathbf B \, \mathrm d \boldsymbol y
  = 0
  \qquad (\ell > 0).
\]
Since \(\Pi_B(\boldsymbol y) + \frac{1}{2} |\mathbf B(\boldsymbol y)|^2 \to 0\) as \(|\boldsymbol y| \to \infty\), we obtain
\[
  \sup_{\mathbb R^n} \Phi_B
  \ge \lim_{\ell \to \infty} \frac{1}{|\partial B_\ell|}
      \int_{\partial B_\ell}
        \left(\Pi_B + \frac{1}{2} |\mathbf B|^2 - \frac{1}{2} \boldsymbol y \cdot \mathbf B\right) \mathrm d \sigma
  = 0.
\]

Finally, setting
\[
  p_R(t) := \sum_{i=1}^n \sum_{j=1}^n \mathcal R_i \mathcal R_j \bigl((\mathbf u_R)_i (\mathbf u_R)_j\bigr),
\]
we see that the smooth \(L^{n+1}\) mild solution \(\mathbf u_R\) solves \eqref{eq:NS} classically on \((0, 1] \times \mathbb R^n\) \cite[Theorem~4]{Giga1986}.
Replacing \(-\Delta \mathbf B + (\mathbf B \cdot \nabla) \mathbf B + \nabla \Pi_B\) by \(-\partial_t \mathbf u_R(1)\) at \(t = 1\) and recalling \(\mathbf S_R(1) = 2 \partial_t \mathbf u_R(1) + \boldsymbol y \cdot \nabla \mathbf B + \mathbf B\) establishes \eqref{eq:bgresidual}.
\end{proof}

\section{Existence and smoothness of the profile}\label{sec:profile}

\subsection{Weighted artificial compressibility}\label{subsec:penalized}

In view of the gradient bound \(\|\nabla \mathbf u_R(1)\|_\infty \le C R^{-2}\) in \eqref{eq:bgbounds}, we henceforth fix a cutoff radius \(R \ge 1\) satisfying \(R^{-2} \le T_0\) large enough that \(\mathbf B := \mathbf u_R(1)\) obeys
\begin{equation}\label{eq:coercivity}
  \|\nabla \mathbf B\|_\infty \le \frac{n - 2}{8},
  \qquad\text{so that}\qquad
  \frac{n - 2}{4} |\mathbf w|^2 + ((\mathbf w \cdot \nabla) \mathbf B) \cdot \mathbf w
  \ge \frac{n - 2}{8} |\mathbf w|^2
\end{equation}
pointwise on \(\mathbb R^n\) for all \(\mathbf w \in \mathbb R^n\).
All constants below may depend on this fixed radius \(R\).

Under the vorticity convention \((\Omega_{\mathbf v})_{ij} = \partial_j v_i - \partial_i v_j\) from Section~\ref{sec:intro}, for any \(\mathbf v \in H^1_{\mathrm{loc}}(\mathbb R^n; \mathbb R^n)\) direct componentwise expansion gives
\[
  \begin{aligned}
  (\Omega_{\mathbf v} \mathbf v)_i
  &= \sum_{j=1}^n (\partial_j v_i - \partial_i v_j) v_j
   = \mathbf v \cdot \nabla v_i - \frac{1}{2} \partial_i |\mathbf v|^2, \\
  -\frac{1}{2} (\Omega_{\mathbf v} \boldsymbol y)_i
  &= -\frac{1}{2} \sum_{j=1}^n (\partial_j v_i - \partial_i v_j) y_j
   = -\frac{1}{2} \boldsymbol y \cdot \nabla v_i
     + \frac{1}{2} \partial_i (\boldsymbol y \cdot \mathbf v) - \frac{1}{2} v_i,
  \end{aligned}
\]
which sum to Lamb's identity \cite{FR94maximum,Lamb1932}
\begin{equation}\label{eq:lambidentity}
  -\frac{1}{2} \Omega_{\mathbf v} \boldsymbol y + \Omega_{\mathbf v} \mathbf v
  = -\frac{1}{2} \boldsymbol y \cdot \nabla \mathbf v - \frac{1}{2} \mathbf v
    + (\mathbf v \cdot \nabla) \mathbf v
    + \nabla\left(\frac{1}{2} \boldsymbol y \cdot \mathbf v - \frac{1}{2} |\mathbf v|^2\right),
\end{equation}
while differentiating \((\Omega_{\mathbf v})_{ij}\) with respect to \(y_j\) yields the viscous identity
\begin{equation}\label{eq:curlcurl}
  -\operatorname{div} \Omega_{\mathbf v}
  = -\Delta \mathbf v + \nabla \operatorname{div} \mathbf v.
\end{equation}
Moreover, by antisymmetry of \(\Omega_{\mathbf v}\), we have in \(\mathscr{D}'(\mathbb R^n)\) (even when \(\operatorname{div} \mathbf v \ne 0\))
\[
  \operatorname{div}(\operatorname{div} \Omega_{\mathbf v})
  = \sum_{i=1}^n \sum_{j=1}^n \partial_i \partial_j (\Omega_{\mathbf v})_{ij} = 0.
\]

In the classical artificial-compressibility method on bounded domains \cite{Chorin67,Temam68}, one relaxes \(\operatorname{div} \mathbf W_\varepsilon = 0\) via a constant-coefficient pressure penalty \(-\varepsilon^{-2} \operatorname{div} \mathbf W_\varepsilon\).
On \(\mathbb R^n\), however, testing the rotational similarity drift \(-\frac{1}{2} \Omega_{\mathbf W_\varepsilon} \boldsymbol y\) against \(\mathbf W_\varepsilon\) when \(\operatorname{div} \mathbf W_\varepsilon \ne 0\) produces the far-field divergence cross term
\[
  -\frac{1}{2} \int_{\mathbb R^n} (\boldsymbol y \cdot \mathbf W_\varepsilon) \operatorname{div} \mathbf W_\varepsilon \, \mathrm d \boldsymbol y.
\]
Since the unweighted energy controls only \(\|\mathbf W_\varepsilon\|_2\) while \(|\boldsymbol y|\) grows linearly, absorbing this cross term by Cauchy--Schwarz requires a quadratically growing penalty weight \(\varepsilon^{-2} \langle \boldsymbol y \rangle^2\).

Accordingly, for \(0 < \varepsilon \le 1\) we set \(\mathbf U_\varepsilon := \mathbf B + \mathbf W_\varepsilon\), define the penalized Bernoulli function
\begin{equation}\label{eq:penaltydef}
  \Phi_\varepsilon
  := \Phi_B - \left(1 + \frac{\langle \boldsymbol y \rangle^2}{\varepsilon^2}\right) \operatorname{div} \mathbf W_\varepsilon,
\end{equation}
and seek \(\mathbf W_\varepsilon\) solving the penalized rotational system in \(\mathscr{D}'(\mathbb R^n)\):
\begin{equation}\label{eq:penalizedlamb}
  -\operatorname{div} \Omega_{\mathbf U_\varepsilon}
  - \frac{1}{2} \Omega_{\mathbf U_\varepsilon} \boldsymbol y
  + \Omega_{\mathbf U_\varepsilon} \mathbf U_\varepsilon
  + \nabla \Phi_\varepsilon = \boldsymbol 0.
\end{equation}
By \eqref{eq:curlcurl} and \(\operatorname{div} \mathbf B = 0\), the unweighted \(\nabla \operatorname{div} \mathbf W_\varepsilon\) terms cancel in \eqref{eq:penalizedlamb}:
\begin{equation}\label{eq:viscouscancel}
  -\operatorname{div} \Omega_{\mathbf U_\varepsilon} + \nabla \Phi_\varepsilon
  = -\Delta \mathbf B + \nabla \Phi_B
    - \Delta \mathbf W_\varepsilon
    - \nabla\left(\frac{\langle \boldsymbol y \rangle^2}{\varepsilon^2} \operatorname{div} \mathbf W_\varepsilon\right).
\end{equation}
Expanding \(\Omega_{\mathbf U_\varepsilon} = \Omega_{\mathbf B} + \Omega_{\mathbf W_\varepsilon}\) in the rotational drift and convective terms gives
\[
  -\frac{1}{2} \Omega_{\mathbf U_\varepsilon} \boldsymbol y
  + \Omega_{\mathbf U_\varepsilon} \mathbf U_\varepsilon
  = \left(-\frac{1}{2} \Omega_{\mathbf B} \boldsymbol y + \Omega_{\mathbf B} \mathbf B\right)
    - \frac{1}{2} \Omega_{\mathbf W_\varepsilon} \boldsymbol y
    + \Omega_{\mathbf W_\varepsilon} (\mathbf B + \mathbf W_\varepsilon)
    + \Omega_{\mathbf B} \mathbf W_\varepsilon,
\]
and applying \eqref{eq:lambidentity} to \(\mathbf v = \mathbf B\) together with \(\Phi_B = \Pi_B + \frac{1}{2} |\mathbf B|^2 - \frac{1}{2} \boldsymbol y \cdot \mathbf B\) and the background residual identity \eqref{eq:bgresidual} yields
\[
  -\Delta \mathbf B + \nabla \Phi_B
  - \frac{1}{2} \Omega_{\mathbf B} \boldsymbol y + \Omega_{\mathbf B} \mathbf B
  = -\Delta \mathbf B - \frac{1}{2} \boldsymbol y \cdot \nabla \mathbf B - \frac{1}{2} \mathbf B
    + (\mathbf B \cdot \nabla) \mathbf B + \nabla \Pi_B
  = -\frac{1}{2} \mathbf S_R(1).
\]
Substituting these two identities and \eqref{eq:viscouscancel} into \eqref{eq:penalizedlamb} shows that \eqref{eq:penalizedlamb} is equivalent in \(\mathscr{D}'(\mathbb R^n)\) to the remainder equation
\begin{equation}\label{eq:rempenalized}
  -\Delta \mathbf W_\varepsilon
  - \frac{1}{2} \Omega_{\mathbf W_\varepsilon} \boldsymbol y
  + \Omega_{\mathbf W_\varepsilon} (\mathbf B + \mathbf W_\varepsilon)
  + \Omega_{\mathbf B} \mathbf W_\varepsilon
  - \nabla\left(\frac{\langle \boldsymbol y \rangle^2}{\varepsilon^2} \operatorname{div} \mathbf W_\varepsilon\right)
  = \frac{1}{2} \mathbf S_R(1).
\end{equation}

\begin{proposition}\label{prop:approx}
Let \(n \ge 5\), let \(\mathbf u_0\) satisfy the hypotheses of Theorem~\ref{thm:main}, and let \(R \ge 1\) with \(R^{-2} \le T_0\) be fixed so that \(\mathbf B = \mathbf u_R(1)\) satisfies \eqref{eq:coercivity}.
Set \(\varepsilon_0 := \bigl[2(1 + \|\mathbf B\|_\infty)\bigr]^{-1} \in (0, 1/2]\).
Then for every \(0 < \varepsilon \le \varepsilon_0\), there exists \(\mathbf W_\varepsilon \in H^1(\mathbb R^n; \mathbb R^n)\) solving \eqref{eq:penalizedlamb} (equivalently \eqref{eq:rempenalized}) in \(\mathscr{D}'(\mathbb R^n)\) and satisfying the uniform bound
\begin{equation}\label{eq:approxbounds}
  \|\mathbf W_\varepsilon\|_{H^1}^2
  + \frac{1}{\varepsilon^2} \int_{\mathbb R^n} \langle \boldsymbol y \rangle^2 |\operatorname{div} \mathbf W_\varepsilon(\boldsymbol y)|^2 \, \mathrm d \boldsymbol y
  \le C,
\end{equation}
so that \(\|\langle \boldsymbol y \rangle \operatorname{div} \mathbf W_\varepsilon\|_2 \le C \varepsilon\), together with the penalized energy inequality
\begin{equation}\label{eq:approxenergy}
  \begin{aligned}
  &\|\nabla \mathbf W_\varepsilon\|_2^2
   + \frac{n - 2}{4} \|\mathbf W_\varepsilon\|_2^2
   + \int_{\mathbb R^n} ((\mathbf W_\varepsilon \cdot \nabla) \mathbf B) \cdot \mathbf W_\varepsilon \, \mathrm d \boldsymbol y
   + \frac{1}{\varepsilon^2} \|\langle \boldsymbol y \rangle \operatorname{div} \mathbf W_\varepsilon\|_2^2 \\
  {}+{}& \int_{\mathbb R^n}
           \left(\left(\mathbf B - \frac{1}{2} \boldsymbol y\right) \cdot \mathbf W_\varepsilon\right)
           \operatorname{div} \mathbf W_\varepsilon \, \mathrm d \boldsymbol y \\
   \le{}& \frac{1}{2} \int_{\mathbb R^n} \mathbf S_R(1) \cdot \mathbf W_\varepsilon \, \mathrm d \boldsymbol y.
  \end{aligned}
\end{equation}
\end{proposition}

\begin{proof}
Fix \(0 < \varepsilon \le \varepsilon_0\).
Following the finite-dimensional Galerkin strategy around a background profile (cf.\ \cite[Section~5.2]{BradshawTsai2017} in the divergence-free setting), choose a countable sequence \(\{\mathbf e_j\}_{j \ge 1} \subset C_c^\infty(\mathbb R^n; \mathbb R^n)\) containing, for each integer \(\ell \ge 1\), a \(C^1(\overline{B_\ell})\)-dense subset of \(C_c^\infty(B_\ell; \mathbb R^n)\), and define the finite-dimensional subspaces \(V_N := \operatorname{span}\{\mathbf e_1, \ldots, \mathbf e_N\}\) (\(N \ge 1\)) endowed with the \(L^2(\mathbb R^n)\) inner product.
Define the continuous map \(P_N : V_N \to V_N\) by pairing the left-hand side of \eqref{eq:rempenalized} (evaluated at \(\mathbf w \in V_N\)) minus \(\frac{1}{2} \mathbf S_R(1)\) against test functions in \(V_N\).

For every \(\mathbf w \in V_N\), antisymmetry of \(\Omega_{\mathbf w}\) and \(\Omega_{\mathbf B}\) gives the pointwise orthogonalities
\[
  (\Omega_{\mathbf w} \mathbf w) \cdot \mathbf w = 0,
  \qquad
  (\Omega_{\mathbf B} \mathbf w) \cdot \mathbf w = 0.
\]
Using the product-rule identities
\[
  \partial_i (\boldsymbol y \cdot \mathbf w)
  = w_i + \sum_{j=1}^n y_j \partial_i w_j,
  \qquad
  \partial_i (\mathbf B \cdot \mathbf w)
  = \sum_{j=1}^n \bigl((\partial_i B_j) w_j + B_j \partial_i w_j\bigr),
\]
we expand the remaining two rotational pairings pointwise as
\[
  \begin{aligned}
  -\frac{1}{2} (\Omega_{\mathbf w} \boldsymbol y) \cdot \mathbf w
  &= -\frac{1}{2} \sum_{i=1}^n \sum_{j=1}^n (\partial_j w_i - \partial_i w_j) y_j w_i \\
  &= -\frac{1}{4} \boldsymbol y \cdot \nabla |\mathbf w|^2
     + \frac{1}{2} \mathbf w \cdot \nabla(\boldsymbol y \cdot \mathbf w)
     - \frac{1}{2} |\mathbf w|^2, \\
  (\Omega_{\mathbf w} \mathbf B) \cdot \mathbf w
  &= \sum_{i=1}^n \sum_{j=1}^n (\partial_j w_i - \partial_i w_j) B_j w_i \\
  &= \frac{1}{2} \mathbf B \cdot \nabla |\mathbf w|^2
     - \mathbf w \cdot \nabla(\mathbf B \cdot \mathbf w)
     + ((\mathbf w \cdot \nabla) \mathbf B) \cdot \mathbf w.
  \end{aligned}
\]
Using compact support of \(\mathbf w\), \(\operatorname{div} \boldsymbol y = n\), and \(\operatorname{div} \mathbf B = 0\), integration by parts over \(\mathbb R^n\) yields
\[
  \begin{aligned}
  -\frac{1}{2} \int_{\mathbb R^n} (\Omega_{\mathbf w} \boldsymbol y) \cdot \mathbf w \, \mathrm d \boldsymbol y
  &= \frac{n - 2}{4} \|\mathbf w\|_2^2
     - \frac{1}{2} \int_{\mathbb R^n} (\boldsymbol y \cdot \mathbf w) \operatorname{div} \mathbf w \, \mathrm d \boldsymbol y, \\
  \int_{\mathbb R^n} (\Omega_{\mathbf w} \mathbf B) \cdot \mathbf w \, \mathrm d \boldsymbol y
  &= \int_{\mathbb R^n} ((\mathbf w \cdot \nabla) \mathbf B) \cdot \mathbf w \, \mathrm d \boldsymbol y
     + \int_{\mathbb R^n} (\mathbf B \cdot \mathbf w) \operatorname{div} \mathbf w \, \mathrm d \boldsymbol y.
  \end{aligned}
\]
Integrating the viscous and penalty terms by parts and completing the square on half of the weighted divergence norm gives
\begin{equation}\label{eq:Qcoercive}
  \begin{aligned}
  & \left \langle P_N(\mathbf w), \mathbf w \right \rangle _{L^2}
    + \frac{1}{2} \int_{\mathbb R^n} \mathbf S_R(1) \cdot \mathbf w \, \mathrm d \boldsymbol y \\
  ={}& \|\nabla \mathbf w\|_2^2 + \frac{n - 2}{4} \|\mathbf w\|_2^2
     + \int_{\mathbb R^n} ((\mathbf w \cdot \nabla) \mathbf B) \cdot \mathbf w \, \mathrm d \boldsymbol y
     + \frac{1}{\varepsilon^2} \|\langle \boldsymbol y \rangle \operatorname{div} \mathbf w\|_2^2 \\
  &+ \int_{\mathbb R^n}
         \left(\left(\mathbf B - \frac{1}{2} \boldsymbol y\right) \cdot \mathbf w\right)
         \operatorname{div} \mathbf w \, \mathrm d \boldsymbol y \\
={}& \|\nabla \mathbf w\|_2^2
     + \frac{1}{2\varepsilon^2} \|\langle \boldsymbol y \rangle \operatorname{div} \mathbf w\|_2^2
     + \frac{1}{2\varepsilon^2}
       \left\|\langle \boldsymbol y \rangle \operatorname{div} \mathbf w
             + \frac{\varepsilon^2}{\langle \boldsymbol y \rangle}
               \left(\mathbf B - \frac{1}{2} \boldsymbol y\right) \cdot \mathbf w\right\|_2^2 \\
  &+ \int_{\mathbb R^n}
       \left[\frac{n - 2}{4} |\mathbf w|^2
         + ((\mathbf w \cdot \nabla) \mathbf B) \cdot \mathbf w
         - \frac{\varepsilon^2}{2 \langle \boldsymbol y \rangle^2}
           \left|\left(\mathbf B - \frac{1}{2} \boldsymbol y\right) \cdot \mathbf w\right|^2\right] \mathrm d \boldsymbol y.
  \end{aligned}
\end{equation}

Since
\[
  \left|\mathbf B(\boldsymbol y) - \frac{1}{2} \boldsymbol y\right|
  \le (1 + \|\mathbf B\|_\infty) \langle \boldsymbol y \rangle
  \qquad\text{and}\qquad
  0 < \varepsilon \le \varepsilon_0 = \frac{1}{2(1 + \|\mathbf B\|_\infty)},
\]
we have
\[
  \frac{\varepsilon^2}{2 \langle \boldsymbol y \rangle^2}
  \left|\left(\mathbf B - \frac{1}{2} \boldsymbol y\right) \cdot \mathbf w\right|^2
  \le \frac{1}{8} |\mathbf w|^2.
\]
Together with \eqref{eq:coercivity} and \(n \ge 5\), the bracketed quadratic integrand in \eqref{eq:Qcoercive} is bounded below pointwise by \(\frac{n - 3}{8} |\mathbf w|^2 \ge \frac{1}{4} |\mathbf w|^2\), whence
\[
  \langle P_N(\mathbf w), \mathbf w \rangle_{L^2}
  \ge \|\nabla \mathbf w\|_2^2 + \frac{1}{4} \|\mathbf w\|_2^2
      + \frac{1}{2\varepsilon^2} \|\langle \boldsymbol y \rangle \operatorname{div} \mathbf w\|_2^2
      - \frac{1}{2} \|\mathbf S_R(1)\|_2 \|\mathbf w\|_2.
\]
Since \(\mathbf S_R(1) \in L^2(\mathbb R^n)\) by Proposition~\ref{prop:defect}, we have \(\langle P_N(\mathbf w), \mathbf w \rangle_{L^2} > 0\) on the sphere \(\|\mathbf w\|_2 = 2 \|\mathbf S_R(1)\|_2 + 1\) in \(V_N\).
By a standard corollary of Brouwer's fixed-point theorem \cite[Chapter~II, Lemma~1.4]{Temam2001}, there exists \(\mathbf W_N \in V_N\) solving \(P_N(\mathbf W_N) = \boldsymbol 0\), which satisfies
\[
  \|\nabla \mathbf W_N\|_2^2 + \frac{1}{8} \|\mathbf W_N\|_2^2
  + \frac{1}{2\varepsilon^2} \|\langle \boldsymbol y \rangle \operatorname{div} \mathbf W_N\|_2^2
  \le \frac{1}{2} \|\mathbf S_R(1)\|_2^2
\]
uniformly in \(N\) and \(\varepsilon \in (0, \varepsilon_0]\).

Passing to a subsequence as \(N \to \infty\) at fixed \(\varepsilon\), we obtain a limit \(\mathbf W_\varepsilon \in H^1(\mathbb R^n; \mathbb R^n)\) such that
\[
  \begin{aligned}
  \mathbf W_N &\rightharpoonup \mathbf W_\varepsilon \quad \text{weakly in } H^1(\mathbb R^n; \mathbb R^n), \\
  \mathbf W_N &\longrightarrow \mathbf W_\varepsilon \quad \text{strongly in } L^2(B_\ell; \mathbb R^n) \ \text{for every } \ell > 0, \\
  \langle \boldsymbol y \rangle \operatorname{div} \mathbf W_N &\rightharpoonup \langle \boldsymbol y \rangle \operatorname{div} \mathbf W_\varepsilon \quad \text{weakly in } L^2(\mathbb R^n).
  \end{aligned}
\]
Consequently, \(\Omega_{\mathbf W_N} \mathbf W_N \rightharpoonup \Omega_{\mathbf W_\varepsilon} \mathbf W_\varepsilon\) in \(L^1_{\mathrm{loc}}(\mathbb R^n)\) and \(\langle \boldsymbol y \rangle^2 \operatorname{div} \mathbf W_N \rightharpoonup \langle \boldsymbol y \rangle^2 \operatorname{div} \mathbf W_\varepsilon\) in \(L^2_{\mathrm{loc}}(\mathbb R^n)\).
By the \(C^1(\overline{B_\ell})\)-density of compactly supported elements of \(\bigcup_{N \ge 1} V_N\) in \(C_c^\infty(B_\ell; \mathbb R^n)\), the limit \(\mathbf W_\varepsilon\) solves \eqref{eq:rempenalized} (and hence \eqref{eq:penalizedlamb}) in \(\mathscr{D}'(\mathbb R^n)\).

On the completed-square right-hand side of \eqref{eq:Qcoercive}, the first three terms are squared \(L^2(\mathbb R^n)\) norms of continuous linear images of \((\mathbf w, \nabla \mathbf w, \langle \boldsymbol y \rangle \operatorname{div} \mathbf w) \in L^2\), while the fourth integrand is a pointwise positive-definite quadratic form bounded below by \(\frac{1}{4} |\mathbf w|^2\).
All four terms are therefore convex and weakly lower semicontinuous.
Evaluating \eqref{eq:Qcoercive} at \(\mathbf w = \mathbf W_N\) (where \(P_N(\mathbf W_N) = \boldsymbol 0\)), taking \(\liminf_{N \to \infty}\), and noting that
\[
  \int_{\mathbb R^n} \mathbf S_R(1) \cdot \mathbf W_N \, \mathrm d \boldsymbol y
  \longrightarrow \int_{\mathbb R^n} \mathbf S_R(1) \cdot \mathbf W_\varepsilon \, \mathrm d \boldsymbol y,
\]
we expand the completed square back out at \(\mathbf W_\varepsilon\) to obtain \eqref{eq:approxenergy} and \eqref{eq:approxbounds}.
\end{proof}

We next establish an \(\varepsilon\)-independent upper bound on the penalized Bernoulli function \(\Phi_\varepsilon\) for the approximate solutions constructed in Proposition~\ref{prop:approx}, preparing for the incompressible limit \(\varepsilon \downarrow 0\).

\begin{proposition}\label{prop:bernoullimax}
Let \(n \ge 5\), and let \(\mathbf B \in C^\infty(\mathbb R^n; \mathbb R^n) \cap L^\infty(\mathbb R^n)\) and \(\Phi_B \in C^\infty(\mathbb R^n) \cap L^\infty(\mathbb R^n)\) be as in Proposition~\ref{prop:defect}.
Suppose \(0 < \varepsilon \le 1\) and \(\mathbf W_\varepsilon \in H^1(\mathbb R^n; \mathbb R^n)\) with \(\langle \boldsymbol y \rangle \operatorname{div} \mathbf W_\varepsilon \in L^2(\mathbb R^n)\) solves \eqref{eq:penalizedlamb} in \(\mathscr{D}'(\mathbb R^n)\), where \(\mathbf U_\varepsilon = \mathbf B + \mathbf W_\varepsilon\) and \(\Phi_\varepsilon\) is given by \eqref{eq:penaltydef}.
Then
\begin{equation}\label{eq:epsmax}
  \Phi_\varepsilon(\boldsymbol y)
  \le \sup_{\mathbb R^n} \Phi_B
  \quad \text{for a.e. } \boldsymbol y \in \mathbb R^n.
\end{equation}
\end{proposition}

\begin{proof}
\noindent\emph{Step 1.}
Set \(\mathbf v_\varepsilon := \mathbf U_\varepsilon - \frac{1}{2} \boldsymbol y \in H^1_{\mathrm{loc}}(\mathbb R^n; \mathbb R^n)\).
Since \(\Omega_{\boldsymbol y} = \boldsymbol 0\), we have \(\Omega_{\mathbf v_\varepsilon} = \Omega_{\mathbf U_\varepsilon}\) and \(-\frac{1}{2} \Omega_{\mathbf U_\varepsilon} \boldsymbol y + \Omega_{\mathbf U_\varepsilon} \mathbf U_\varepsilon = \Omega_{\mathbf v_\varepsilon} \mathbf v_\varepsilon\), so \eqref{eq:penalizedlamb} takes the form in \(\mathscr{D}'(\mathbb R^n)\):
\begin{equation}\label{eq:driftlamb}
  -\operatorname{div} \Omega_{\mathbf v_\varepsilon}
  + \Omega_{\mathbf v_\varepsilon} \mathbf v_\varepsilon
  + \nabla \Phi_\varepsilon = \boldsymbol 0.
\end{equation}

Since \(\mathbf v_\varepsilon \in H^1_{\mathrm{loc}}(\mathbb R^n; \mathbb R^n)\) gives \(\Omega_{\mathbf v_\varepsilon} \mathbf v_\varepsilon \in L^{n/(n-1)}_{\mathrm{loc}}(\mathbb R^n) \subset L^1_{\mathrm{loc}}(\mathbb R^n)\) and does not in general justify testing \eqref{eq:driftlamb} directly against \(\mathbf v_\varepsilon \varphi\) when \(n \ge 5\), we first apply a collinear velocity truncation.
Fix a standard nonnegative radial mollifier \(\zeta_1 \in C_c^\infty(B_1)\) with unit integral, set \(\zeta_\delta(\boldsymbol z) := \delta^{-n} \zeta_1(\boldsymbol z / \delta)\), and for any \(\varphi \in C_c^\infty(B_\ell)\) and \(K > 0\) define
\[
  \mathbf v_{\varepsilon,K}
  := \min\left\{1, \, \frac{K}{|\mathbf v_\varepsilon|}\right\} \mathbf v_\varepsilon.
\]
Since \(\boldsymbol z \mapsto \min\{1, K / |\boldsymbol z|\} \boldsymbol z\) is the metric projection onto the closed ball \(\overline{B_K} \subset \mathbb R^n\), it is globally \(1\)-Lipschitz.
Hence \(\mathbf v_{\varepsilon,K} \in H^1_{\mathrm{loc}}(\mathbb R^n; \mathbb R^n) \cap L^\infty_{\mathrm{loc}}(\mathbb R^n; \mathbb R^n)\) satisfies almost everywhere
\[
  |\mathbf v_{\varepsilon,K}| \le |\mathbf v_\varepsilon|,
  \qquad
  |\nabla \mathbf v_{\varepsilon,K}| \le |\nabla \mathbf v_\varepsilon|.
\]
Because \(\mathbf v_{\varepsilon,K}\) is collinear with \(\mathbf v_\varepsilon\), the pointwise orthogonality \((\Omega_{\mathbf v_\varepsilon} \mathbf v_\varepsilon) \cdot \mathbf v_{\varepsilon,K} = 0\) remains valid.
Testing \eqref{eq:driftlamb} with \((\mathbf v_{\varepsilon,K} \varphi) * \zeta_\delta\) and sending \(\delta \downarrow 0\) at fixed \(K > 0\), we obtain by dominated convergence
\[
  \lim_{\delta \downarrow 0} \int_{\mathbb R^n}
    (\Omega_{\mathbf v_\varepsilon} \mathbf v_\varepsilon) \cdot
    \bigl[(\mathbf v_{\varepsilon,K} \varphi) * \zeta_\delta\bigr] \, \mathrm d \boldsymbol y
  = \int_{\mathbb R^n} (\Omega_{\mathbf v_\varepsilon} \mathbf v_\varepsilon)
    \cdot (\mathbf v_{\varepsilon,K} \varphi) \, \mathrm d \boldsymbol y
  = 0.
\]
In the remaining linear pairings with \(\Omega_{\mathbf v_\varepsilon}, \Phi_\varepsilon \in L^2(B_\ell)\), note that \(\mathbf v_{\varepsilon,K} - \mathbf v_\varepsilon = \boldsymbol 0\) on \(\{|\mathbf v_\varepsilon| \le K\}\).
Applying Stampacchia's lemma on level sets componentwise gives \(\nabla \mathbf v_{\varepsilon,K} = \nabla \mathbf v_\varepsilon\) almost everywhere on \(\{|\mathbf v_\varepsilon| \le K\}\).
Combined with the dominating bounds above and \(\bigl|\{|\mathbf v_\varepsilon| > K\} \cap B_\ell\bigr| \to 0\), dominated convergence yields \(\mathbf v_{\varepsilon,K} \varphi \to \mathbf v_\varepsilon \varphi\) strongly in \(H^1_0(B_\ell; \mathbb R^n)\) as \(K \to \infty\), so
\[
  \int_{\mathbb R^n} \Omega_{\mathbf v_\varepsilon} : \nabla(\mathbf v_\varepsilon \varphi) \, \mathrm d \boldsymbol y
  - \int_{\mathbb R^n} \Phi_\varepsilon \operatorname{div}(\mathbf v_\varepsilon \varphi) \, \mathrm d \boldsymbol y
  = 0.
\]
Using antisymmetry of \(\Omega_{\mathbf v_\varepsilon}\), the identity \(|\Omega_{\mathbf v_\varepsilon}|^2 = |\Omega_{\mathbf U_\varepsilon}|^2\), and relabeling \(i \leftrightarrow j\) in the second sum below, we have almost everywhere on \(B_\ell\)
\[
  \begin{aligned}
  \Omega_{\mathbf v_\varepsilon} : \nabla(\mathbf v_\varepsilon \varphi)
  &= \sum_{i=1}^n \sum_{j=1}^n (\partial_j v_{\varepsilon,i} - \partial_i v_{\varepsilon,j}) \partial_j v_{\varepsilon,i} \, \varphi
     + \sum_{i=1}^n \sum_{j=1}^n (\partial_j v_{\varepsilon,i} - \partial_i v_{\varepsilon,j}) v_{\varepsilon,i} \partial_j \varphi \\
  &= \frac{1}{2} |\Omega_{\mathbf U_\varepsilon}|^2 \varphi
     - \sum_{i=1}^n \sum_{j=1}^n (\partial_j v_{\varepsilon,i} - \partial_i v_{\varepsilon,j}) v_{\varepsilon,j} \partial_i \varphi \\
  &= \frac{1}{2} |\Omega_{\mathbf U_\varepsilon}|^2 \varphi
     - (\Omega_{\mathbf v_\varepsilon} \mathbf v_\varepsilon) \cdot \nabla \varphi,
  \end{aligned}
\]
whereas testing \eqref{eq:driftlamb} against \(\nabla \varphi\) and using \(\operatorname{div}(\operatorname{div} \Omega_{\mathbf v_\varepsilon}) = 0\) gives
\[
  \int_{\mathbb R^n} (\Omega_{\mathbf v_\varepsilon} \mathbf v_\varepsilon) \cdot \nabla \varphi \, \mathrm d \boldsymbol y
  = \int_{\mathbb R^n} \Phi_\varepsilon \Delta \varphi \, \mathrm d \boldsymbol y.
\]
Substituting both identities above together with
\[
  \operatorname{div}(\mathbf v_\varepsilon \varphi)
  = \mathbf v_\varepsilon \cdot \nabla \varphi + \left(\operatorname{div} \mathbf W_\varepsilon - \frac{n}{2}\right) \varphi,
\]
we obtain in \(\mathscr{D}'(\mathbb R^n)\):
\begin{equation}\label{eq:weakbernoulli}
  -\Delta \Phi_\varepsilon
  + \operatorname{div}(\mathbf v_\varepsilon \Phi_\varepsilon)
  - \left(\operatorname{div} \mathbf W_\varepsilon - \frac{n}{2}\right) \Phi_\varepsilon
  = -\frac{1}{2} |\Omega_{\mathbf U_\varepsilon}|^2.
\end{equation}

\smallskip
\noindent\emph{Step 2.}
Since \(\Phi_\varepsilon \in L^2_{\mathrm{loc}}(\mathbb R^n)\) and \(\mathbf v_\varepsilon \in H^1_{\mathrm{loc}}(\mathbb R^n; \mathbb R^n)\), we can renormalize \eqref{eq:weakbernoulli} via Friedrichs' commutator lemma \cite[Lemma~II.1]{DL89}.
For every bounded open set \(D \Subset D' \Subset \mathbb R^n\), the commutator
\[
  c_\delta := \left[\operatorname{div}(\mathbf v_\varepsilon \Phi_\varepsilon)
    - \left(\operatorname{div} \mathbf W_\varepsilon - \frac{n}{2}\right) \Phi_\varepsilon\right] * \zeta_\delta
    - \mathbf v_\varepsilon \cdot \nabla (\Phi_\varepsilon * \zeta_\delta)
\]
satisfies
\[
  \|c_\delta\|_{L^1(D)} \le C \|\mathbf v_\varepsilon\|_{H^1(D')} \|\Phi_\varepsilon\|_{L^2(D')},
  \qquad
  \lim_{\delta \downarrow 0} \|c_\delta\|_{L^1(D)} = 0.
\]
Convolving \eqref{eq:weakbernoulli} with \(\zeta_\delta \ge 0\) gives
\[
  -\Delta (\Phi_\varepsilon * \zeta_\delta)
  + \mathbf v_\varepsilon \cdot \nabla (\Phi_\varepsilon * \zeta_\delta)
  \le |c_\delta|.
\]
Consider any smooth convex nondecreasing function \(\theta\) vanishing on \((-\infty, \sup_{\mathbb R^n} \Phi_B]\) with \(0 \le \theta' \le 1\) and \(\theta'' \ge 0\); for example, for \(\kappa > 0\),
\[
  \theta_\kappa(s)
  := \mathbf 1_{\{s > \sup_{\mathbb R^n} \Phi_B\}}
     \bigl(s - \sup_{\mathbb R^n} \Phi_B\bigr)
     \exp\!\left(-\frac{\kappa}{s - \sup_{\mathbb R^n} \Phi_B}\right).
\]
Then we have pointwise
\[
  \begin{aligned}
  -\theta'(\Phi_\varepsilon * \zeta_\delta) \Delta (\Phi_\varepsilon * \zeta_\delta)
  &= -\Delta \theta(\Phi_\varepsilon * \zeta_\delta)
     + \theta''(\Phi_\varepsilon * \zeta_\delta) |\nabla (\Phi_\varepsilon * \zeta_\delta)|^2
   \ge -\Delta \theta(\Phi_\varepsilon * \zeta_\delta), \\
  \theta'(\Phi_\varepsilon * \zeta_\delta) \mathbf v_\varepsilon \cdot \nabla (\Phi_\varepsilon * \zeta_\delta)
  &= \mathbf v_\varepsilon \cdot \nabla \theta(\Phi_\varepsilon * \zeta_\delta) \\
  &= \operatorname{div}\bigl(\mathbf v_\varepsilon \, \theta(\Phi_\varepsilon * \zeta_\delta)\bigr)
     + \left(\frac{n}{2} - \operatorname{div} \mathbf W_\varepsilon\right) \theta(\Phi_\varepsilon * \zeta_\delta),
  \end{aligned}
\]
so multiplying the mollified inequality by \(\theta'(\Phi_\varepsilon * \zeta_\delta)\) yields
\[
  -\Delta \theta(\Phi_\varepsilon * \zeta_\delta)
  + \operatorname{div}\bigl(\mathbf v_\varepsilon \, \theta(\Phi_\varepsilon * \zeta_\delta)\bigr)
  + \left(\frac{n}{2} - \operatorname{div} \mathbf W_\varepsilon\right) \theta(\Phi_\varepsilon * \zeta_\delta)
  \le |c_\delta|.
\]
Since \(0 \le \theta' \le 1\) implies \(\theta(\Phi_\varepsilon * \zeta_\delta) \to \theta(\Phi_\varepsilon)\) strongly in \(L^2_{\mathrm{loc}}(\mathbb R^n)\) while \(\mathbf v_\varepsilon, \operatorname{div} \mathbf W_\varepsilon \in L^2_{\mathrm{loc}}(\mathbb R^n)\) and \(\|c_\delta\|_{L^1(D)} \to 0\), sending \(\delta \downarrow 0\) and then taking \(\theta = \theta_\kappa\) with \(\kappa \downarrow 0\) yields in \(\mathscr{D}'(\mathbb R^n)\):
\begin{equation}\label{eq:posbernoulli}
  \begin{aligned}
  -\Delta \bigl(\Phi_\varepsilon - \sup_{\mathbb R^n} \Phi_B\bigr)_+
   + \operatorname{div}\left[\left(\mathbf U_\varepsilon - \frac{1}{2} \boldsymbol y\right)
     \bigl(\Phi_\varepsilon - \sup_{\mathbb R^n} \Phi_B\bigr)_+\right]& \\
  + \left(\frac{n}{2} - \operatorname{div} \mathbf W_\varepsilon\right)
     \bigl(\Phi_\varepsilon - \sup_{\mathbb R^n} \Phi_B\bigr)_+ &\le 0.
  \end{aligned}
\end{equation}

\smallskip
\noindent\emph{Step 3.}
By \eqref{eq:penaltydef}, \(\Phi_B \le \sup_{\mathbb R^n} \Phi_B\), and \(1 + \varepsilon^{-2} \langle \boldsymbol y \rangle^2 \le 2 \varepsilon^{-2} \langle \boldsymbol y \rangle^2\) for \(0 < \varepsilon \le 1\), we have almost everywhere on \(\mathbb R^n\):
\begin{align}
  -(\operatorname{div} \mathbf W_\varepsilon)
  \bigl(\Phi_\varepsilon - \sup_{\mathbb R^n} \Phi_B\bigr)_+
  &= \frac{\Phi_\varepsilon - \Phi_B}{1 + \frac{\langle \boldsymbol y \rangle^2}{\varepsilon^2}}
     \bigl(\Phi_\varepsilon - \sup_{\mathbb R^n} \Phi_B\bigr)_+ \nonumber \\
  &\ge \frac{\varepsilon^2 \bigl(\Phi_\varepsilon - \sup_{\mathbb R^n} \Phi_B\bigr)_+^2}
            {2 \langle \boldsymbol y \rangle^2}. \label{eq:divsign}
\end{align}
We test \eqref{eq:posbernoulli} against a radial barrier \(\langle \alpha \boldsymbol y \rangle^{-\gamma}\) with \(\gamma > 0\) and \(0 < \alpha \le 1\).
Integrating the linear drift by parts applies \(\frac{1}{2} \boldsymbol y \cdot \nabla + \frac{n}{2}\) to \(\langle \alpha \boldsymbol y \rangle^{-\gamma}\), which is bounded below by a positive multiple of \(\langle \alpha \boldsymbol y \rangle^{-\gamma}\) only when \(\gamma < n\).
On the other hand, \eqref{eq:penaltydef} and \(\langle \boldsymbol y \rangle \operatorname{div} \mathbf W_\varepsilon \in L^2(\mathbb R^n)\) give
\[
  \frac{\bigl(\Phi_\varepsilon - \sup_{\mathbb R^n} \Phi_B\bigr)_+}{\langle \boldsymbol y \rangle}
  \le \frac{2}{\varepsilon^2} \langle \boldsymbol y \rangle |\operatorname{div} \mathbf W_\varepsilon|
  \in L^2(\mathbb R^n),
\]
so ensuring integrability of \(\bigl(\Phi_\varepsilon - \sup_{\mathbb R^n} \Phi_B\bigr)_+ \langle \alpha \boldsymbol y \rangle^{-\gamma}\) by Cauchy--Schwarz requires \(\langle \boldsymbol y \rangle \langle \alpha \boldsymbol y \rangle^{-\gamma} \in L^2(\mathbb R^n)\), that is, \(\gamma > (n + 2)/2\).
Since \(n \ge 5\) implies \((n + 2)/2 < n - 1 < n\), we choose \(\gamma = n - 1\) and work with \(\langle \alpha \boldsymbol y \rangle^{1-n}\).

Direct differentiation gives
\[
  \begin{aligned}
  \nabla_{\boldsymbol y} \left(\frac{1}{\langle \alpha \boldsymbol y \rangle^{n-1}}\right)
  &= -\frac{(n - 1) \alpha^2 \boldsymbol y}{1 + \alpha^2 |\boldsymbol y|^2}
     \frac{1}{\langle \alpha \boldsymbol y \rangle^{n-1}}, \\
  -\Delta_{\boldsymbol y} \left(\frac{1}{\langle \alpha \boldsymbol y \rangle^{n-1}}\right)
  &= (n - 1) \alpha^2 \frac{n - \alpha^2 |\boldsymbol y|^2}{(1 + \alpha^2 |\boldsymbol y|^2)^2}
     \frac{1}{\langle \alpha \boldsymbol y \rangle^{n-1}}
   \ge -\frac{(n - 1) \alpha^2}{\langle \alpha \boldsymbol y \rangle^{n-1}}, \\
  \left(\frac{1}{2} \boldsymbol y \cdot \nabla_{\boldsymbol y} + \frac{n}{2}\right)
    \left(\frac{1}{\langle \alpha \boldsymbol y \rangle^{n-1}}\right)
  &= \left(\frac{n}{2} - \frac{n - 1}{2} \frac{\alpha^2 |\boldsymbol y|^2}{1 + \alpha^2 |\boldsymbol y|^2}\right)
     \frac{1}{\langle \alpha \boldsymbol y \rangle^{n-1}}
   \ge \frac{1}{2 \langle \alpha \boldsymbol y \rangle^{n-1}}, \\
  \left|\mathbf B \cdot \nabla_{\boldsymbol y} \left(\frac{1}{\langle \alpha \boldsymbol y \rangle^{n-1}}\right)\right|
  &\le \frac{(n - 1) \alpha \|\mathbf B\|_\infty}{\langle \alpha \boldsymbol y \rangle^{n-1}},
  \end{aligned}
\]
whence for all sufficiently small \(\alpha \in (0, 1]\),
\[
  \left(-\Delta + \left(\frac{1}{2} \boldsymbol y - \mathbf B\right) \cdot \nabla + \frac{n}{2}\right)
  \left(\frac{1}{\langle \alpha \boldsymbol y \rangle^{n-1}}\right)
  \ge \frac{1}{4 \langle \alpha \boldsymbol y \rangle^{n-1}}
  \quad \text{on } \mathbb R^n.
\]

Observe that \(\langle \boldsymbol y \rangle |\nabla \langle \alpha \boldsymbol y \rangle^{1-n}| \le (n - 1) \langle \alpha \boldsymbol y \rangle^{1-n}\) and \(\langle \boldsymbol y \rangle \langle \alpha \boldsymbol y \rangle^{1-n} \in L^2(\mathbb R^n)\).
By Cauchy--Schwarz, the nonnegative function \(\bigl(\Phi_\varepsilon - \sup_{\mathbb R^n} \Phi_B\bigr)_+\) is therefore integrable on \(\mathbb R^n\) against \(\langle \alpha \boldsymbol y \rangle^{1-n}\), against \(|\operatorname{div} \mathbf W_\varepsilon| \langle \alpha \boldsymbol y \rangle^{1-n}\), and against \(|\mathbf W_\varepsilon| |\nabla \langle \alpha \boldsymbol y \rangle^{1-n}|\).
Fix \(\eta \in C_c^\infty(B_2)\) with \(0 \le \eta \le 1\) and \(\eta|_{B_1} = 1\), split \(\mathbf U_\varepsilon - \frac{1}{2} \boldsymbol y = (\mathbf B - \frac{1}{2} \boldsymbol y) + \mathbf W_\varepsilon\) in \eqref{eq:posbernoulli}, and test against \(\langle \alpha \boldsymbol y \rangle^{1-n} \eta(\boldsymbol y / \ell)\) (\(\ell \ge 1\)).
On the annulus \(\{\ell \le |\boldsymbol y| \le 2\ell\}\), we have
\[
  |\nabla \eta(\boldsymbol y / \ell)| \le \frac{C}{\ell} \le \frac{3 C}{\langle \boldsymbol y \rangle},
  \qquad
  |\boldsymbol y| |\nabla \eta(\boldsymbol y / \ell)| \le C,
  \qquad
  |\Delta \eta(\boldsymbol y / \ell)| \le \frac{C}{\ell^2},
\]
so all terms in which a derivative hits \(\eta(\boldsymbol y / \ell)\) are bounded by
\[
  \begin{aligned}
  &C \int_{\ell \le |\boldsymbol y| \le 2\ell}
     \frac{\bigl(\Phi_\varepsilon - \sup_{\mathbb R^n} \Phi_B\bigr)_+}
          {\langle \alpha \boldsymbol y \rangle^{n-1}} \, \mathrm d \boldsymbol y \\
  {}+{}& C \left(\int_{\ell \le |\boldsymbol y| \le 2\ell}
             \frac{\bigl(\Phi_\varepsilon - \sup_{\mathbb R^n} \Phi_B\bigr)_+^2}
                  {\langle \boldsymbol y \rangle^2} \, \mathrm d \boldsymbol y\right)^{1/2}
           \left(\int_{\ell \le |\boldsymbol y| \le 2\ell} |\mathbf W_\varepsilon|^2 \, \mathrm d \boldsymbol y\right)^{1/2}
  \longrightarrow 0
  \end{aligned}
\]
as \(\ell \to \infty\).
Sending \(\ell \to \infty\) and substituting the linear barrier lower bound above together with \eqref{eq:divsign} and the explicit formula for \(\nabla \langle \alpha \boldsymbol y \rangle^{1-n}\) yields
\[
  \begin{aligned}
  &\frac{1}{4} \int_{\mathbb R^n}
     \frac{\bigl(\Phi_\varepsilon - \sup_{\mathbb R^n} \Phi_B\bigr)_+}
          {\langle \alpha \boldsymbol y \rangle^{n-1}} \, \mathrm d \boldsymbol y
   + \frac{\varepsilon^2}{2} \int_{\mathbb R^n}
     \frac{\bigl(\Phi_\varepsilon - \sup_{\mathbb R^n} \Phi_B\bigr)_+^2}
          {\langle \boldsymbol y \rangle^2 \langle \alpha \boldsymbol y \rangle^{n-1}} \, \mathrm d \boldsymbol y \\
  \le{}& \int_{\mathbb R^n}
       \bigl(\Phi_\varepsilon - \sup_{\mathbb R^n} \Phi_B\bigr)_+ \,
       \mathbf W_\varepsilon \cdot \nabla_{\boldsymbol y}
         \left(\frac{1}{\langle \alpha \boldsymbol y \rangle^{n-1}}\right) \mathrm d \boldsymbol y \\
  \le{}& (n - 1) \int_{\mathbb R^n}
       \frac{\alpha^2 |\boldsymbol y|}{1 + \alpha^2 |\boldsymbol y|^2} |\mathbf W_\varepsilon|
       \frac{\bigl(\Phi_\varepsilon - \sup_{\mathbb R^n} \Phi_B\bigr)_+}
            {\langle \alpha \boldsymbol y \rangle^{n-1}} \, \mathrm d \boldsymbol y.
  \end{aligned}
\]

To absorb the right-hand side into half of the quadratic superlevel integral at fixed \(\varepsilon > 0\) before sending \(\alpha \downarrow 0\), we factor the integrand via \(\langle \boldsymbol y \rangle^{-1} \cdot \langle \boldsymbol y \rangle\) and apply Young's inequality:
\[
  \begin{aligned}
  &(n - 1) \int_{\mathbb R^n}
     \frac{\alpha^2 |\boldsymbol y|}{1 + \alpha^2 |\boldsymbol y|^2} |\mathbf W_\varepsilon|
     \frac{\bigl(\Phi_\varepsilon - \sup_{\mathbb R^n} \Phi_B\bigr)_+}
          {\langle \alpha \boldsymbol y \rangle^{n-1}} \, \mathrm d \boldsymbol y \\
  ={}& \int_{\mathbb R^n}
     \left[\frac{\varepsilon}{\sqrt 2}
           \frac{\bigl(\Phi_\varepsilon - \sup_{\mathbb R^n} \Phi_B\bigr)_+}{\langle \boldsymbol y \rangle}\right]
     \left[\frac{\sqrt 2 (n - 1)}{\varepsilon}
           \frac{\alpha^2 |\boldsymbol y| \langle \boldsymbol y \rangle}{1 + \alpha^2 |\boldsymbol y|^2}
           |\mathbf W_\varepsilon|\right]
     \frac{\mathrm d \boldsymbol y}{\langle \alpha \boldsymbol y \rangle^{n-1}} \\
  \le{}& \frac{\varepsilon^2}{4} \int_{\mathbb R^n}
       \frac{\bigl(\Phi_\varepsilon - \sup_{\mathbb R^n} \Phi_B\bigr)_+^2}
            {\langle \boldsymbol y \rangle^2 \langle \alpha \boldsymbol y \rangle^{n-1}} \, \mathrm d \boldsymbol y
     + \frac{(n - 1)^2}{\varepsilon^2} \int_{\mathbb R^n}
       \frac{\alpha^4 |\boldsymbol y|^2 \langle \boldsymbol y \rangle^2}{(1 + \alpha^2 |\boldsymbol y|^2)^2}
       \frac{|\mathbf W_\varepsilon|^2}{\langle \alpha \boldsymbol y \rangle^{n-1}} \, \mathrm d \boldsymbol y.
  \end{aligned}
\]
Because \(\langle \boldsymbol y \rangle^{-1} (\Phi_\varepsilon - \sup_{\mathbb R^n} \Phi_B)_+ \in L^2(\mathbb R^n)\) keeps the quadratic superlevel integral finite at fixed \(\varepsilon > 0\), absorbing that term into the left-hand side gives
\[
  \begin{aligned}
  &\frac{1}{4} \int_{\mathbb R^n}
     \frac{\bigl(\Phi_\varepsilon - \sup_{\mathbb R^n} \Phi_B\bigr)_+}
          {\langle \alpha \boldsymbol y \rangle^{n-1}} \, \mathrm d \boldsymbol y
   + \frac{\varepsilon^2}{4} \int_{\mathbb R^n}
     \frac{\bigl(\Phi_\varepsilon - \sup_{\mathbb R^n} \Phi_B\bigr)_+^2}
          {\langle \boldsymbol y \rangle^2 \langle \alpha \boldsymbol y \rangle^{n-1}} \, \mathrm d \boldsymbol y \\
  \le{}& \frac{(n - 1)^2}{\varepsilon^2} \int_{\mathbb R^n}
       \frac{\alpha^4 |\boldsymbol y|^2 \langle \boldsymbol y \rangle^2}{(1 + \alpha^2 |\boldsymbol y|^2)^2}
       \frac{|\mathbf W_\varepsilon|^2}{\langle \alpha \boldsymbol y \rangle^{n-1}} \, \mathrm d \boldsymbol y.
  \end{aligned}
\]
For all \(0 < \alpha \le 1\) and \(\boldsymbol y \in \mathbb R^n\), we have \(\alpha^2 |\boldsymbol y|^2 \le 1 + \alpha^2 |\boldsymbol y|^2\) and \(\alpha^2 \langle \boldsymbol y \rangle^2 = \alpha^2 + \alpha^2 |\boldsymbol y|^2 \le 1 + \alpha^2 |\boldsymbol y|^2\), so the factor multiplying \(|\mathbf W_\varepsilon|^2\) is bounded uniformly by \(1\) and vanishes pointwise as \(\alpha \downarrow 0\).
Holding \(\varepsilon > 0\) fixed and sending \(\alpha \downarrow 0\), we see that the right-hand side vanishes by dominated convergence since \(\mathbf W_\varepsilon \in L^2(\mathbb R^n)\), whereas monotone convergence on the left-hand side yields \(\bigl(\Phi_\varepsilon - \sup_{\mathbb R^n} \Phi_B\bigr)_+ = 0\) almost everywhere on \(\mathbb R^n\), proving \eqref{eq:epsmax}.
\end{proof}

Propositions~\ref{prop:approx} and~\ref{prop:bernoullimax} provide approximate solutions satisfying \eqref{eq:approxbounds}--\eqref{eq:approxenergy} together with the uniform upper bound \eqref{eq:epsmax}.
Because \eqref{eq:penaltydef} and \eqref{eq:approxbounds} bound \(\Phi_\varepsilon\) in \(L^2_{\mathrm{loc}}(\mathbb R^n)\) only as \(O(\varepsilon^{-1})\), passing \eqref{eq:epsmax} to the incompressible limit \(\varepsilon \downarrow 0\) requires strong local convergence of \(\Pi_\varepsilon\) and \(\Phi_\varepsilon\), which we obtain by rewriting \eqref{eq:penalizedlamb} via \eqref{eq:lambidentity} in momentum form:
\begin{equation}\label{eq:approxmomentum}
  \begin{gathered}
  -\operatorname{div} \Omega_{\mathbf U_\varepsilon}
  - \frac{1}{2} \boldsymbol y \cdot \nabla \mathbf U_\varepsilon - \frac{1}{2} \mathbf U_\varepsilon
  + (\mathbf U_\varepsilon \cdot \nabla) \mathbf U_\varepsilon
  + \nabla \Pi_\varepsilon = \boldsymbol 0, \\
  \Pi_\varepsilon := \Phi_\varepsilon - \frac{1}{2} |\mathbf U_\varepsilon|^2
                     + \frac{1}{2} \boldsymbol y \cdot \mathbf U_\varepsilon.
  \end{gathered}
\end{equation}

\begin{proposition}\label{prop:limitprofile}
Let \(n \ge 5\), let \(\mathbf u_0\) satisfy the hypotheses of Theorem~\ref{thm:main}, and let \(R \ge 1\) with \(R^{-2} \le T_0\) be fixed so that \(\mathbf B = \mathbf u_R(1)\) satisfies \eqref{eq:coercivity}.
Then there exists \(\mathbf W \in H^1_\sigma(\mathbb R^n)\) such that \(\mathbf U = \mathbf B + \mathbf W\) and the canonical pressure \(\Pi\) defined in \eqref{eq:mainpressure} solve \eqref{eq:profile} in \(\mathscr{D}'(\mathbb R^n)\) and satisfy both the energy inequality \eqref{eq:mainenergy} and the modified Bernoulli bound \eqref{eq:mainmax} almost everywhere on \(\mathbb R^n\).
\end{proposition}

\begin{proof}
\noindent\emph{Step 1.}
For \(0 < \varepsilon \le \varepsilon_0\), let \(\mathbf W_\varepsilon \in H^1(\mathbb R^n; \mathbb R^n)\) be an approximate solution from Proposition~\ref{prop:approx}; by Proposition~\ref{prop:bernoullimax}, the associated function \(\Phi_\varepsilon\) also satisfies \eqref{eq:epsmax}.
We take the distributional divergence of \eqref{eq:approxmomentum} using \(\operatorname{div}(\operatorname{div} \Omega_{\mathbf U_\varepsilon}) = 0\) and \(\operatorname{div} \mathbf U_\varepsilon = \operatorname{div} \mathbf W_\varepsilon\).
To absorb one derivative from \(\boldsymbol y \cdot \nabla \operatorname{div} \mathbf W_\varepsilon\) into an outer divergence, we expand
\[
  \begin{aligned}
  \operatorname{div}(\boldsymbol y \cdot \nabla \mathbf U_\varepsilon)
  &= \sum_{i=1}^n \sum_{j=1}^n \partial_i (y_j \partial_j U_{\varepsilon,i})
   = \operatorname{div} \mathbf W_\varepsilon + \boldsymbol y \cdot \nabla \operatorname{div} \mathbf W_\varepsilon \\
  &= \operatorname{div}(\boldsymbol y \operatorname{div} \mathbf W_\varepsilon) - (n - 1) \operatorname{div} \mathbf W_\varepsilon,
  \end{aligned}
\]
which combines with \(\operatorname{div}(-\frac{1}{2} \mathbf U_\varepsilon) = -\frac{1}{2} \operatorname{div} \mathbf W_\varepsilon\) and the convective identity
\[
  \begin{aligned}
  \operatorname{div}\bigl((\mathbf U_\varepsilon \cdot \nabla) \mathbf U_\varepsilon\bigr)
  &= \sum_{i=1}^n \sum_{j=1}^n \partial_i \bigl(U_{\varepsilon,j} \partial_j U_{\varepsilon,i}\bigr) \\
  &= \sum_{i=1}^n \sum_{j=1}^n \partial_i \partial_j (U_{\varepsilon,i} U_{\varepsilon,j})
     - \operatorname{div}(\mathbf U_\varepsilon \operatorname{div} \mathbf W_\varepsilon)
  \end{aligned}
\]
to show that \(\Pi_\varepsilon\) satisfies in \(\mathscr{D}'(\mathbb R^n)\):
\begin{equation}\label{eq:divPeeps}
  -\Delta \Pi_\varepsilon
  = \sum_{i=1}^n \sum_{j=1}^n \partial_i \partial_j (U_{\varepsilon,i} U_{\varepsilon,j})
    - \operatorname{div}\left[\left(\frac{1}{2} \boldsymbol y + \mathbf U_\varepsilon\right)
                        \operatorname{div} \mathbf W_\varepsilon\right]
    + \frac{n - 2}{2} \operatorname{div} \mathbf W_\varepsilon.
\end{equation}

As the middle term in \eqref{eq:divPeeps} is in divergence form and \((-\Delta)^{-1} \partial_i = \mathcal R_i (-\Delta)^{-1/2}\), applying \((-\Delta)^{-1}\) to the last two terms of \eqref{eq:divPeeps} leads to the divergence-error potential
\[
  E_\varepsilon
  := -\sum_{i=1}^n \mathcal R_i (-\Delta)^{-1/2}
       \left[\left(\frac{1}{2} y_i + U_{\varepsilon,i}\right) \operatorname{div} \mathbf W_\varepsilon\right]
     + \frac{n - 2}{2} (-\Delta)^{-1} \operatorname{div} \mathbf W_\varepsilon.
\]
Splitting \(\frac{1}{2} \boldsymbol y + \mathbf U_\varepsilon = (\frac{1}{2} \boldsymbol y + \mathbf B) + \mathbf W_\varepsilon\), we obtain from \eqref{eq:approxbounds}, \(|\frac{1}{2} \boldsymbol y + \mathbf B| \le C \langle \boldsymbol y \rangle\), H\"older's inequality, and the Sobolev embedding \(H^1(\mathbb R^n) \hookrightarrow L^{2n/(n-2)}(\mathbb R^n)\) that
\[
  \begin{gathered}
  \|\operatorname{div} \mathbf W_\varepsilon\|_2
  + \left\|\left(\frac{1}{2} \boldsymbol y + \mathbf B\right) \operatorname{div} \mathbf W_\varepsilon\right\|_2
  \le C \varepsilon, \\
  \|\mathbf W_\varepsilon \operatorname{div} \mathbf W_\varepsilon\|_{n/(n-1)}
  \le \|\mathbf W_\varepsilon\|_{2n/(n-2)} \|\operatorname{div} \mathbf W_\varepsilon\|_2
  \le C \varepsilon.
  \end{gathered}
\]
By boundedness of the Riesz transforms \(\mathcal R_i\) on \(L^q(\mathbb R^n)\) (\(1 < q < \infty\)) and the Hardy--Littlewood--Sobolev inequality
\[
  \|(-\Delta)^{-\alpha/2} f\|_q \le C_{n,\alpha,p} \|f\|_p
  \qquad \left(\alpha \in (0, n), \quad 1 < p < \frac{n}{\alpha}, \quad \frac{1}{q} = \frac{1}{p} - \frac{\alpha}{n}\right),
\]
applied successively with \((\alpha, p) = (1, 2)\), \((1, n / (n - 1))\), and \((2, 2)\) (yielding \(q = 2n / (n - 2)\), \(n / (n - 2)\), and \(2n / (n - 4)\), where the last case uses \(n \ge 5 > 4\)), we obtain \(E_\varepsilon \in L^{2n/(n-2)}(\mathbb R^n) + L^{n/(n-2)}(\mathbb R^n) + L^{2n/(n-4)}(\mathbb R^n)\) with
\begin{equation}\label{eq:Eepsbound}
  \|E_\varepsilon\|_{L^{2n/(n-2)} + L^{n/(n-2)} + L^{2n/(n-4)}}
  \le C \varepsilon,
  \qquad
  \|E_\varepsilon\|_{L^{n/(n-2)}(B_{\ell_0})} \le C_{\ell_0} \varepsilon
  \quad (\ell_0 > 0).
\end{equation}

\smallskip
\noindent\emph{Step 2.}
Since \(-\Delta(\mathcal R_i \mathcal R_j f) = \partial_i \partial_j f\) in \(\mathscr{D}'(\mathbb R^n)\), combining \eqref{eq:divPeeps} with Weyl's lemma for harmonic distributions gives, in \(\mathscr{D}'(\mathbb R^n)\) and almost everywhere on \(\mathbb R^n\),
\begin{equation}\label{eq:Piexact}
  \Pi_\varepsilon
  = \Pi_B + \sum_{i=1}^n \sum_{j=1}^n \mathcal R_i \mathcal R_j
      (B_i W_{\varepsilon,j} + W_{\varepsilon,i} B_j + W_{\varepsilon,i} W_{\varepsilon,j})
    + E_\varepsilon + H_\varepsilon
\end{equation}
for a smooth harmonic function \(H_\varepsilon \in C^\infty(\mathbb R^n)\).
We now show that \(H_\varepsilon \equiv 0\).
Expanding \(\Pi_\varepsilon - \Pi_B\) from \eqref{eq:penaltydef} and \eqref{eq:approxmomentum} cancels \(\Pi_B\), \(\pm \frac{1}{2} |\mathbf B|^2\), and \(\mp \frac{1}{2} \boldsymbol y \cdot \mathbf B\):
\[
  \Pi_\varepsilon - \Pi_B
  = -\mathbf B \cdot \mathbf W_\varepsilon - \frac{1}{2} |\mathbf W_\varepsilon|^2
    + \frac{1}{2} \boldsymbol y \cdot \mathbf W_\varepsilon
    - \left(1 + \frac{\langle \boldsymbol y \rangle^2}{\varepsilon^2}\right) \operatorname{div} \mathbf W_\varepsilon.
\]
By \eqref{eq:approxbounds}, the last two terms belong to \(\langle \boldsymbol y \rangle L^2(\mathbb R^n)\).
For every fixed \(\boldsymbol y_0 \in \mathbb R^n\), Cauchy--Schwarz then shows that their \(L^1\) average over \(B_\ell(\boldsymbol y_0)\) vanishes as \(\ell \to \infty\):
\[
  \begin{aligned}
  &\frac{1}{|B_\ell(\boldsymbol y_0)|} \int_{B_\ell(\boldsymbol y_0)}
    \left|\frac{1}{2} \boldsymbol y \cdot \mathbf W_\varepsilon
          - \left(1 + \frac{\langle \boldsymbol y \rangle^2}{\varepsilon^2}\right) \operatorname{div} \mathbf W_\varepsilon\right| \mathrm d \boldsymbol y \\
  \le{}& \frac{C (\ell + |\boldsymbol y_0| + 1)}{\ell^{n/2}}
       \left(\|\mathbf W_\varepsilon\|_2
             + \left(1 + \frac{1}{\varepsilon^2}\right) \|\langle \boldsymbol y \rangle \operatorname{div} \mathbf W_\varepsilon\|_2\right)
   \longrightarrow 0.
  \end{aligned}
\]
Solving \eqref{eq:Piexact} for \(H_\varepsilon\) gives
\[
  H_\varepsilon
  = (\Pi_\varepsilon - \Pi_B)
    - \sum_{i=1}^n \sum_{j=1}^n \mathcal R_i \mathcal R_j
        (B_i W_{\varepsilon,j} + W_{\varepsilon,i} B_j + W_{\varepsilon,i} W_{\varepsilon,j})
    - E_\varepsilon.
\]
Beyond the two terms already estimated above, every remaining contribution to \(H_\varepsilon\) (namely \(-\mathbf B \cdot \mathbf W_\varepsilon - \frac{1}{2} |\mathbf W_\varepsilon|^2\), the double Riesz transforms, and \(E_\varepsilon\)) belongs to a finite sum of \(L^s(\mathbb R^n)\) spaces with \(1 < s < \infty\), and therefore also has vanishing \(L^1\) average on \(B_\ell(\boldsymbol y_0)\) as \(\ell \to \infty\).
By the mean-value property of harmonic functions, letting \(\ell \to \infty\) yields
\[
  H_\varepsilon(\boldsymbol y_0)
  = \frac{1}{|B_\ell(\boldsymbol y_0)|} \int_{B_\ell(\boldsymbol y_0)} H_\varepsilon(\boldsymbol y) \, \mathrm d \boldsymbol y
  \longrightarrow 0
\]
for every \(\boldsymbol y_0 \in \mathbb R^n\), whence \(H_\varepsilon \equiv 0\) on \(\mathbb R^n\).

\smallskip
\noindent\emph{Step 3.}
By \eqref{eq:approxbounds}, along a subsequence \(\varepsilon \downarrow 0\) we have \(\mathbf W_\varepsilon \rightharpoonup \mathbf W\) weakly in \(H^1(\mathbb R^n; \mathbb R^n)\) and strongly in \(L^q_{\mathrm{loc}}(\mathbb R^n)\) for all \(1 \le q < 2n / (n - 2)\), with \(\operatorname{div} \mathbf W = 0\).
Fix \(1 < p < n / (n - 2)\) and a ball \(B_{\ell_0}\).
For any \(\ell > 2\ell_0\) and smooth cutoff \(\eta_\ell \in C_c^\infty(B_{2\ell})\) with \(\eta_\ell|_{B_\ell} = 1\), the Rellich--Kondrachov compactness theorem gives strong \(L^p(\mathbb R^n)\) convergence of \(\eta_\ell (B_i W_{\varepsilon,j} + W_{\varepsilon,i} B_j + W_{\varepsilon,i} W_{\varepsilon,j})\) since \(2p < 2n / (n - 2)\).
On the exterior \(\{|\boldsymbol z| \ge \ell > 2\ell_0\}\), every \(\boldsymbol y \in B_{\ell_0}\) satisfies \(|\boldsymbol y - \boldsymbol z| \ge |\boldsymbol z| / 2\), so the Calder\'on--Zygmund kernel bound \(O(|\boldsymbol y - \boldsymbol z|^{-n})\) yields
\[
  \begin{aligned}
  &\Bigl\|\mathcal R_i \mathcal R_j \bigl[(1 - \eta_\ell)
     (B_i W_{\varepsilon,j} + W_{\varepsilon,i} B_j + W_{\varepsilon,i} W_{\varepsilon,j})\bigr]\Bigr\|_{L^\infty(B_{\ell_0})} \\
  \le{}& \frac{C}{\ell^n} \|\mathbf W_\varepsilon\|_2^2
       + \frac{C}{\ell^{n/2}} \|\mathbf B \otimes \mathbf W_\varepsilon\|_2
   \le C \left(\frac{1}{\ell^n} + \frac{1}{\ell^{n/2}}\right).
  \end{aligned}
\]
Sending \(\varepsilon \downarrow 0\) and then \(\ell \to \infty\) in \eqref{eq:Piexact} (with \(H_\varepsilon \equiv 0\)) and using \eqref{eq:Eepsbound}, we obtain \(\Pi_\varepsilon \to \Pi\) strongly in \(L^p(B_{\ell_0})\).
Expressing \(\Phi_\varepsilon\) via \eqref{eq:approxmomentum} as
\[
  \Phi_\varepsilon
  = \Pi_\varepsilon + \frac{1}{2} |\mathbf U_\varepsilon|^2 - \frac{1}{2} \boldsymbol y \cdot \mathbf U_\varepsilon
\]
and using \(\mathbf U_\varepsilon \to \mathbf U\) strongly in \(L^{2p}(B_{\ell_0}; \mathbb R^n)\), we obtain \(\Phi_\varepsilon \to \Phi\) strongly in \(L^p(B_{\ell_0})\) as well.

\smallskip
\noindent\emph{Step 4.}
In \eqref{eq:approxmomentum}, we rewrite the viscous and convective terms via \eqref{eq:curlcurl} as
\[
  -\operatorname{div} \Omega_{\mathbf U_\varepsilon}
  = -\Delta \mathbf U_\varepsilon + \nabla \operatorname{div} \mathbf W_\varepsilon,
  \qquad
  (\mathbf U_\varepsilon \cdot \nabla) \mathbf U_\varepsilon
  = \operatorname{div}(\mathbf U_\varepsilon \otimes \mathbf U_\varepsilon)
    - \mathbf U_\varepsilon \operatorname{div} \mathbf W_\varepsilon,
\]
so passing to the distributional limit shows that \((\mathbf U, \Pi)\) solves \eqref{eq:profile}.
Extracting a further subsequence along which \(\Phi_\varepsilon \to \Phi\) almost everywhere passes \eqref{eq:epsmax} to \eqref{eq:mainmax}.
Finally, to pass to the limit in \eqref{eq:approxenergy}, we drop the nonnegative penalty term and bound the divergence cross term using \(|\mathbf B(\boldsymbol y) - \frac{1}{2} \boldsymbol y| \le (1 + \|\mathbf B\|_\infty) \langle \boldsymbol y \rangle\):
\[
  \left|\int_{\mathbb R^n}
    \left(\left(\mathbf B - \frac{1}{2} \boldsymbol y\right) \cdot \mathbf W_\varepsilon\right)
    \operatorname{div} \mathbf W_\varepsilon \, \mathrm d \boldsymbol y\right|
  \le (1 + \|\mathbf B\|_\infty) \|\mathbf W_\varepsilon\|_2
      \|\langle \boldsymbol y \rangle \operatorname{div} \mathbf W_\varepsilon\|_2
  \le C \varepsilon \longrightarrow 0.
\]
Combining strong \(L^2_{\mathrm{loc}}\) convergence of \(\mathbf W_\varepsilon\) with \(|\nabla \mathbf B(\boldsymbol y)| \le C \langle \boldsymbol y \rangle^{-2}\) gives
\[
  \int_{\mathbb R^n} ((\mathbf W_\varepsilon \cdot \nabla) \mathbf B) \cdot \mathbf W_\varepsilon \, \mathrm d \boldsymbol y
  \longrightarrow \int_{\mathbb R^n} ((\mathbf W \cdot \nabla) \mathbf B) \cdot \mathbf W \, \mathrm d \boldsymbol y,
\]
and weak lower semicontinuity of \(\|\nabla \mathbf W_\varepsilon\|_2^2 + \frac{n-2}{4} \|\mathbf W_\varepsilon\|_2^2\) yields \eqref{eq:mainenergy}.
\end{proof}

\subsection{Global smoothness via uniaxial straining lift}\label{subsec:regularity}

To prove that the weak profile from Proposition~\ref{prop:limitprofile} is smooth, we lift \eqref{eq:profile} to an unforced stationary Navier--Stokes system in \(\mathbb R^{n+1}\) by superimposing an irrotational uniaxial straining flow.
Global smoothness and the strong \(H^1_{\mathrm{loc}}\) initial trace then follow from the quantitative stationary head-pressure regularity result in Appendix~\ref{app:stationary} (Proposition~\ref{thm:stationaryreg}; cf.\ \cite{FR94maximum,FR94bounded,FR98criterion,LY22,GLLWX26transverse}).

\begin{proposition}\label{prop:smoothprofile}
Let \(n \ge 5\), let \(\mathbf u_0\) satisfy the hypotheses of Theorem~\ref{thm:main}, and let \(\mathbf B = \mathbf u_R(1)\) and \(\Phi_B\) be as in Proposition~\ref{prop:defect} with \(R \ge 1\) and \(R^{-2} \le T_0\).
Suppose \(\mathbf W \in H^1_\sigma(\mathbb R^n)\) is such that \(\mathbf U = \mathbf B + \mathbf W\) and \(\Pi\) defined in \eqref{eq:mainpressure} solve \eqref{eq:profile} in \(\mathscr{D}'(\mathbb R^n)\) and satisfy \eqref{eq:mainmax} almost everywhere on \(\mathbb R^n\).
Then \((\mathbf U, \Pi)\) admits a representative in \(C^\infty(\mathbb R^n; \mathbb R^n) \times C^\infty(\mathbb R^n)\) solving \eqref{eq:profile} classically and satisfying \eqref{eq:mainmax} everywhere on \(\mathbb R^n\).
Moreover, the pair \((\mathbf u, p)\) defined by \eqref{eq:similarity} belongs to \(C^\infty((0, \infty) \times \mathbb R^n)\), solves \eqref{eq:NS} classically for \(t > 0\), and satisfies \eqref{eq:initialtrace} for every \(\ell > 0\).
\end{proposition}

\begin{proof}
The harmonic polynomial \(-\frac{1}{4} |\boldsymbol y|^2 + \frac{n}{4} y_{n+1}^2\) on \(\mathbb R^{n+1}\) generates the divergence-free uniaxial straining flow \(\mathbf V_0(\boldsymbol y, y_{n+1}) := (-\frac{1}{2} \boldsymbol y, \frac{n}{2} y_{n+1})\), which solves the unforced stationary Navier--Stokes equations with pressure \(\widetilde \Pi_0 := -\frac{1}{2} |\mathbf V_0|^2\).
Superimposing the \(y_{n+1}\)-independent profile \((\mathbf U(\boldsymbol y), 0)\) onto \(\mathbf V_0\), we define on \(\mathbb R^{n+1} = \mathbb R^n \times \mathbb R\):
\begin{equation}\label{eq:lift}
  \begin{aligned}
  \mathbf V(\boldsymbol y, y_{n+1})
  &:= \left(\mathbf U(\boldsymbol y) - \frac{1}{2} \boldsymbol y, \, \frac{n}{2} y_{n+1}\right), \\
  \widetilde \Pi(\boldsymbol y, y_{n+1})
  &:= \Pi(\boldsymbol y) - \frac{1}{8} |\boldsymbol y|^2 - \frac{n^2}{8} y_{n+1}^2.
  \end{aligned}
\end{equation}
By Fubini's theorem, \(\mathbf V \in H^1_{\mathrm{loc}}(\mathbb R^{n+1}; \mathbb R^{n+1})\) and \(\widetilde \Pi \in L^1_{\mathrm{loc}}(\mathbb R^{n+1})\), and \(\mathbf V\) is divergence-free:
\[
  \operatorname{div}_{\boldsymbol y, y_{n+1}} \mathbf V
  = \operatorname{div} \mathbf U - \frac{n}{2} + \frac{n}{2} = 0.
\]
As the first \(n\) components \(\mathbf U(\boldsymbol y) - \frac{1}{2} \boldsymbol y\) of \(\mathbf V\) are independent of \(y_{n+1}\), direct expansion gives
\[
  \begin{aligned}
  &\left(\left(\mathbf U - \frac{1}{2} \boldsymbol y\right) \cdot \nabla_{\boldsymbol y}\right)
   \left(\mathbf U - \frac{1}{2} \boldsymbol y\right)
   + \nabla_{\boldsymbol y} \widetilde \Pi \\
  ={}& (\mathbf U \cdot \nabla_{\boldsymbol y}) \mathbf U
     - \frac{1}{2} \boldsymbol y \cdot \nabla_{\boldsymbol y} \mathbf U
     - \frac{1}{2} \mathbf U + \frac{1}{4} \boldsymbol y
     + \nabla_{\boldsymbol y} \Pi - \frac{1}{4} \boldsymbol y \\
  ={}& -\frac{1}{2} \boldsymbol y \cdot \nabla_{\boldsymbol y} \mathbf U
     - \frac{1}{2} \mathbf U
     + (\mathbf U \cdot \nabla_{\boldsymbol y}) \mathbf U
     + \nabla_{\boldsymbol y} \Pi,
  \end{aligned}
\]
so the first \(n\) momentum equations of \(-\Delta \mathbf V + (\mathbf V \cdot \nabla) \mathbf V + \nabla \widetilde \Pi = \boldsymbol 0\) on \(\mathbb R^{n+1}\) reduce to \eqref{eq:profile}.
In the \((n + 1)\)-th component, we have \(\Delta V_{n+1} = 0\), \(\nabla_{\boldsymbol y} V_{n+1} = \boldsymbol 0\), and
\[
  V_{n+1} \partial_{y_{n+1}} V_{n+1} + \partial_{y_{n+1}} \widetilde \Pi
  = \frac{n^2}{4} y_{n+1} - \frac{n^2}{4} y_{n+1} = 0,
\]
so the \((n + 1)\)-th momentum equation holds identically.
Moreover, expanding
\[
  \frac{1}{2} |\mathbf V|^2
  = \frac{1}{2} |\mathbf U|^2 - \frac{1}{2} \boldsymbol y \cdot \mathbf U
    + \frac{1}{8} |\boldsymbol y|^2 + \frac{n^2}{8} y_{n+1}^2
\]
shows that the quadratic straining terms cancel in the steady total head pressure:
\begin{equation}\label{eq:lifthead}
  \begin{aligned}
  \widetilde \Pi(\boldsymbol y, y_{n+1}) + \frac{1}{2} |\mathbf V(\boldsymbol y, y_{n+1})|^2
  &= \Pi(\boldsymbol y) + \frac{1}{2} |\mathbf U(\boldsymbol y)|^2
     - \frac{1}{2} \boldsymbol y \cdot \mathbf U(\boldsymbol y) \\
  &= \Phi(\boldsymbol y) \le \sup_{\mathbb R^n} \Phi_B.
  \end{aligned}
\end{equation}
By \eqref{eq:lifthead}, applying Proposition~\ref{thm:stationaryreg} with \(m = n + 1 \ge 6\) yields \(\mathbf V \in C^\infty(\mathbb R^{n+1}; \mathbb R^{n+1})\) and \(\widetilde \Pi \in C^\infty(\mathbb R^{n+1})\), whence \((\mathbf U, \Pi) \in C^\infty(\mathbb R^n; \mathbb R^n) \times C^\infty(\mathbb R^n)\) solves \eqref{eq:profile} classically and satisfies \eqref{eq:mainmax} everywhere on \(\mathbb R^n\).
Consequently, \((\mathbf u, p)\) defined by \eqref{eq:similarity} belongs to \(C^\infty((0, \infty) \times \mathbb R^n)\) and solves \eqref{eq:NS} classically for \(t > 0\).

To verify \eqref{eq:initialtrace}, we express \(\mathbf B - \mathbf u_{0,R}\) via \eqref{eq:rescale}, \eqref{eq:mild}, and the fundamental theorem of calculus for \(s \mapsto e^{s\Delta} \mathbf u_{0,R}\) on \([0, 1]\): for \(k \in \{0, 1\}\),
\[
  \nabla^k (\mathbf B - \mathbf u_{0,R})
  = \int_0^1 e^{s\Delta} \nabla^k \Delta \mathbf u_{0,R} \, \mathrm d s
    - \int_0^1 \mathcal K_{1-s} * \nabla^k (\mathbf u_R \otimes \mathbf u_R)(s) \, \mathrm d s.
\]
By \eqref{eq:cutoffbounds} and \eqref{eq:bgbounds}, we have \(|\nabla^k \Delta \mathbf u_{0,R}(\boldsymbol y)| \le C \langle \boldsymbol y \rangle^{-3-k}\) (with \(k = 1\) using three derivatives of \(\mathbf u_0\)) and \(|\nabla^k (\mathbf u_R \otimes \mathbf u_R)(s, \boldsymbol y)| \le C \langle \boldsymbol y \rangle^{-2-k}\) uniformly for \(0 \le s \le 1\).
Since \(2 + k < 3 + k \le 4 < n + 1\) for \(n \ge 5\), \eqref{eq:kernelweighted} in Lemma~\ref{lem:kernel} yields
\begin{equation}\label{eq:bgdiff}
  |\nabla^k (\mathbf B - \mathbf u_{0,R})(\boldsymbol y)|
  \le \frac{C}{\langle \boldsymbol y \rangle^{k+2}},
  \qquad k \in \{0, 1\}.
\end{equation}
Moreover, for \(k = 0\) the decay in \eqref{eq:bgdiff} improves to
\begin{equation}\label{eq:bgdiff3}
  |\mathbf B(\boldsymbol y) - \mathbf u_{0,R}(\boldsymbol y)|
  \le \frac{C}{\langle \boldsymbol y \rangle^3},
  \qquad \boldsymbol y \in \mathbb R^n.
\end{equation}
Indeed, \eqref{eq:kernelweighted} at \(q = 3\) already bounds \(\int_0^1 |e^{s\Delta} \Delta \mathbf u_{0,R}| \, \mathrm d s\) by \(C \langle \boldsymbol y \rangle^{-3}\).
For the nonlinear term at \(s \in [0, 1)\) with \(\tau := 1 - s \in (0, 1]\), factoring \(\partial_j\) out of \eqref{eq:kerneldef} and integrating by parts using \(\operatorname{div}(\mathbf u_R \otimes \mathbf u_R) = (\mathbf u_R \cdot \nabla) \mathbf u_R\) gives
\[
  \begin{aligned}
  \bigl(\mathcal K_\tau * (\mathbf u_R \otimes \mathbf u_R)(s)\bigr)_i(\boldsymbol y)
  ={}& \bigl(e^{\tau\Delta} ((\mathbf u_R \cdot \nabla) \mathbf u_R)_i(s)\bigr)(\boldsymbol y) \\
     &+ \sum_{k=1}^n \int_{\mathbb R^n}
          \left(\int_\tau^\infty \partial_i \partial_k G_\sigma(\boldsymbol y - \boldsymbol z) \, \mathrm d \sigma\right)
          \bigl((\mathbf u_R \cdot \nabla) \mathbf u_R\bigr)_k(s, \boldsymbol z) \, \mathrm d \boldsymbol z.
  \end{aligned}
\]
By \eqref{eq:bgbounds} (\(k = 0, 1, 2\)), we have uniformly in \(s \in [0, 1]\),
\[
  \bigl|(\mathbf u_R \cdot \nabla) \mathbf u_R(s, \boldsymbol z)\bigr|
  \le \frac{C}{\langle \boldsymbol z \rangle^3},
  \qquad
  \bigl|\nabla \bigl((\mathbf u_R \cdot \nabla) \mathbf u_R\bigr)(s, \boldsymbol z)\bigr|
  \le \frac{C}{\langle \boldsymbol z \rangle^4},
\]
so \eqref{eq:kernelweighted} bounds the first term by \(C \langle \boldsymbol y \rangle^{-3}\).
For the second term, the Gaussian derivative bound \(|\partial_i \partial_k G_\sigma(\boldsymbol z)| \le C (|\boldsymbol z| + \sqrt \sigma)^{-n-2}\) yields
\[
  \left|\int_\tau^\infty \partial_i \partial_k G_\sigma(\boldsymbol z) \, \mathrm d \sigma\right|
  \le \frac{C}{(|\boldsymbol z| + \sqrt \tau)^n},
\]
while Fubini's theorem, radial symmetry of \(G_\sigma\), and \(\Delta G_\sigma = \partial_\sigma G_\sigma\) give, for every \(r > 0\),
\[
  \begin{aligned}
  \int_{B_r} \int_\tau^\infty \partial_i \partial_k G_\sigma(\boldsymbol z) \, \mathrm d \sigma \, \mathrm d \boldsymbol z
  &= \frac{\delta_{ik}}{n} \int_\tau^\infty \frac{\mathrm d}{\mathrm d \sigma} \left(\int_{B_r} G_\sigma(\boldsymbol z) \, \mathrm d \boldsymbol z\right) \mathrm d \sigma \\
  &= -\frac{\delta_{ik}}{n} \int_{B_r} G_\tau(\boldsymbol z) \, \mathrm d \boldsymbol z
   \in \left[-\frac{1}{n}, \, 0\right].
  \end{aligned}
\]
Now fix \(|\boldsymbol y| \ge 2\) and split the \(\boldsymbol z\)-integral in the second term into three regions: \(B_{|\boldsymbol y|/2}(\boldsymbol y)\), \(B_{2|\boldsymbol y|} \setminus B_{|\boldsymbol y|/2}(\boldsymbol y)\), and \(\mathbb R^n \setminus B_{2|\boldsymbol y|}\).
On \(B_{|\boldsymbol y|/2}(\boldsymbol y)\), we subtract the value of \((\mathbf u_R \cdot \nabla) \mathbf u_R\) at \(\boldsymbol y\) using the uniform ball-integral bound above and the \(O(|\boldsymbol y|^{-4})\) gradient bound; on \(B_{2|\boldsymbol y|} \setminus B_{|\boldsymbol y|/2}(\boldsymbol y)\) we use \(|\boldsymbol y - \boldsymbol z| \ge |\boldsymbol y| / 2\) and \(n \ge 5 > 3\); on \(\mathbb R^n \setminus B_{2|\boldsymbol y|}\) we use \(|\boldsymbol y - \boldsymbol z| \ge |\boldsymbol z| / 2\).
This gives, uniformly in \(\tau \in (0, 1]\),
\[
  \begin{aligned}
  &\left|\int_{B_{|\boldsymbol y|/2}(\boldsymbol y)}
     \left(\int_\tau^\infty \partial_i \partial_k G_\sigma(\boldsymbol y - \boldsymbol z) \, \mathrm d \sigma\right)
     \bigl((\mathbf u_R \cdot \nabla) \mathbf u_R\bigr)_k(s, \boldsymbol z) \, \mathrm d \boldsymbol z\right| \\
  \le{}& \frac{C}{|\boldsymbol y|^3}
         + \frac{C}{|\boldsymbol y|^4} \int_{B_{|\boldsymbol y|/2}} \frac{\mathrm d \boldsymbol z}{|\boldsymbol z|^{n-1}}
   \le \frac{C}{|\boldsymbol y|^3},
  \end{aligned}
\]
together with
\[
  \begin{aligned}
  &\left|\int_{B_{2|\boldsymbol y|} \setminus B_{|\boldsymbol y|/2}(\boldsymbol y)}
     \left(\int_\tau^\infty \partial_i \partial_k G_\sigma(\boldsymbol y - \boldsymbol z) \, \mathrm d \sigma\right)
     \bigl((\mathbf u_R \cdot \nabla) \mathbf u_R\bigr)_k(s, \boldsymbol z) \, \mathrm d \boldsymbol z\right| \\
  \le{}& \frac{C}{|\boldsymbol y|^n} \int_{B_{2|\boldsymbol y|}} \frac{\mathrm d \boldsymbol z}{\langle \boldsymbol z \rangle^3}
   \le \frac{C}{|\boldsymbol y|^3},
  \end{aligned}
\]
and
\[
  \begin{aligned}
  &\left|\int_{\mathbb R^n \setminus B_{2|\boldsymbol y|}}
     \left(\int_\tau^\infty \partial_i \partial_k G_\sigma(\boldsymbol y - \boldsymbol z) \, \mathrm d \sigma\right)
     \bigl((\mathbf u_R \cdot \nabla) \mathbf u_R\bigr)_k(s, \boldsymbol z) \, \mathrm d \boldsymbol z\right| \\
  \le{}& C \int_{|\boldsymbol z| \ge 2|\boldsymbol y|} \frac{\mathrm d \boldsymbol z}{|\boldsymbol z|^{n+3}}
   \le \frac{C}{|\boldsymbol y|^3}.
  \end{aligned}
\]
Combining these three bounds with \eqref{eq:bgdiff} on \(|\boldsymbol y| < 2\) and integrating over \(s \in [0, 1]\) establishes \eqref{eq:bgdiff3}.
For any \(\boldsymbol x \ne \boldsymbol 0\) and \(0 < t < |\boldsymbol x|^2 / (4R^2)\), we have \(\mathbf u_{0,R}(\boldsymbol x / \sqrt t) = \mathbf u_0(\boldsymbol x / \sqrt t) = \sqrt t \, \mathbf u_0(\boldsymbol x)\), so \eqref{eq:bgdiff} gives, for \(k \in \{0, 1\}\),
\[
  \left|\nabla_{\boldsymbol x}^k
    \left[\frac{1}{\sqrt t} \mathbf B\left(\frac{\boldsymbol x}{\sqrt t}\right)
          - \mathbf u_0(\boldsymbol x)\right]\right|
  \le \frac{C \sqrt t}{(\sqrt t + |\boldsymbol x|)^{k+2}}
  \le \frac{C \sqrt t}{|\boldsymbol x|^{k+2}} \longrightarrow 0
  \quad \text{as } t \downarrow 0.
\]
Moreover, for all \(t > 0\) and \(\boldsymbol x \ne \boldsymbol 0\), combining \eqref{eq:bgbounds} with the \((-1)\)-homogeneity of \(\mathbf u_0\) gives the uniform integrable majorant
\[
  \begin{aligned}
  \left|\nabla_{\boldsymbol x}^k
    \left[\frac{1}{\sqrt t} \mathbf B\left(\frac{\boldsymbol x}{\sqrt t}\right)
          - \mathbf u_0(\boldsymbol x)\right]\right|
  &\le \frac{1}{t^{(k+1)/2}} \left|(\nabla^k \mathbf B)\left(\frac{\boldsymbol x}{\sqrt t}\right)\right|
       + |\nabla^k \mathbf u_0(\boldsymbol x)| \\
  &\le \frac{C}{(\sqrt t + |\boldsymbol x|)^{k+1}} + \frac{C}{|\boldsymbol x|^{k+1}}
   \le \frac{C}{|\boldsymbol x|^{k+1}} \in L^2(B_\ell)
  \end{aligned}
\]
since \(n \ge 5\).
As \(|\mathbf u_0(\boldsymbol x)| \le C |\boldsymbol x|^{-1}\) and \(n \ge 5\), the classical gradient of \(\mathbf u_0\) on \(\mathbb R^n \setminus \{\boldsymbol 0\}\) coincides with its distributional gradient on \(\mathbb R^n\), and dominated convergence yields \(t^{-1/2} \mathbf B(\boldsymbol x / \sqrt t) \to \mathbf u_0(\boldsymbol x)\) strongly in \(H^1(B_\ell; \mathbb R^n)\).
For \(\mathbf W \in H^1_\sigma(\mathbb R^n)\), changing variables gives as \(t \downarrow 0\) (using \(n - 4 > 0\)):
\[
  \begin{aligned}
  \left\|\frac{1}{\sqrt t} \mathbf W\left(\frac{\cdot}{\sqrt t}\right)\right\|_{L^2(B_\ell)}
  &\le t^{(n-2)/4} \|\mathbf W\|_2 \longrightarrow 0, \\
  \left\|\nabla_{\boldsymbol x} \left[\frac{1}{\sqrt t} \mathbf W\left(\frac{\cdot}{\sqrt t}\right)\right]\right\|_{L^2(B_\ell)}
  &\le t^{(n-4)/4} \|\nabla \mathbf W\|_2 \longrightarrow 0,
  \end{aligned}
\]
which establishes \eqref{eq:initialtrace}.
\end{proof}

\begin{remark}\label{rem:liftliouville}
The lifted solution \((\mathbf V, \widetilde \Pi)\) in \eqref{eq:lift} grows linearly along \(y_{n+1}\) as a perturbation of the uniaxial straining flow \(\mathbf V_0 = (-\frac{1}{2} \boldsymbol y, \frac{n}{2} y_{n+1})\), so decay-based stationary Liouville theorems on \(\mathbb R^{n+1}\) do not apply, whereas Proposition~\ref{thm:stationaryreg} is a purely local interior estimate.
Moreover, the quadratic straining pressure in \eqref{eq:lift} cancels \(\frac{1}{2} |\mathbf V_0|^2\), so the total head pressure \(\widetilde \Pi + \frac{1}{2} |\mathbf V|^2 = \Phi\) is uniformly bounded above on \(\mathbb R^{n+1}\).
\end{remark}

\section{Spatial decay of the profile}\label{sec:decay}

\subsection{Local Morrey and weighted energy estimates}\label{subsec:weighted}

\begin{proposition}\label{prop:subquadmorrey}
Let \(n \ge 5\), let \(\mathbf u_0\) satisfy the hypotheses of Theorem~\ref{thm:main}, and let \(\mathbf B = \mathbf u_R(1)\) and \(\Phi_B\) be as in Proposition~\ref{prop:defect} with \(R \ge 1\) and \(R^{-2} \le T_0\).
Suppose \((\mathbf U, \Pi) \in C^\infty(\mathbb R^n; \mathbb R^n) \times C^\infty(\mathbb R^n)\) solves \eqref{eq:profile} on \(\mathbb R^n\) with \(\mathbf U = \mathbf B + \mathbf W\), \(\mathbf W \in H^1_\sigma(\mathbb R^n)\), \(\Pi\) given by \eqref{eq:mainpressure}, and \(\Phi\) satisfying \eqref{eq:mainmax}.
Then there exists a constant \(C > 0\), depending on \(n\), \(\mathbf u_0\), \(R\), and \(\|\mathbf W\|_{H^1}\) (and on \(n\), \(\mathbf u_0\), \(R\) alone for the profile constructed in Proposition~\ref{prop:limitprofile}, by \eqref{eq:mainenergy} and \eqref{eq:coercivity}), such that for every \(\boldsymbol y_0 \in \mathbb R^n\),
\begin{align}
  \|\mathbf U\|_{L^{2,n-2}(B_1(\boldsymbol y_0))}^2
  &\le C \langle \boldsymbol y_0 \rangle^{2 - 4/n}, \label{eq:subquadmorrey} \\
  \int_{B_{1/4}(\boldsymbol y_0)} |\mathbf W(\boldsymbol y)|^4 \, \mathrm d \boldsymbol y
  &\le C \langle \boldsymbol y_0 \rangle^{2 - 4/n}
      \int_{B_{1/2}(\boldsymbol y_0)}
        \bigl(|\nabla \mathbf W(\boldsymbol y)|^2 + |\mathbf W(\boldsymbol y)|^2\bigr) \, \mathrm d \boldsymbol y. \label{eq:unitballL4}
\end{align}
\end{proposition}

\begin{proof}
In the pressure decomposition \eqref{eq:mainpressure},
\[
  \Pi - \Pi_B
  = \sum_{i=1}^n \sum_{j=1}^n \mathcal R_i \mathcal R_j (B_i W_j + W_i B_j)
    + \sum_{i=1}^n \sum_{j=1}^n \mathcal R_i \mathcal R_j (W_i W_j),
\]
the linear remainder terms lie in \(L^2(\mathbb R^n)\) and the quadratic term lies in \(L^{n/(n-2)}(\mathbb R^n)\).
Combined with \(\Pi_B \in L^\infty(\mathbb R^n)\) and \(\mathbf U = \mathbf B + \mathbf W \in L^\infty(\mathbb R^n) + L^2(\mathbb R^n)\), H\"older's inequality yields
\begin{equation}\label{eq:uniformB2}
  \sup_{\boldsymbol y_0 \in \mathbb R^n}
  \int_{B_2(\boldsymbol y_0)} \bigl(|\mathbf U(\boldsymbol y)|^2 + |\Pi(\boldsymbol y)|\bigr) \, \mathrm d \boldsymbol y
  \le C.
\end{equation}

Fix \(\boldsymbol y_0 \in \mathbb R^n\).
Since \(\operatorname{div} \mathbf U = 0\) implies \(\operatorname{div}(\boldsymbol y \cdot \nabla \mathbf U + \mathbf U) = (\boldsymbol y \cdot \nabla + 2) \operatorname{div} \mathbf U = 0\), taking the divergence of \eqref{eq:profile} gives
\[
  -\Delta \Pi = \sum_{i=1}^n \sum_{j=1}^n \partial_i \partial_j (U_i U_j).
\]
Hence \((\mathbf U, \Pi)\) satisfies \eqref{eq:weakpreseq} on \(B_2(\boldsymbol y_0)\) with steady total head pressure \(\varpi := \Pi + \frac{1}{2} |\mathbf U|^2\).
By \eqref{eq:mainmax} and \(\sup_{\mathbb R^n} \Phi_B \ge 0\), on \(B_2(\boldsymbol y_0)\) we have
\[
  \varpi_+(\boldsymbol y)
  = \left(\Pi(\boldsymbol y) + \frac{1}{2} |\mathbf U(\boldsymbol y)|^2\right)_+
  = \left(\Phi(\boldsymbol y) + \frac{1}{2} \boldsymbol y \cdot \mathbf U(\boldsymbol y)\right)_+
  \le \sup_{\mathbb R^n} \Phi_B + C \langle \boldsymbol y_0 \rangle |\mathbf U(\boldsymbol y)|.
\]

Let \(1 \le r_1 < r_2 \le 2\).
Applying \eqref{eq:prestestmorrey} in Lemma~\ref{lem:prestest} to \((\mathbf U, \Pi)\) on \(B_{r_1}(\boldsymbol y_0) \subset B_{r_2}(\boldsymbol y_0)\) together with \eqref{eq:uniformB2} yields
\begin{equation}\label{eq:Thetaraw}
  \|\mathbf U\|_{L^{2,n-2}(B_{r_1}(\boldsymbol y_0))}^2
  \le \frac{C}{(r_2 - r_1)^{n-2}} 
  + C \langle \boldsymbol y_0 \rangle
        \sup_{\boldsymbol z \in B_{(r_1+r_2)/2}(\boldsymbol y_0)}
        \int_{B_{r_2}(\boldsymbol y_0)}
          \frac{|\mathbf U(\boldsymbol y)|}{|\boldsymbol y - \boldsymbol z|^{n-2}} \, \mathrm d \boldsymbol y.
\end{equation}

For any \(\boldsymbol z \in B_{(r_1+r_2)/2}(\boldsymbol y_0)\) and any splitting radius \(0 < \delta \le (r_2 - r_1)/8\), we have \(B_\delta(\boldsymbol z) \subset B_{r_2}(\boldsymbol y_0)\).
Splitting the integral in \eqref{eq:Thetaraw} at \(|\boldsymbol y - \boldsymbol z| = \delta\), decomposing \(B_\delta(\boldsymbol z)\) into dyadic annuli, and applying Cauchy--Schwarz together with \eqref{eq:uniformB2} outside \(B_\delta(\boldsymbol z)\) gives
\[
  \begin{aligned}
  \int_{B_\delta(\boldsymbol z)}
     \frac{|\mathbf U(\boldsymbol y)|}{|\boldsymbol y - \boldsymbol z|^{n-2}} \, \mathrm d \boldsymbol y
   &= \sum_{k=0}^\infty
      \int_{ \frac{\delta}{2^{k+1}} \le |\boldsymbol y - \boldsymbol z| < \frac{\delta}{2^{k}} }
        \frac{|\mathbf U(\boldsymbol y)|}{|\boldsymbol y - \boldsymbol z|^{n-2}} \, \mathrm d \boldsymbol y
    \le C \delta \|\mathbf U\|_{L^{2,n-2}(B_{r_2}(\boldsymbol y_0))}, \\
  \int_{B_{r_2}(\boldsymbol y_0) \setminus B_\delta(\boldsymbol z)}
     \frac{|\mathbf U(\boldsymbol y)|}{|\boldsymbol y - \boldsymbol z|^{n-2}} \, \mathrm d \boldsymbol y
   &\le \|\mathbf U\|_{L^2(B_2(\boldsymbol y_0))}
        \left(\int_{|\boldsymbol y - \boldsymbol z| \ge \delta}
          \frac{\mathrm d \boldsymbol y}{|\boldsymbol y - \boldsymbol z|^{2n - 4}} \right)^{1/2}
    \le \frac{C}{\delta^{(n-4)/2}}.
  \end{aligned}
\]
In \eqref{eq:Thetaraw}, applying Young's inequality to the inner-ball contribution produces a term of order \(\langle \boldsymbol y_0 \rangle^2 \delta^2\) together with the outer-ball contribution \(\langle \boldsymbol y_0 \rangle \delta^{-(n-4)/2}\).
To balance these two powers of \(\langle \boldsymbol y_0 \rangle\) via \(\delta^{n/2} \sim \langle \boldsymbol y_0 \rangle^{-1}\), we choose
\[
  \delta := \frac{r_2 - r_1}{8 \langle \boldsymbol y_0 \rangle^{2/n}}
  \in \left(0, \, \frac{r_2 - r_1}{8}\right].
\]
Substituting both bounds with this choice of \(\delta\) into \eqref{eq:Thetaraw} and applying Young's inequality gives, for every \(\theta_0 > 0\) and \(1 \le r_1 < r_2 \le 2\),
\[
  \|\mathbf U\|_{L^{2,n-2}(B_{r_1}(\boldsymbol y_0))}^2
  \le \theta_0 \|\mathbf U\|_{L^{2,n-2}(B_{r_2}(\boldsymbol y_0))}^2
      + \frac{C_{\theta_0}}{(r_2 - r_1)^{n-2}} \langle \boldsymbol y_0 \rangle^{2 - 4/n}.
\]
Choosing \(\theta_0 = 2^{-(n-1)}\) and iterating along the concentric radii \(r_j = 2 - 2^{-j}\) (\(j \ge 0\), with \(r_0 = 1\) and \(r_{j+1} - r_j = 2^{-(j+1)}\)) gives, for every \(k \ge 1\),
\[
  \begin{aligned}
  \|\mathbf U\|_{L^{2,n-2}(B_1(\boldsymbol y_0))}^2
  &\le \theta_0^k \|\mathbf U\|_{L^{2,n-2}(B_2(\boldsymbol y_0))}^2
       + C \langle \boldsymbol y_0 \rangle^{2 - 4/n}
         \sum_{j=0}^{k-1} \frac{2^{(j+1)(n-2)}}{2^{(n-1)j}} \\
  &= \frac{1}{2^{(n-1)k}} \|\mathbf U\|_{L^{2,n-2}(B_2(\boldsymbol y_0))}^2
     + C \, 2^{n-2} \langle \boldsymbol y_0 \rangle^{2 - 4/n}
       \sum_{j=0}^{k-1} \frac{1}{2^j}.
  \end{aligned}
\]
Because \(\|\mathbf U\|_{L^{2,n-2}(B_2(\boldsymbol y_0))} < \infty\) by smoothness of \(\mathbf U\), sending \(k \to \infty\) proves \eqref{eq:subquadmorrey}.

Since \(\mathbf W = \mathbf U - \mathbf B\), \(B_{1/2}(\boldsymbol y_0) \subset B_1(\boldsymbol y_0)\), and \(\mathbf B \in L^\infty(\mathbb R^n)\), \eqref{eq:subquadmorrey} implies
\[
  \|\mathbf W\|_{L^{2,n-2}(B_{1/2}(\boldsymbol y_0))}^2
  \le 2 \|\mathbf U\|_{L^{2,n-2}(B_1(\boldsymbol y_0))}^2 + C \|\mathbf B\|_\infty^2
  \le C \langle \boldsymbol y_0 \rangle^{2 - 4/n}.
\]
Applying the concentric-ball estimate in Lemma~\ref{lem:morreyl4} to \(\mathbf W\) on \(B_{1/4}(\boldsymbol y_0) \subset B_{1/2}(\boldsymbol y_0)\) yields \eqref{eq:unitballL4}.
\end{proof}

To propagate weighted \(H^1(\mathbb R^n)\) bounds for \(\langle \boldsymbol y \rangle^\beta \mathbf W\) across \(0 \le \beta < (n + 2)/2\) starting only from \(\mathbf W \in H^1(\mathbb R^n)\), we test the remainder equation with the bounded, radially nondecreasing weight \(\langle \boldsymbol y \rangle^{2\beta} \langle \alpha \boldsymbol y \rangle^{-2\beta}\) (\(0 < \alpha \le 1\)) before letting \(\alpha \downarrow 0\).
For \(0 \le \beta < (n + 2)/2\) and \(0 < \alpha \le 1\), direct differentiation gives uniformly on \(\mathbb R^n\):
\begin{equation}\label{eq:walphaderiv}
  \begin{gathered}
  \nabla\left(\frac{\langle \boldsymbol y \rangle^{2\beta}}{\langle \alpha \boldsymbol y \rangle^{2\beta}}\right)
  = 2\beta (1 - \alpha^2) \frac{\langle \boldsymbol y \rangle^{2\beta - 2}}{\langle \alpha \boldsymbol y \rangle^{2\beta + 2}} \boldsymbol y, \qquad
  0 \le \boldsymbol y \cdot \nabla\left(\frac{\langle \boldsymbol y \rangle^{2\beta}}{\langle \alpha \boldsymbol y \rangle^{2\beta}}\right)
  \le 2\beta \frac{\langle \boldsymbol y \rangle^{2\beta}}{\langle \alpha \boldsymbol y \rangle^{2\beta}}, \\
  \left|\nabla\left(\frac{\langle \boldsymbol y \rangle^{2\beta}}{\langle \alpha \boldsymbol y \rangle^{2\beta}}\right)\right|
  \le 2\beta \frac{\langle \boldsymbol y \rangle^{2\beta - 1}}{\langle \alpha \boldsymbol y \rangle^{2\beta}}, \qquad
  \left|\Delta\left(\frac{\langle \boldsymbol y \rangle^{2\beta}}{\langle \alpha \boldsymbol y \rangle^{2\beta}}\right)\right|
  \le C_\beta \frac{\langle \boldsymbol y \rangle^{2\beta - 2}}{\langle \alpha \boldsymbol y \rangle^{2\beta}}.
  \end{gathered}
\end{equation}
Moreover, whenever \(|\boldsymbol y - \boldsymbol z| \le 1\) and \(0 < \alpha \le 1\), we have
\[
  \frac{1}{\sqrt 3} \langle \boldsymbol z \rangle \le \langle \boldsymbol y \rangle \le \sqrt 3 \langle \boldsymbol z \rangle,
  \qquad
  \frac{1}{\sqrt 3} \langle \alpha \boldsymbol z \rangle \le \langle \alpha \boldsymbol y \rangle \le \sqrt 3 \langle \alpha \boldsymbol z \rangle.
\]
Cover \(\{|\boldsymbol y| > \ell\}\) (\(\ell \ge 4\)) by the balls \(\{B_{1/4}(\boldsymbol y_j)\}\) from a cubic lattice that intersect \(\{|\boldsymbol y| > \ell\}\), so that \(|\boldsymbol y_j| > \ell - 1/4\) and the doubled balls \(\{B_{1/2}(\boldsymbol y_j)\}\) have uniformly finite overlap.
Multiplying \eqref{eq:unitballL4} by \(\langle \boldsymbol y_j \rangle^{2\beta - 2} \langle \alpha \boldsymbol y_j \rangle^{-2\beta}\) and summing over \(j\) using \(\langle \boldsymbol y_j \rangle^{-4/n} \le C \ell^{-4/n}\) bounds the weighted quartic tail below; factoring \(|\mathbf W|^3\) via Cauchy--Schwarz then bounds the cubic tail as well:
\begin{equation}\label{eq:weightedtails}
  \begin{aligned}
  \int_{|\boldsymbol y| > \ell}
    \frac{ \langle \boldsymbol y \rangle^{2\beta - 2}}{\langle \alpha \boldsymbol y \rangle^{2\beta}} |\mathbf W|^4 \, \mathrm d \boldsymbol y
  &\le \frac{C_\beta}{\ell^{4/n}} \int_{\mathbb R^n}
       \frac{\langle \boldsymbol y \rangle^{2\beta}}{\langle \alpha \boldsymbol y \rangle^{2\beta}}
       \bigl(|\nabla \mathbf W|^2 + |\mathbf W|^2\bigr) \, \mathrm d \boldsymbol y, \\
  \int_{|\boldsymbol y| > \ell}
    \frac{\langle \boldsymbol y \rangle^{2\beta - 1}}{\langle \alpha \boldsymbol y \rangle^{2\beta}} |\mathbf W|^3 \, \mathrm d \boldsymbol y
  &\le \left(\int_{|\boldsymbol y| > \ell}
             \frac{\langle \boldsymbol y \rangle^{2\beta - 2}}{\langle \alpha \boldsymbol y \rangle^{2\beta}} |\mathbf W|^4 \, \mathrm d \boldsymbol y\right)^{1/2}
       \left(\int_{|\boldsymbol y| > \ell}
             \frac{\langle \boldsymbol y \rangle^{2\beta}}{\langle \alpha \boldsymbol y \rangle^{2\beta}} |\mathbf W|^2 \, \mathrm d \boldsymbol y\right)^{1/2} \\
  &\le \frac{C_\beta}{\ell^{2/n}} \int_{\mathbb R^n}
       \frac{\langle \boldsymbol y \rangle^{2\beta}}{\langle \alpha \boldsymbol y \rangle^{2\beta}}
       \bigl(|\nabla \mathbf W|^2 + |\mathbf W|^2\bigr) \, \mathrm d \boldsymbol y.
  \end{aligned}
\end{equation}

We next verify that the pressure weight \(\langle \boldsymbol y \rangle^{2\beta - 2} \langle \alpha \boldsymbol y \rangle^{-2\beta}\) satisfies Muckenhoupt's \(A_2(\mathbb R^n)\) condition \cite{Muckenhoupt1972} uniformly in \(\alpha \in (0, 1]\).

\begin{lemma}\label{lem:weightedriesz}
Let \(n \ge 5\), \(0 \le \beta < (n + 2)/2\), and \(0 < \alpha \le 1\).
For every tensor field \(\mathbf F = (F_{ij}) \in L^2(\mathbb R^n; \mathbb R^{n\times n}) + L^{n/(n-2)}(\mathbb R^n; \mathbb R^{n\times n})\) for which the right-hand integral below is finite, there holds
\begin{equation}\label{eq:weightedriesz}
  \int_{\mathbb R^n} \frac{\langle \boldsymbol y \rangle^{2\beta - 2}}{ \langle \alpha \boldsymbol y \rangle^{2\beta}}
  \left|\sum_{i=1}^n \sum_{j=1}^n \mathcal R_i \mathcal R_j F_{ij}(\boldsymbol y)\right|^2 \, \mathrm d \boldsymbol y
  \le C_\beta \int_{\mathbb R^n} \frac{\langle \boldsymbol y \rangle^{2\beta - 2}}{\langle \alpha \boldsymbol y \rangle^{2\beta}}
      |\mathbf F(\boldsymbol y)|^2 \, \mathrm d \boldsymbol y,
\end{equation}
with a constant \(C_\beta > 0\) independent of \(\alpha\).
\end{lemma}

\begin{proof}
For \(r > 0\), since \(0 \le \beta < (n + 2)/2\) and \(0 < \alpha \le 1\) imply \(0 \le \frac{\alpha^2 r^2}{1 + \alpha^2 r^2} \le \frac{r^2}{1 + r^2} < 1\), direct differentiation gives
\[
  r \frac{\mathrm d}{\mathrm d r} \log\left(\frac{\langle r \rangle^{2\beta - 2}}{\langle \alpha r \rangle^{2\beta}}\right)
  = 2\beta \left(\frac{r^2}{1 + r^2} - \frac{\alpha^2 r^2}{1 + \alpha^2 r^2}\right) - \frac{2 r^2}{1 + r^2}
  \in \bigl[-2, \, (2\beta - 2)_+\bigr],
\]
with \((2\beta - 2)_+ < n\).
Consequently, for all \(r > 0\),
\[
  \begin{aligned}
  \frac{\mathrm d}{\mathrm d r}\left(r^n \frac{\langle r \rangle^{2\beta - 2}}{\langle \alpha r \rangle^{2\beta}}\right)
  &\ge (n - 2) \, r^{n-1} \frac{\langle r \rangle^{2\beta - 2}}{\langle \alpha r \rangle^{2\beta}}, \\
  \frac{\mathrm d}{\mathrm d r}\left(r^n \frac{\langle \alpha r \rangle^{2\beta}}{\langle r \rangle^{2\beta - 2}}\right)
  &\ge \bigl(n - (2\beta - 2)_+\bigr) \, r^{n-1} \frac{\langle \alpha r \rangle^{2\beta}}{\langle r \rangle^{2\beta - 2}}.
  \end{aligned}
\]
Integrating both inequalities over \(r \in (0, \rho)\) (\(\rho > 0\)) yields, on every centered ball \(B_\rho\),
\[
  \left(\frac{1}{|B_\rho|} \int_{B_\rho}
    \frac{\langle \boldsymbol y \rangle^{2\beta - 2}}{\langle \alpha \boldsymbol y \rangle^{2\beta}} \, \mathrm d \boldsymbol y\right)
  \left(\frac{1}{|B_\rho|} \int_{B_\rho}
    \frac{\langle \alpha \boldsymbol y \rangle^{2\beta}}{\langle \boldsymbol y \rangle^{2\beta - 2}} \, \mathrm d \boldsymbol y\right)
  \le \frac{n^2}{(n - 2)\bigl(n - (2\beta - 2)_+\bigr)}.
\]
For any ball \(B_\rho(\boldsymbol y_0)\) (\(\boldsymbol y_0 \in \mathbb R^n\), \(\rho > 0\)), if \(|\boldsymbol y_0| < 2\rho\) then \(B_\rho(\boldsymbol y_0) \subset B_{3\rho}\) with \(|B_{3\rho}| = 3^n |B_\rho|\), so the corresponding product of averages over \(B_\rho(\boldsymbol y_0)\) is bounded by \(9^n\) times the preceding right-hand side.
If \(|\boldsymbol y_0| \ge 2\rho\), then every \(\boldsymbol y, \boldsymbol z \in B_\rho(\boldsymbol y_0)\) satisfies \(|\boldsymbol y_0| - \rho \le |\boldsymbol y|, |\boldsymbol z| \le 3(|\boldsymbol y_0| - \rho)\).
Since \(\max\{2, (2\beta - 2)_+\} \le n\), integrating the logarithmic derivative bound above between \(|\boldsymbol z|\) and \(|\boldsymbol y|\) gives
\[
  \frac{\langle \boldsymbol y \rangle^{2\beta - 2}}{\langle \alpha \boldsymbol y \rangle^{2\beta}}
  \le 3^n \frac{\langle \boldsymbol z \rangle^{2\beta - 2}}{\langle \alpha \boldsymbol z \rangle^{2\beta}},
\]
which bounds the product of averages over \(B_\rho(\boldsymbol y_0)\) by \(3^n\).

Hence the \(A_2(\mathbb R^n)\) constant of \(\langle \boldsymbol y \rangle^{2\beta - 2} \langle \alpha \boldsymbol y \rangle^{-2\beta}\) is bounded uniformly in \(\alpha \in (0, 1]\).
By the weighted \(L^2\) boundedness of Calder\'on--Zygmund singular integrals \cite{Muckenhoupt1972}, \cite[Chapter~V, Section~4.2]{Stein1993}, each Riesz transform \(\mathcal R_i\) is bounded on \(L^2\bigl(\mathbb R^n, \langle \boldsymbol y \rangle^{2\beta - 2} \langle \alpha \boldsymbol y \rangle^{-2\beta} \, \mathrm d \boldsymbol y\bigr)\) with operator norm uniform in \(\alpha \in (0, 1]\).
Since \(\langle \boldsymbol y \rangle^{2\beta - 2} \langle \alpha \boldsymbol y \rangle^{-2\beta} \le \alpha^{-2\beta}\), this gives \eqref{eq:weightedriesz} for all \(\mathbf F \in L^2(\mathbb R^n; \mathbb R^{n\times n})\).
For general \(\mathbf F \in L^2(\mathbb R^n) + L^{n/(n-2)}(\mathbb R^n)\) with finite right-hand side in \eqref{eq:weightedriesz}, positivity of the weight on each ball \(B_k\) gives \(\mathbf 1_{B_k} \mathbf F \in L^2(\mathbb R^n; \mathbb R^{n\times n})\) for every \(k \ge 1\).
As \(k \to \infty\), \(\mathbf 1_{B_k} \mathbf F \to \mathbf F\) strongly both in the weighted \(L^2\) space and in \(L^2(\mathbb R^n) + L^{n/(n-2)}(\mathbb R^n)\), so Calder\'on--Zygmund boundedness of \(\mathcal R_i \mathcal R_j\) and Fatou's lemma along an almost-everywhere convergent subsequence give \eqref{eq:weightedriesz}.
\end{proof}

\begin{proposition}\label{prop:weightedH1}
Under the hypotheses of Proposition~\ref{prop:subquadmorrey}, assume further that \(R \ge 1\) is chosen so that \(\mathbf B = \mathbf u_R(1)\) satisfies \eqref{eq:coercivity}.
Then for every \(0 \le \beta < (n + 2)/2\), there exists a constant \(C_\beta > 0\), depending on \(n\), \(\mathbf u_0\), \(R\), \(\beta\), and \(\|\mathbf W\|_{H^1}\) (via the uniform compact-ball bounds from Proposition~\ref{thm:stationaryreg}, and hence on \(n\), \(\mathbf u_0\), \(R\), and \(\beta\) alone for the profile of Proposition~\ref{prop:limitprofile}), such that
\begin{equation}\label{eq:weightedH1}
  \int_{\mathbb R^n} \langle \boldsymbol y \rangle^{2\beta}
    \bigl(|\nabla \mathbf W(\boldsymbol y)|^2 + |\mathbf W(\boldsymbol y)|^2\bigr) \, \mathrm d \boldsymbol y
  + \int_{\mathbb R^n} \langle \boldsymbol y \rangle^{2\beta - 2}
    |\Pi(\boldsymbol y) - \Pi_B(\boldsymbol y)|^2 \, \mathrm d \boldsymbol y
  \le C_\beta.
\end{equation}
\end{proposition}

\begin{proof}
Let \(0 \le \beta < (n + 2)/2\) and \(0 < \alpha \le 1\).
Since \(\mathbf U, \Pi\) are smooth by hypothesis and \(\mathbf B, \Pi_B, \mathbf S_R(1)\) are smooth by Proposition~\ref{prop:defect}, \(\mathbf W = \mathbf U - \mathbf B\) satisfies the smooth remainder equation
\[
  -\Delta \mathbf W - \frac{1}{2} \boldsymbol y \cdot \nabla \mathbf W - \frac{1}{2} \mathbf W
  + ((\mathbf B + \mathbf W) \cdot \nabla) \mathbf W
  + (\mathbf W \cdot \nabla) \mathbf B + \nabla (\Pi - \Pi_B)
  = \frac{1}{2} \mathbf S_R(1).
\]
At fixed \(\alpha \in (0, 1]\), the damped weight satisfies
\[
  \frac{\langle \boldsymbol y \rangle^\beta}{\langle \alpha \boldsymbol y \rangle^\beta} \le \frac{1}{\alpha^\beta},
  \qquad
  \frac{\langle \boldsymbol y \rangle^{\beta - 1}}{\langle \alpha \boldsymbol y \rangle^\beta} \le \frac{1}{\alpha^\beta}.
\]
Combining \(\mathbf W \in H^1(\mathbb R^n)\) and \(|\mathbf B(\boldsymbol y)| \le C \langle \boldsymbol y \rangle^{-1}\) with \eqref{eq:weightedtails} and Lemma~\ref{lem:weightedriesz} shows that
\[
  \frac{\langle \boldsymbol y \rangle^{2\beta - 1}}{\langle \alpha \boldsymbol y \rangle^{2\beta}} |\mathbf W|^3 \in L^1(\mathbb R^n),
  \qquad
  \frac{\langle \boldsymbol y \rangle^{\beta - 1}}{\langle \alpha \boldsymbol y \rangle^\beta} (\Pi - \Pi_B) \in L^2(\mathbb R^n).
\]

Choose \(\eta \in C_c^\infty(B_2)\) with \(0 \le \eta \le 1\) and \(\eta|_{B_1} = 1\), and for \(\ell \ge 1\) test the smooth remainder equation against the compactly supported vector field \(\eta(\boldsymbol y / \ell) \langle \boldsymbol y \rangle^{2\beta} \langle \alpha \boldsymbol y \rangle^{-2\beta} \mathbf W\).
Integrating by parts using \(\operatorname{div} \mathbf B = \operatorname{div} \mathbf W = 0\) and \(\operatorname{div} \boldsymbol y = n\) gives, for each \(\ell \ge 1\),
\[
  \begin{aligned}
  &-\int_{\mathbb R^n} \eta(\boldsymbol y / \ell)
      \frac{\langle \boldsymbol y \rangle^{2\beta}}{\langle \alpha \boldsymbol y \rangle^{2\beta}}
      \Delta \mathbf W \cdot \mathbf W \, \mathrm d \boldsymbol y \\
  ={}& \int_{\mathbb R^n} \eta(\boldsymbol y / \ell)
         \frac{\langle \boldsymbol y \rangle^{2\beta}}{\langle \alpha \boldsymbol y \rangle^{2\beta}}
         |\nabla \mathbf W|^2 \, \mathrm d \boldsymbol y
       - \frac{1}{2} \int_{\mathbb R^n}
           \Delta\left(\eta(\boldsymbol y / \ell)
                       \frac{\langle \boldsymbol y \rangle^{2\beta}}{\langle \alpha \boldsymbol y \rangle^{2\beta}}\right)
           |\mathbf W|^2 \, \mathrm d \boldsymbol y,
  \end{aligned}
\]
together with
\[
  \begin{aligned}
  &-\int_{\mathbb R^n} \eta(\boldsymbol y / \ell)
      \frac{\langle \boldsymbol y \rangle^{2\beta}}{\langle \alpha \boldsymbol y \rangle^{2\beta}}
      \left(\frac{1}{2} \boldsymbol y \cdot \nabla \mathbf W + \frac{1}{2} \mathbf W\right) \cdot \mathbf W \, \mathrm d \boldsymbol y \\
  ={}& \frac{n - 2}{4} \int_{\mathbb R^n} \eta(\boldsymbol y / \ell)
         \frac{\langle \boldsymbol y \rangle^{2\beta}}{\langle \alpha \boldsymbol y \rangle^{2\beta}}
         |\mathbf W|^2 \, \mathrm d \boldsymbol y \\
     &+ \frac{1}{4} \int_{\mathbb R^n}
          \left[\boldsymbol y \cdot \nabla\left(\eta(\boldsymbol y / \ell)
                                                \frac{\langle \boldsymbol y \rangle^{2\beta}}{\langle \alpha \boldsymbol y \rangle^{2\beta}}\right)\right]
          |\mathbf W|^2 \, \mathrm d \boldsymbol y,
  \end{aligned}
\]
and
\[
  \begin{aligned}
  &\int_{\mathbb R^n} \eta(\boldsymbol y / \ell)
     \frac{\langle \boldsymbol y \rangle^{2\beta}}{\langle \alpha \boldsymbol y \rangle^{2\beta}}
     \bigl[((\mathbf B + \mathbf W) \cdot \nabla) \mathbf W + \nabla(\Pi - \Pi_B)\bigr] \cdot \mathbf W \, \mathrm d \boldsymbol y \\
  ={}& -\int_{\mathbb R^n}
         \nabla\left(\eta(\boldsymbol y / \ell)
                     \frac{\langle \boldsymbol y \rangle^{2\beta}}{\langle \alpha \boldsymbol y \rangle^{2\beta}}\right)
         \cdot \left[\frac{1}{2} |\mathbf W|^2 (\mathbf B + \mathbf W) + (\Pi - \Pi_B) \mathbf W\right] \mathrm d \boldsymbol y.
  \end{aligned}
\]
On the support \(\{\ell \le |\boldsymbol y| \le 2\ell\}\) of \(\nabla \eta(\boldsymbol y / \ell)\), we have \(\langle \boldsymbol y \rangle |\nabla \eta(\boldsymbol y / \ell)| + \langle \boldsymbol y \rangle^2 |\Delta \eta(\boldsymbol y / \ell)| \le C\).
Combining these bounds with \eqref{eq:walphaderiv} and the whole-space integrability at fixed \(\alpha\) established above shows that every integral in which a derivative hits \(\eta(\boldsymbol y / \ell)\) vanishes as \(\ell \to \infty\), while all remaining integrals converge by dominated convergence.
Passing to the limit \(\ell \to \infty\), dropping the nonnegative radial-drift integral via \eqref{eq:walphaderiv}, and applying \eqref{eq:coercivity} together with the derivative bounds in \eqref{eq:walphaderiv} yields
\begin{equation}\label{eq:walphaidentity}
  \begin{aligned}
  &\int_{\mathbb R^n}
     \frac{\langle \boldsymbol y \rangle^{2\beta}}{\langle \alpha \boldsymbol y \rangle^{2\beta}}
     \left(|\nabla \mathbf W|^2 + \frac{n - 2}{8} |\mathbf W|^2\right) \mathrm d \boldsymbol y \\
  \le{}& \frac{1}{2} \int_{\mathbb R^n}
        \frac{\langle \boldsymbol y \rangle^{2\beta}}{\langle \alpha \boldsymbol y \rangle^{2\beta}}
        \mathbf S_R(1) \cdot \mathbf W \, \mathrm d \boldsymbol y 
  + C_\beta \int_{\mathbb R^n}
        \frac{\langle \boldsymbol y \rangle^{2\beta - 2}}{\langle \alpha \boldsymbol y \rangle^{2\beta}}
        \bigl(1 + \langle \boldsymbol y \rangle |\mathbf B + \mathbf W|\bigr) |\mathbf W|^2 \, \mathrm d \boldsymbol y \\
  &+ 2\beta \int_{\mathbb R^n}
        \frac{\langle \boldsymbol y \rangle^{\beta - 1}}{\langle \alpha \boldsymbol y \rangle^\beta} |\Pi - \Pi_B| \,
        \frac{\langle \boldsymbol y \rangle^\beta}{\langle \alpha \boldsymbol y \rangle^\beta} |\mathbf W| \, \mathrm d \boldsymbol y.
  \end{aligned}
\end{equation}

Now let \(\ell \ge 4\) be a splitting radius to be chosen shortly.
For the source term, since \(2\beta - 2n - 2 < -n\), \eqref{eq:defectdecay} and Cauchy--Schwarz give
\[
  \frac{1}{2} \int_{\mathbb R^n}
    \frac{\langle \boldsymbol y \rangle^{2\beta}}{\langle \alpha \boldsymbol y \rangle^{2\beta}}
    \mathbf S_R(1) \cdot \mathbf W \, \mathrm d \boldsymbol y
  \le C_\beta \left(\int_{\mathbb R^n}
        \frac{\langle \boldsymbol y \rangle^{2\beta}}{\langle \alpha \boldsymbol y \rangle^{2\beta}}
        |\mathbf W|^2 \, \mathrm d \boldsymbol y\right)^{1/2}.
\]
In the second integral on the right-hand side of \eqref{eq:walphaidentity}, using \(|\mathbf B(\boldsymbol y)| \le C \langle \boldsymbol y \rangle^{-1}\) gives
\[
  1 + \langle \boldsymbol y \rangle |\mathbf B + \mathbf W|
  \le C + \langle \boldsymbol y \rangle |\mathbf W|.
\]
Splitting \(\mathbb R^n = B_\ell \cup (\mathbb R^n \setminus B_\ell)\), smoothness of \(\mathbf W\) and \(\langle \alpha \boldsymbol y \rangle^{-2\beta} \le 1\) bound the integral over \(B_\ell\) by a constant \(C_{\beta,\ell}\) independent of \(\alpha\); on \(\{|\boldsymbol y| > \ell\}\), using \(\langle \boldsymbol y \rangle^{-2} \le \ell^{-2}\) and \eqref{eq:weightedtails} yields
\[
  \begin{aligned}
  &\int_{\mathbb R^n}
     \frac{\langle \boldsymbol y \rangle^{2\beta - 2}}{\langle \alpha \boldsymbol y \rangle^{2\beta}}
     \bigl(1 + \langle \boldsymbol y \rangle |\mathbf B + \mathbf W|\bigr) |\mathbf W|^2 \, \mathrm d \boldsymbol y \\
  \le{}& C_{\beta,\ell}
       + C_\beta \int_{|\boldsymbol y| > \ell}
           \frac{\langle \boldsymbol y \rangle^{2\beta - 2}}{\langle \alpha \boldsymbol y \rangle^{2\beta}}
           |\mathbf W|^2 \, \mathrm d \boldsymbol y
       + C_\beta \int_{|\boldsymbol y| > \ell}
           \frac{\langle \boldsymbol y \rangle^{2\beta - 1}}{\langle \alpha \boldsymbol y \rangle^{2\beta}}
           |\mathbf W|^3 \, \mathrm d \boldsymbol y \\
  \le{}& C_{\beta,\ell}
       + C_\beta \left(\frac{1}{\ell^2} + \frac{1}{\ell^{2/n}}\right)
         \int_{\mathbb R^n}
           \frac{\langle \boldsymbol y \rangle^{2\beta}}{\langle \alpha \boldsymbol y \rangle^{2\beta}}
           \bigl(|\nabla \mathbf W|^2 + |\mathbf W|^2\bigr) \, \mathrm d \boldsymbol y.
  \end{aligned}
\]
Similarly, in the pressure representation \eqref{eq:mainpressure}, we have
\[
  |\mathbf B \otimes \mathbf W + \mathbf W \otimes \mathbf B + \mathbf W \otimes \mathbf W|^2
  \le C \langle \boldsymbol y \rangle^{-2} |\mathbf W|^2 + C |\mathbf W|^4.
\]
Applying Lemma~\ref{lem:weightedriesz}, splitting \(\mathbb R^n = B_\ell \cup (\mathbb R^n \setminus B_\ell)\), and bounding the quartic exterior tail via \eqref{eq:weightedtails} gives
\begin{equation}\label{eq:pressuresplit}
  \begin{aligned}
  &\int_{\mathbb R^n} \frac{\langle \boldsymbol y \rangle^{2\beta - 2}}{\langle \alpha \boldsymbol y \rangle^{2\beta}} |\Pi - \Pi_B|^2 \, \mathrm d \boldsymbol y \\
  \le{}& C_\beta \int_{\mathbb R^n} \frac{\langle \boldsymbol y \rangle^{2\beta - 2}}{\langle \alpha \boldsymbol y \rangle^{2\beta}}
           \bigl(\langle \boldsymbol y \rangle^{-2} |\mathbf W|^2 + |\mathbf W|^4\bigr) \, \mathrm d \boldsymbol y \\
  \le{}& C_{\beta,\ell}
       + C_\beta \int_{|\boldsymbol y| > \ell}
           \frac{\langle \boldsymbol y \rangle^{2\beta - 4}}{\langle \alpha \boldsymbol y \rangle^{2\beta}}
           |\mathbf W|^2 \, \mathrm d \boldsymbol y
       + C_\beta \int_{|\boldsymbol y| > \ell}
           \frac{\langle \boldsymbol y \rangle^{2\beta - 2}}{\langle \alpha \boldsymbol y \rangle^{2\beta}}
           |\mathbf W|^4 \, \mathrm d \boldsymbol y \\
  \le{}& C_{\beta,\ell}
       + C_\beta \left(\frac{1}{\ell^4} + \frac{1}{\ell^{4/n}}\right)
         \int_{\mathbb R^n}
           \frac{\langle \boldsymbol y \rangle^{2\beta}}{\langle \alpha \boldsymbol y \rangle^{2\beta}}
           \bigl(|\nabla \mathbf W|^2 + |\mathbf W|^2\bigr) \, \mathrm d \boldsymbol y.
  \end{aligned}
\end{equation}

Since \(\frac{n - 2}{8} \ge \frac{3}{8}\) for \(n \ge 5\), combining \eqref{eq:pressuresplit} and the two preceding bounds with \eqref{eq:walphaidentity} via Young's inequality and choosing \(\ell \ge 4\) large enough to absorb all \(O(\ell^{-2/n})\) terms into the left-hand side gives
\[
  \sup_{0 < \alpha \le 1}
  \int_{\mathbb R^n} \frac{\langle \boldsymbol y \rangle^{2\beta}}{\langle \alpha \boldsymbol y \rangle^{2\beta}}
    \bigl(|\nabla \mathbf W|^2 + |\mathbf W|^2\bigr) \, \mathrm d \boldsymbol y
  \le C_\beta.
\]
Substituting this uniform bound back into \eqref{eq:pressuresplit} and sending \(\alpha \downarrow 0\) via monotone convergence establishes \eqref{eq:weightedH1}.
\end{proof}

\subsection{Pointwise decay and proofs of \texorpdfstring{Theorems~\ref*{thm:main} and~\ref*{thm:structure}}{Theorems 1.1 and 1.2}}\label{subsec:bootstrap}

Consider the linear similarity operator
\[
  \mathcal L := -\Delta - \frac{1}{2} \boldsymbol y \cdot \nabla - \frac{1}{2}.
\]
Following the logarithmic-time heat semigroup representation for \(\mathcal L\) \cite{Brandolese2009,CannonePlanchon1996}, for \(\tau \ge 0\) we define
\begin{equation}\label{eq:Lsemigroup}
  (e^{-\tau \mathcal L} \mathbf f)(\boldsymbol y)
  := e^{\tau / 2} \bigl(e^{(e^\tau - 1)\Delta} \mathbf f\bigr)(e^{\tau / 2} \boldsymbol y).
\end{equation}
For any Schwartz field \(\mathbf f\), differentiating in \(\tau\) via the chain rule and \(e^\tau (\Delta \mathbf g)(e^{\tau / 2} \boldsymbol y) = \Delta_{\boldsymbol y} \bigl[\mathbf g(e^{\tau / 2} \boldsymbol y)\bigr]\) gives
\[
  \partial_\tau (e^{-\tau \mathcal L} \mathbf f)(\boldsymbol y)
  = \left(\frac{1}{2} + \frac{1}{2} \boldsymbol y \cdot \nabla_{\boldsymbol y} + \Delta_{\boldsymbol y}\right)
    (e^{-\tau \mathcal L} \mathbf f)(\boldsymbol y)
  = -(\mathcal L e^{-\tau \mathcal L} \mathbf f)(\boldsymbol y)
  = -(e^{-\tau \mathcal L} \mathcal L \mathbf f)(\boldsymbol y),
\]
while substituting \(\boldsymbol z = e^{\tau / 2} \boldsymbol y\) (\(\mathrm d \boldsymbol y = e^{-\frac{n}{2} \tau} \, \mathrm d \boldsymbol z\)) together with \(L^2\) contractivity of the heat semigroup yields
\[
  \|e^{-\tau \mathcal L} \mathbf f\|_2
  = e^{\frac{\tau}{2} - \frac{n}{4} \tau} \|e^{(e^\tau - 1)\Delta} \mathbf f\|_2
  = e^{-\frac{n - 2}{4} \tau} \|e^{(e^\tau - 1)\Delta} \mathbf f\|_2
  \le e^{-\frac{n - 2}{4} \tau} \|\mathbf f\|_2.
\]

\begin{lemma}\label{lem:Linverse}
Let \(n \ge 5\), and let \(\mathcal L = -\Delta - \frac{1}{2} \boldsymbol y \cdot \nabla - \frac{1}{2}\), with semigroup \((e^{-\tau \mathcal L})_{\tau \ge 0}\) defined by \eqref{eq:Lsemigroup}.
Then the integral operators
\[
  \mathcal L^{-1} \mathbf f := \int_0^\infty e^{-\tau \mathcal L} \mathbf f \, \mathrm d \tau,
  \qquad
  \mathcal L^{-1} \mathbb P \operatorname{div} \mathbf F
  := \int_0^\infty e^{-\tau \mathcal L} \mathbb P \operatorname{div} \mathbf F \, \mathrm d \tau
\]
define bounded linear operators \(\mathcal L^{-1} : L^2(\mathbb R^n; \mathbb R^n) \to L^2(\mathbb R^n; \mathbb R^n)\) and \(\mathcal L^{-1} \mathbb P \operatorname{div} : L^2(\mathbb R^n; \mathbb R^{n\times n}) \to L^2(\mathbb R^n; \mathbb R^n)\), and satisfy for all \(\mathbf f \in L^2(\mathbb R^n; \mathbb R^n)\), \(\mathbf F \in L^2(\mathbb R^n; \mathbb R^{n\times n})\), and \(0 < q \le n + 1\) the weighted estimates
\begin{equation}\label{eq:Linverseweighted}
  \begin{aligned}
  \sup_{\boldsymbol y \in \mathbb R^n} \langle \boldsymbol y \rangle^{n+1} |\mathcal L^{-1} \mathbf f(\boldsymbol y)|
  &\le C \sup_{\boldsymbol y \in \mathbb R^n} \langle \boldsymbol y \rangle^{n+1} |\mathbf f(\boldsymbol y)|, \\
  \sup_{\boldsymbol y \in \mathbb R^n} \langle \boldsymbol y \rangle^q |\mathcal L^{-1} \mathbb P \operatorname{div} \mathbf F(\boldsymbol y)|
  &\le C_q \sup_{\boldsymbol y \in \mathbb R^n} \langle \boldsymbol y \rangle^q |\mathbf F(\boldsymbol y)|,
  \end{aligned}
\end{equation}
where whenever the right-hand side is finite, the time-integrated convolution provides a pointwise representative agreeing with the \(L^2\) Bochner integral.
Moreover, if \(\mathbf w \in C^\infty(\mathbb R^n; \mathbb R^n) \cap H^1(\mathbb R^n; \mathbb R^n)\) satisfies \(\mathcal L \mathbf w = \mathbf f - \mathbb P \operatorname{div} \mathbf F\) on \(\mathbb R^n\) for some \(\mathbf f \in L^2(\mathbb R^n; \mathbb R^n)\) and \(\mathbf F \in H^1(\mathbb R^n; \mathbb R^{n\times n})\), then \(\mathbf w = \mathcal L^{-1} \mathbf f - \mathcal L^{-1} \mathbb P \operatorname{div} \mathbf F\).
\end{lemma}

\begin{proof}
By \eqref{eq:Lsemigroup} and \eqref{eq:kerneldef}, we have
\[
  (e^{-\tau \mathcal L} \mathbb P \operatorname{div} \mathbf F)(\boldsymbol y)
  = e^{\tau / 2} (\mathcal K_{e^\tau - 1} * \mathbf F)(e^{\tau / 2} \boldsymbol y).
\]
Changing variables and applying Young's convolution inequality with \(\|\mathcal K_{e^\tau - 1}\|_1 \le C (e^\tau - 1)^{-1/2}\) gives
\[
  \|e^{-\tau \mathcal L} \mathbb P \operatorname{div} \mathbf F\|_2
  = e^{-\frac{n-2}{4} \tau} \|\mathcal K_{e^\tau - 1} * \mathbf F\|_2
  \le C \frac{e^{-\frac{n-2}{4} \tau}}{\sqrt{e^\tau - 1}} \|\mathbf F\|_2,
\]
which is integrable over \(\tau \in (0, \infty)\) and proves \(L^2(\mathbb R^n)\) boundedness.

To prove \eqref{eq:Linverseweighted}, note that for every \(\tau > 0\) and \(\boldsymbol y \in \mathbb R^n\),
\[
  \bigl\langle\bigl(e^{\tau / 2} \boldsymbol y, \, \sqrt{e^\tau - 1}\bigr)\bigr\rangle
  = \bigl(1 + e^\tau |\boldsymbol y|^2 + e^\tau - 1\bigr)^{1/2}
  = e^{\tau / 2} \langle \boldsymbol y \rangle.
\]
For \(0 < q \le n + 1\), evaluating the two bounds in \eqref{eq:kernelalltime} at time \(e^\tau - 1 > 0\) and spatial point \(e^{\tau / 2} \boldsymbol y\) (with \(q = n + 1\) in the second bound) and using \(\langle \boldsymbol y \rangle \ge 1\) together with \(\max\{q, \, n + 1/2\} \ge n + 1/2\) gives
\[
  \begin{aligned}
  &\langle \boldsymbol y \rangle^q
    \bigl|(e^{-\tau \mathcal L} \mathbb P \operatorname{div} \mathbf F)(\boldsymbol y)\bigr| \\
  \le{}& C_q \langle \boldsymbol y \rangle^q e^{\tau / 2}
       \left(\frac{1}{\sqrt{e^\tau - 1} \, (e^{\tau / 2} \langle \boldsymbol y \rangle)^q}
             + \frac{1}{(e^{\tau / 2} \langle \boldsymbol y \rangle)^{\max\{q, \, n + 1/2\}}}\right)
       \sup_{\boldsymbol z \in \mathbb R^n} \langle \boldsymbol z \rangle^q |\mathbf F(\boldsymbol z)| \\
  \le{}& C_q \left(\frac{e^{\frac{1 - q}{2} \tau}}{\sqrt{e^\tau - 1}}
                 + e^{-\frac{2n - 1}{4} \tau}\right)
       \sup_{\boldsymbol z \in \mathbb R^n} \langle \boldsymbol z \rangle^q |\mathbf F(\boldsymbol z)|,
  \end{aligned}
\]
and similarly
\[
  \begin{aligned}
  \langle \boldsymbol y \rangle^{n+1}
    \bigl|(e^{-\tau \mathcal L} \mathbf f)(\boldsymbol y)\bigr|
  &\le C \langle \boldsymbol y \rangle^{n+1} e^{\tau / 2}
       \frac{1 + \sqrt{e^\tau - 1}}{(e^{\tau / 2} \langle \boldsymbol y \rangle)^{n+1}}
       \sup_{\boldsymbol z \in \mathbb R^n} \langle \boldsymbol z \rangle^{n+1} |\mathbf f(\boldsymbol z)| \\
  &\le C e^{-\frac{n - 1}{2} \tau}
       \sup_{\boldsymbol z \in \mathbb R^n} \langle \boldsymbol z \rangle^{n+1} |\mathbf f(\boldsymbol z)|.
  \end{aligned}
\]
Because \(q > 0\) and \(n \ge 5\), the factor \(e^{\frac{1 - q}{2} \tau} (e^\tau - 1)^{-1/2}\) behaves as \(\tau^{-1/2}\) as \(\tau \downarrow 0\) and decays as \(e^{-\frac{q}{2} \tau}\) as \(\tau \to \infty\); hence both time factors belong to \(L^1(0, \infty)\), and integrating over \(\tau \in (0, \infty)\) establishes \eqref{eq:Linverseweighted}.

Finally, suppose \(\mathbf w \in C^\infty(\mathbb R^n; \mathbb R^n) \cap H^1(\mathbb R^n; \mathbb R^n)\) satisfies \(\mathcal L \mathbf w = \mathbf f - \mathbb P \operatorname{div} \mathbf F\) with \(\mathbf f \in L^2(\mathbb R^n; \mathbb R^n)\) and \(\mathbf F \in H^1(\mathbb R^n; \mathbb R^{n\times n})\), so that \(\mathcal L \mathbf w \in L^2(\mathbb R^n; \mathbb R^n)\).
Fix \(\eta \in C_c^\infty(B_2)\) with \(\eta|_{B_1} = 1\) and set \(\eta_\ell(\boldsymbol y) := \eta(\boldsymbol y / \ell)\) for \(\ell \ge 1\).
Then \(\eta_\ell \mathbf w \in C_c^\infty(\mathbb R^n; \mathbb R^n)\) satisfies
\[
  \mathcal L(\eta_\ell \mathbf w)
  = \eta_\ell \mathcal L \mathbf w
    - 2 (\nabla \eta_\ell \cdot \nabla) \mathbf w
    - (\Delta \eta_\ell) \mathbf w
    - \frac{1}{2} (\boldsymbol y \cdot \nabla \eta_\ell) \mathbf w,
\]
which converges to \(\mathcal L \mathbf w\) in \(L^2(\mathbb R^n)\) as \(\ell \to \infty\) since \(|\nabla \eta_\ell| \le C \ell^{-1}\), \(|\Delta \eta_\ell| \le C \ell^{-2}\), and \(|\boldsymbol y \cdot \nabla \eta_\ell| \le C \mathbf 1_{\{\ell \le |\boldsymbol y| \le 2\ell\}}\).
Integrating \(\partial_\tau \bigl(e^{-\tau \mathcal L} (\eta_\ell \mathbf w)\bigr) = -e^{-\tau \mathcal L} \mathcal L(\eta_\ell \mathbf w)\) over \([0, T]\) gives
\[
  \eta_\ell \mathbf w - e^{-T \mathcal L}(\eta_\ell \mathbf w)
  = \int_0^T e^{-\tau \mathcal L} \mathcal L(\eta_\ell \mathbf w) \, \mathrm d \tau.
\]
Sending \(\ell \to \infty\) and then \(T \to \infty\) using \(\|e^{-\tau \mathcal L}\|_{L^2 \to L^2} \le e^{-\frac{n-2}{4} \tau}\) yields \(\mathbf w = \mathcal L^{-1} \mathbf f - \mathcal L^{-1} \mathbb P \operatorname{div} \mathbf F\).
\end{proof}

\begin{proposition}\label{prop:sharpdecay}
Under the hypotheses of Proposition~\ref{prop:weightedH1}, there exists a constant \(C > 0\) (with the same dependence as in Proposition~\ref{prop:weightedH1}) such that for all \(\boldsymbol y \in \mathbb R^n\) and every \(q \in (n / 2, \infty)\),
\[
  |\mathbf W(\boldsymbol y)| \le \frac{C}{\langle \boldsymbol y \rangle^{n+1}},
  \qquad
  |\mathbf U(\boldsymbol y)| \le \frac{C}{\langle \boldsymbol y \rangle},
  \qquad
  \Pi \in L^q(\mathbb R^n),
\]
and moreover \(|\mathbf U(\boldsymbol y) - \mathbf u_0(\boldsymbol y)| \le C |\boldsymbol y|^{-3}\) for all \(|\boldsymbol y| \ge 1\).
\end{proposition}

\begin{proof}
\noindent\emph{Step 1.}
Fix \(\beta := (2n + 3) / 4 \in (n / 2, (n + 2) / 2)\) in Proposition~\ref{prop:weightedH1}, so that \(2\beta - n = 3 / 2 > 0\).
For any \(\boldsymbol y_0 \in \mathbb R^n\) with \(r := |\boldsymbol y_0| \ge 4\), rescale the \((n + 1)\)-dimensional stationary solution \eqref{eq:lift} onto \(B_2^{n+1} \subset \mathbb R^{n+1}\):
\[
  \mathbf v_{\boldsymbol y_0}(\boldsymbol \xi', \xi_{n+1})
  := \frac{1}{r} \mathbf V\left(\boldsymbol y_0 + \frac{\boldsymbol \xi'}{r}, \, \frac{\xi_{n+1}}{r}\right),
  \qquad
  \pi_{\boldsymbol y_0}(\boldsymbol \xi', \xi_{n+1})
  := \frac{1}{r^2} \widetilde \Pi\left(\boldsymbol y_0 + \frac{\boldsymbol \xi'}{r}, \, \frac{\xi_{n+1}}{r}\right).
\]
Then \((\mathbf v_{\boldsymbol y_0}, \pi_{\boldsymbol y_0})\) solves \eqref{eq:steadyNS} on \(B_2^{n+1}\).
Since \(r \ge 4 \ge 1\) and \(\sup_{\mathbb R^n} \Phi_B \ge 0\), \eqref{eq:lifthead} gives
\[
  \pi_{\boldsymbol y_0} + \frac{1}{2} |\mathbf v_{\boldsymbol y_0}|^2
  = \frac{1}{r^2} \Phi\left(\boldsymbol y_0 + \frac{\boldsymbol \xi'}{r}\right)
  \le \frac{1}{r^2} \sup_{\mathbb R^n} \Phi_B
  \le \sup_{\mathbb R^n} \Phi_B
  \qquad \text{on } B_2^{n+1}.
\]

Define the rescaled remainder
\[
  \mathbf w_{\boldsymbol y_0}(\boldsymbol \xi')
  := \frac{1}{r} \mathbf W\left(\boldsymbol y_0 + \frac{\boldsymbol \xi'}{r}\right),
  \qquad \boldsymbol \xi' \in B_2.
\]
Since \(|\boldsymbol y_0 + \boldsymbol \xi' / r| \ge \frac{1}{2} r\) on \(B_2\), changing variables in \eqref{eq:weightedH1} yields
\begin{equation}\label{eq:fL2}
  \begin{aligned}
  \|\mathbf w_{\boldsymbol y_0}\|_{L^2(B_2)}^2
  &= r^{n-2} \int_{B_{2/r}(\boldsymbol y_0)} |\mathbf W(\boldsymbol y)|^2 \, \mathrm d \boldsymbol y
   \le \frac{C}{r^{2\beta + 2 - n}} = \frac{C}{r^{7/2}}, \\
  \|\nabla_{\boldsymbol \xi'} \mathbf w_{\boldsymbol y_0}\|_{L^2(B_2)}^2
  &= r^{n-4} \int_{B_{2/r}(\boldsymbol y_0)} |\nabla \mathbf W(\boldsymbol y)|^2 \, \mathrm d \boldsymbol y
   \le \frac{C}{r^{2\beta + 4 - n}} = \frac{C}{r^{11/2}},
  \end{aligned}
\end{equation}
while Cauchy--Schwarz and the weighted pressure bound in \eqref{eq:weightedH1} give
\[
  \begin{aligned}
  \frac{1}{r^2} \int_{B_2} \left|(\Pi - \Pi_B)\left(\boldsymbol y_0 + \frac{\boldsymbol \xi'}{r}\right)\right| \mathrm d \boldsymbol \xi'
  &\le C \left(r^{n-4} \int_{B_{2/r}(\boldsymbol y_0)}
         |\Pi(\boldsymbol y) - \Pi_B(\boldsymbol y)|^2 \, \mathrm d \boldsymbol y\right)^{1/2} \\
  &\le \frac{C}{r^{(2\beta + 2 - n)/2}} = \frac{C}{r^{7/4}}.
  \end{aligned}
\]
By \eqref{eq:lift}, on \(B_2^{n+1}\) (where \(|\boldsymbol y_0 + \boldsymbol \xi' / r| \le \frac{3}{2} r\) and \(|\xi_{n+1}| \le 2\)) the rescaled pair decomposes as
\[
  \begin{aligned}
  \mathbf v_{\boldsymbol y_0}(\boldsymbol \xi', \xi_{n+1})
  &= \bigl(\mathbf w_{\boldsymbol y_0}(\boldsymbol \xi'), \, 0\bigr)
     + \left(\frac{1}{r} \mathbf B\left(\boldsymbol y_0 + \frac{\boldsymbol \xi'}{r}\right)
             - \frac{\boldsymbol y_0}{2 r} - \frac{\boldsymbol \xi'}{2 r^2}, \,
             \frac{n \xi_{n+1}}{2 r^2}\right), \\
  \pi_{\boldsymbol y_0}(\boldsymbol \xi', \xi_{n+1})
  &= \frac{1}{r^2} (\Pi - \Pi_B)\left(\boldsymbol y_0 + \frac{\boldsymbol \xi'}{r}\right)
     + \frac{1}{r^2} \Pi_B\left(\boldsymbol y_0 + \frac{\boldsymbol \xi'}{r}\right)
     - \frac{|\boldsymbol y_0 + \boldsymbol \xi' / r|^2}{8 r^2}
     - \frac{n^2 \xi_{n+1}^2}{8 r^4}.
  \end{aligned}
\]
Since \(\mathbf B, \nabla \mathbf B, \Pi_B \in L^\infty(\mathbb R^n)\) and \(r \ge 4\), the background and straining terms in \(\mathbf v_{\boldsymbol y_0}\) and \(\pi_{\boldsymbol y_0}\) are bounded in \(C^1(B_2^{n+1})\) and \(L^\infty(B_2^{n+1})\), respectively, uniformly in \(\boldsymbol y_0\), so combining them with \eqref{eq:fL2} and the pressure bound above implies
\[
  \sup_{|\boldsymbol y_0| \ge 4}
  \Bigl(\|\mathbf v_{\boldsymbol y_0}\|_{H^1(B_2^{n+1})}
        + \|\pi_{\boldsymbol y_0}\|_{L^1(B_2^{n+1})}\Bigr)
  \le C.
\]

Applying \eqref{eq:quantitativeCk} in Proposition~\ref{thm:stationaryreg} with \(m = n + 1\) and \(k = n\) bounds \(\|\mathbf v_{\boldsymbol y_0}\|_{C^n(B_1^{n+1})} \le C\) uniformly for all \(|\boldsymbol y_0| \ge 4\).
Introduce the rescaled profile
\[
  \mathbf u_{\boldsymbol y_0}(\boldsymbol \xi')
  := \frac{1}{r} \mathbf U\left(\boldsymbol y_0 + \frac{\boldsymbol \xi'}{r}\right)
  = \mathbf w_{\boldsymbol y_0}(\boldsymbol \xi')
    + \frac{1}{r} \mathbf B\left(\boldsymbol y_0 + \frac{\boldsymbol \xi'}{r}\right),
  \qquad \boldsymbol \xi' \in B_1.
\]
Since \(|\mathbf B(\boldsymbol y)| \le C \langle \boldsymbol y \rangle^{-1}\) gives \(\bigl|r^{-1} \mathbf B(\boldsymbol y_0 + \boldsymbol \xi' / r)\bigr| \le C r^{-2}\) on \(B_1\) for \(r \ge 4\), \eqref{eq:fL2} implies \(\|\mathbf u_{\boldsymbol y_0}\|_{L^2(B_1)} \le C r^{-7/4}\).
Writing \(\mathbf v_{\boldsymbol y_0}'\) for the first \(n\) components of \(\mathbf v_{\boldsymbol y_0}\), since the drift term is affine in \(\boldsymbol \xi'\), restricting to \(\xi_{n+1} = 0\) gives
\[
  \nabla_{\boldsymbol \xi'}^n \mathbf u_{\boldsymbol y_0}(\boldsymbol \xi')
  = \nabla_{\boldsymbol \xi'}^n \mathbf v_{\boldsymbol y_0}'(\boldsymbol \xi', 0),
  \qquad
  \sup_{|\boldsymbol y_0| \ge 4} \|\nabla^n \mathbf u_{\boldsymbol y_0}\|_{L^\infty(B_1)} \le C.
\]
Interpolating \(L^\infty(B_1)\) between \(L^2(B_1)\) and \(W^{n,\infty}(B_1)\) on \(B_1 \subset \mathbb R^n\) via the Gagliardo--Nirenberg inequality \cite{Nirenberg1959} (with dimension-independent exponent \(\theta = \frac{n / 2}{n + n / 2} = \frac{1}{3}\)) yields
\[
  \frac{1}{r} |\mathbf U(\boldsymbol y_0)|
  = |\mathbf u_{\boldsymbol y_0}(\boldsymbol 0)|
  \le C \|\mathbf u_{\boldsymbol y_0}\|_{L^2(B_1)}^{2/3}
      \bigl(\|\nabla^n \mathbf u_{\boldsymbol y_0}\|_{L^\infty(B_1)} + \|\mathbf u_{\boldsymbol y_0}\|_{L^2(B_1)}\bigr)^{1/3}
  \le \frac{C}{r^{7/6}}.
\]
Hence \(|\mathbf U(\boldsymbol y_0)| \le C r^{-1/6}\) and \(|\mathbf W(\boldsymbol y_0)| \le |\mathbf U(\boldsymbol y_0)| + |\mathbf B(\boldsymbol y_0)| \le C r^{-1/6}\).
Combining this exterior bound with boundedness of \(\mathbf W\) and \(\mathbf U\) on \(\overline{B_4}\) yields
\begin{equation}\label{eq:seeddecay}
  |\mathbf W(\boldsymbol y)| + |\mathbf U(\boldsymbol y)|
  \le \frac{C}{\langle \boldsymbol y \rangle^{1/6}},
  \qquad \boldsymbol y \in \mathbb R^n.
\end{equation}

\smallskip
\noindent\emph{Step 2.}
Since \(\mathbf W \in C^\infty \cap H^1 \cap L^\infty(\mathbb R^n)\) by \eqref{eq:seeddecay} and \(\mathbf B \in W^{1,\infty}(\mathbb R^n)\) by \eqref{eq:bgbounds}, the tensor field \(\mathbf B \otimes \mathbf W + \mathbf W \otimes \mathbf B + \mathbf W \otimes \mathbf W\) belongs to \(H^1(\mathbb R^n; \mathbb R^{n\times n})\).
Moreover, \eqref{eq:mainpressure} identifies
\[
  \nabla(\Pi - \Pi_B)
  = -(I - \mathbb P) \operatorname{div}
      (\mathbf B \otimes \mathbf W + \mathbf W \otimes \mathbf B + \mathbf W \otimes \mathbf W),
\]
so the remainder equation reads
\[
  \mathcal L \mathbf W
  = \frac{1}{2} \mathbf S_R(1)
    - \mathbb P \operatorname{div}
      (\mathbf B \otimes \mathbf W + \mathbf W \otimes \mathbf B + \mathbf W \otimes \mathbf W),
\]
with \(\mathbf S_R(1) \in L^2(\mathbb R^n; \mathbb R^n)\).
By Lemma~\ref{lem:Linverse}, \(\mathbf W\) satisfies
\begin{equation}\label{eq:Wduhamel}
  \mathbf W
  = \frac{1}{2} \mathcal L^{-1} \mathbf S_R(1)
    - \mathcal L^{-1} \mathbb P \operatorname{div}
      (\mathbf B \otimes \mathbf W + \mathbf W \otimes \mathbf B + \mathbf W \otimes \mathbf W).
\end{equation}

By \eqref{eq:defectdecay}, \eqref{eq:bgbounds}, and \eqref{eq:seeddecay}, we have
\[
  \sup_{\boldsymbol y \in \mathbb R^n} \langle \boldsymbol y \rangle^{n+1} |\mathbf S_R(1, \boldsymbol y)| < \infty,
  \qquad
  \sup_{\boldsymbol y \in \mathbb R^n} \langle \boldsymbol y \rangle |\mathbf B(\boldsymbol y)| < \infty,
  \qquad
  \sup_{\boldsymbol y \in \mathbb R^n} \langle \boldsymbol y \rangle^{q_0} |\mathbf W(\boldsymbol y)| < \infty
\]
with \(q_0 := 1/6\).
Define recursively \(q_{j+1} := \min\{n + 1, \, q_j + 1, \, 2q_j\}\).
Whenever \(|\mathbf W(\boldsymbol y)| \le C_j \langle \boldsymbol y \rangle^{-q_j}\) on \(\mathbb R^n\), combining \(|\mathbf B(\boldsymbol y)| \le C \langle \boldsymbol y \rangle^{-1}\) with the quadratic decay of \(\mathbf W \otimes \mathbf W\) gives
\[
  \sup_{\boldsymbol y \in \mathbb R^n}
  \langle \boldsymbol y \rangle^{q_{j+1}}
  \bigl|\bigl(\mathbf B \otimes \mathbf W + \mathbf W \otimes \mathbf B + \mathbf W \otimes \mathbf W\bigr)(\boldsymbol y)\bigr|
  < \infty,
\]
and applying \eqref{eq:Linverseweighted} to \eqref{eq:Wduhamel} yields \(|\mathbf W(\boldsymbol y)| \le C_{j+1} \langle \boldsymbol y \rangle^{-q_{j+1}}\) on \(\mathbb R^n\).
After \(n + 3\) steps (\(1/6 \to 1/3 \to 2/3 \to 4/3 \to 7/3 \to \cdots \to n + 1\)), we obtain \(q_{n+3} = n + 1\), whence \(|\mathbf W(\boldsymbol y)| \le C \langle \boldsymbol y \rangle^{-n-1}\) and \(|\mathbf U(\boldsymbol y)| \le C \langle \boldsymbol y \rangle^{-1}\) on \(\mathbb R^n\).
Since \(|\mathbf U \otimes \mathbf U(\boldsymbol y)| \le C \langle \boldsymbol y \rangle^{-2} \in L^q(\mathbb R^n)\) for every \(q \in (n / 2, \infty)\), Calder\'on--Zygmund boundedness of \(\mathcal R_i \mathcal R_j\) in \eqref{eq:mainpressure} gives \(\Pi \in L^q(\mathbb R^n)\), establishing \eqref{eq:maindecay} and \eqref{eq:remdecay}.

Finally, for \(|\boldsymbol y| > 2R\) we have \(\mathbf u_{0,R}(\boldsymbol y) = \mathbf u_0(\boldsymbol y)\), so \eqref{eq:bgdiff3} and \eqref{eq:remdecay} give
\[
  |\mathbf U(\boldsymbol y) - \mathbf u_0(\boldsymbol y)|
  \le |\mathbf B(\boldsymbol y) - \mathbf u_{0,R}(\boldsymbol y)| + |\mathbf W(\boldsymbol y)|
  \le \frac{C}{|\boldsymbol y|^3},
\]
while on \(1 \le |\boldsymbol y| \le 2R\) smoothness of \(\mathbf U\) and \(\mathbf u_0\) gives \(|\boldsymbol y|^3 |\mathbf U(\boldsymbol y) - \mathbf u_0(\boldsymbol y)| \le C\).
\end{proof}

\begin{proof}[Proof of Theorems~\ref{thm:main} and~\ref{thm:structure}]
By Lemmas~\ref{lem:cutoff} and~\ref{lem:evolution} and Proposition~\ref{prop:defect}, there exists \(R_0 \ge 1\) with \(R_0^{-2} \le T_0\) such that for every \(R \ge R_0\), the smooth divergence-free datum \(\mathbf u_{0,R}\) generates the unique bounded mild solution \(\mathbf u_R\) on \([0, 1] \times \mathbb R^n\), whose scaling defect satisfies \eqref{eq:defectdecay} and whose time-one slice \(\mathbf B = \mathbf u_R(1)\) satisfies \eqref{eq:coercivity}.
Propositions~\ref{prop:limitprofile} and~\ref{prop:smoothprofile} then yield a smooth solution \((\mathbf U, \Pi)\) of \eqref{eq:profile} satisfying Theorem~\ref{thm:structure}\textup{(i)} together with the classical solvability for \(t > 0\) and strong \(H^1_{\mathrm{loc}}\) initial trace \eqref{eq:initialtrace} in Theorem~\ref{thm:main}\textup{(ii)}, while Proposition~\ref{prop:sharpdecay} establishes Theorem~\ref{thm:main}\textup{(i)} and Theorem~\ref{thm:structure}\textup{(ii)}.

It remains to verify the distributional formulation of \eqref{eq:NS} and the local energy equality up to \(t = 0\) in Theorem~\ref{thm:main}\textup{(ii)}.
Let \(\boldsymbol \varphi \in C_c^\infty([0, \infty) \times \mathbb R^n; \mathbb R^n)\) and \(\phi \in C_c^\infty([0, \infty) \times \mathbb R^n)\) be supported in \([0, T) \times B_\ell\).
By \eqref{eq:similarity} and \eqref{eq:maindecay}, since \(n \ge 5\) and \(3n / 4 \in (n / 2, \infty)\), the velocity satisfies uniformly for \(t > 0\)
\[
  |\mathbf u(t, \boldsymbol x)|
  \le \frac{C}{\sqrt t + |\boldsymbol x|}
  \le \frac{C}{|\boldsymbol x|} \in L^4(B_\ell),
\]
while the pressure satisfies
\[
  \|p(t)\|_{L^{3n/4}(\mathbb R^n)}
  = \frac{\|\Pi\|_{3n/4}}{t^{1/3}} \in L^2(0, T).
\]
In particular, \(\mathbf u \in L^\infty(0, T; L^4(B_\ell))\) implies \(|\mathbf u| + |\mathbf u|^3 \in L^1((0, T) \times B_\ell)\).
Since \(\frac{4}{3n} + \frac{1}{4} < 1\) for \(n \ge 5\), space-time H\"older's inequality pairs \(p \in L^2(0, T; L^{3n/4}(B_\ell))\) with \(\mathbf u \in L^\infty(0, T; L^4(B_\ell))\) to give \(|p| |\mathbf u| \in L^1((0, T) \times B_\ell)\).
Moreover, \(\sup_{0 < t \le T} \|\mathbf u(t)\|_{H^1(B_\ell)} < \infty\) from \eqref{eq:initialtrace} gives \(|\nabla \mathbf u|^2 \in L^1((0, T) \times B_\ell)\).
For every \(t_0 \in (0, T)\), testing the smooth system \eqref{eq:NS} against \(\boldsymbol \varphi\) and \(\mathbf u \phi\) on \([t_0, T] \times \mathbb R^n\) and integrating by parts yields
\[
  \begin{aligned}
  -\int_{\mathbb R^n} \mathbf u(t_0, \boldsymbol x) \cdot \boldsymbol \varphi(t_0, \boldsymbol x) \, \mathrm d \boldsymbol x
  &= \int_{t_0}^T \! \int_{\mathbb R^n}
       \Bigl[\mathbf u \cdot (\partial_t \boldsymbol \varphi + \Delta \boldsymbol \varphi)
             + (\mathbf u \otimes \mathbf u) : \nabla \boldsymbol \varphi
             + p \operatorname{div} \boldsymbol \varphi\Bigr] \, \mathrm d \boldsymbol x \, \mathrm d t, \\
  \int_{t_0}^T \! \int_{\mathbb R^n} |\nabla \mathbf u|^2 \phi \, \mathrm d \boldsymbol x \, \mathrm d t
  &= \frac{1}{2} \int_{\mathbb R^n} |\mathbf u(t_0, \boldsymbol x)|^2 \phi(t_0, \boldsymbol x) \, \mathrm d \boldsymbol x \\
  &\quad + \int_{t_0}^T \! \int_{\mathbb R^n}
             \left[\frac{1}{2} |\mathbf u|^2 (\partial_t \phi + \Delta \phi)
                   + \left(\frac{1}{2} |\mathbf u|^2 + p\right) \mathbf u \cdot \nabla \phi\right]
             \, \mathrm d \boldsymbol x \, \mathrm d t.
  \end{aligned}
\]
Sending \(t_0 \downarrow 0\) via the strong \(L^2(B_\ell)\) convergence in \eqref{eq:initialtrace} on the boundary slices and dominated convergence on the space-time integrals establishes both the distributional formulation of \eqref{eq:NS} and the local energy equality on \([0, \infty) \times \mathbb R^n\).
Finally, for any \(\boldsymbol x_0 \in \mathbb R^n\) and \(T > 0\), changing variables \(\boldsymbol y = \boldsymbol x / \sqrt t\) gives
\[
  \int_{B_1(\boldsymbol x_0)} \Bigl|\nabla_{\boldsymbol x} \bigl[t^{-1/2} \mathbf W(\boldsymbol x / \sqrt t)\bigr]\Bigr|^2 \, \mathrm d \boldsymbol x
  = t^{(n-4)/2} \int_{B_{1/\sqrt t}(\boldsymbol x_0 / \sqrt t)} |\nabla \mathbf W(\boldsymbol y)|^2 \, \mathrm d \boldsymbol y.
\]
Combining this identity with \(|\mathbf u(t, \boldsymbol x)| \le C / |\boldsymbol x|\) and \(|\nabla \mathbf B(\boldsymbol y)| \le C \langle \boldsymbol y \rangle^{-2}\) yields
\[
  \sup_{\boldsymbol x_0 \in \mathbb R^n}
  \left(\sup_{0 < t \le T} \int_{B_1(\boldsymbol x_0)} |\mathbf u(t, \boldsymbol x)|^2 \, \mathrm d \boldsymbol x
        + \int_0^T \! \int_{B_1(\boldsymbol x_0)} |\nabla \mathbf u(t, \boldsymbol x)|^2 \, \mathrm d \boldsymbol x \, \mathrm d t\right)
  < \infty,
\]
and the same bounds vanish as \(|\boldsymbol x_0| \to \infty\) (since \(n \ge 5\) and \(\nabla \mathbf W \in L^2(\mathbb R^n)\)).
Hence \((\mathbf u, p)\) is a local Leray solution attaining \(\mathbf u_0\).
\end{proof}

\appendix

\section{Classical kernel and mild-solution estimates}\label{app:evolution}

This appendix collects classical kernel estimates and local mild-solution theory used in the construction of the parabolic background and the semigroup decay bootstrap.
Pointwise and \(L^1\) bounds for the differentiated Oseen kernel \(\mathcal K_\tau\) representing \(e^{\tau\Delta} \mathbb P \operatorname{div}\) are standard \cite[Section~3, Proposition~1]{Brandolese2009}, \cite[Lemma~2.1]{BrandoleseVigneron2007}, as is the short-time mild-solution theory \cite{Giga1986,Kato1984}.
In Lemma~\ref{lem:kernel} below, \eqref{eq:kernelalltime} gives the all-time weighted convolution estimate applied in Section~\ref{subsec:bootstrap} at heat time \(e^\tau - 1\), whereas \eqref{eq:kernelweighted} records its short-time corollary used in Sections~\ref{sec:evolution} and~\ref{subsec:regularity}.

\begin{lemma}\label{lem:kernel}
Let \(n \ge 5\), and let \(G_\tau(\boldsymbol y) = (4\pi\tau)^{-n/2} e^{-|\boldsymbol y|^2 / (4\tau)}\) (\(\tau > 0\)) denote the Gaussian heat kernel on \(\mathbb R^n\).
For every Schwartz tensor field \(\mathbf F = (F_{kj})_{1 \le k,j \le n}\) and every \(\tau > 0\),
\begin{equation}\label{eq:kerneldef}
  \begin{gathered}
  (e^{\tau\Delta} \mathbb P \operatorname{div} \mathbf F)_i
  = (\mathcal K_\tau * \mathbf F)_i
  := \sum_{k=1}^n \sum_{j=1}^n \mathcal K_{ikj,\tau} * F_{kj}, \\
  \mathcal K_{ikj,\tau}
  = \delta_{ik} \partial_j G_\tau + \int_\tau^\infty \partial_i \partial_k \partial_j G_s \, \mathrm d s,
  \end{gathered}
\end{equation}
where the differentiated Oseen kernel \(\mathcal K_\tau = (\mathcal K_{ikj,\tau})\) satisfies
\begin{equation}\label{eq:kernel}
  |\mathcal K_\tau(\boldsymbol y)|
  \le \frac{C}{(|\boldsymbol y| + \sqrt \tau)^{n+1}},
  \qquad
  \|\mathcal K_\tau\|_1 \le \frac{C}{\sqrt \tau},
  \qquad \tau > 0.
\end{equation}
Furthermore, for all \(0 < q \le n + 1\), \(\tau > 0\), \(\boldsymbol y \in \mathbb R^n\), bounded measurable tensor fields \(\mathbf F\), and bounded measurable functions \(h\),
\begin{equation}\label{eq:kernelalltime}
  \begin{aligned}
  \bigl|(\mathcal K_\tau * \mathbf F)(\boldsymbol y)\bigr|
  &\le C_q \left(\frac{1}{\sqrt \tau \, \langle (\boldsymbol y, \sqrt \tau) \rangle^q}
        + \frac{1}{\langle (\boldsymbol y, \sqrt \tau) \rangle^{\max\{q, \, n + 1/2\}}}\right)
       \sup_{\boldsymbol z \in \mathbb R^n} \langle \boldsymbol z \rangle^q |\mathbf F(\boldsymbol z)|, \\
  \bigl|(e^{\tau\Delta} h)(\boldsymbol y)\bigr|
  &\le C_q \left(\frac{1}{\langle (\boldsymbol y, \sqrt \tau) \rangle^q}
        + \frac{\sqrt \tau}{\langle (\boldsymbol y, \sqrt \tau) \rangle^{\max\{q, \, n + 1/2\}}}\right)
       \sup_{\boldsymbol z \in \mathbb R^n} \langle \boldsymbol z \rangle^q |h(\boldsymbol z)|.
  \end{aligned}
\end{equation}
In particular, for \(0 < \tau \le 1\),
\begin{equation}\label{eq:kernelweighted}
  \begin{aligned}
  \sup_{\boldsymbol y \in \mathbb R^n} \langle \boldsymbol y \rangle^q
    \bigl|(\mathcal K_\tau * \mathbf F)(\boldsymbol y)\bigr|
  &\le \frac{C_q}{\sqrt \tau}
       \sup_{\boldsymbol y \in \mathbb R^n} \langle \boldsymbol y \rangle^q |\mathbf F(\boldsymbol y)|, \\
  \sup_{\boldsymbol y \in \mathbb R^n} \langle \boldsymbol y \rangle^q
    \bigl|(e^{\tau\Delta} h)(\boldsymbol y)\bigr|
  &\le C_q \sup_{\boldsymbol y \in \mathbb R^n} \langle \boldsymbol y \rangle^q |h(\boldsymbol y)|.
  \end{aligned}
\end{equation}
\end{lemma}

\begin{proof}
Combining \(-\Delta \int_\tau^\infty G_s \, \mathrm d s = G_\tau\) in \(\mathscr{D}'(\mathbb R^n)\) with \(\mathbb P_{ik} = \delta_{ik} + \partial_i \partial_k (-\Delta)^{-1}\) yields \eqref{eq:kerneldef} (cf.\ \cite[Section~3, Proposition~1]{Brandolese2009}).
From the Gaussian derivative bound \(|\nabla^3 G_s(\boldsymbol y)| \le C s^{-(n+3)/2} e^{-|\boldsymbol y|^2 / (8s)}\), we have \(\|\nabla^3 G_s\|_1 \le C s^{-3/2}\), so integrating over \(s \in [\tau, \infty)\) gives \(\|\mathcal K_\tau\|_1 \le C \tau^{-1/2}\).
Pointwise, for \(\boldsymbol y \ne \boldsymbol 0\) the substitution \(u = |\boldsymbol y|^2 / (8s)\) yields
\[
  \int_\tau^\infty |\nabla^3 G_s(\boldsymbol y)| \, \mathrm d s
  \le C \int_\tau^\infty \frac{e^{-|\boldsymbol y|^2 / (8s)}}{s^{(n+3)/2}} \, \mathrm d s
  \le C \min\left\{\frac{1}{\tau^{(n+1)/2}}, \, \frac{1}{|\boldsymbol y|^{n+1}}\right\}
  \le \frac{C}{(|\boldsymbol y| + \sqrt \tau)^{n+1}},
\]
which remains valid at \(\boldsymbol y = \boldsymbol 0\) and combines with \(|\nabla G_\tau(\boldsymbol y)| \le C (|\boldsymbol y| + \sqrt \tau)^{-n-1}\) to give \eqref{eq:kernel}.

Fix \(0 < q \le n + 1\) and \(\tau > 0\), and set \(M := \sup_{\boldsymbol z \in \mathbb R^n} \langle \boldsymbol z \rangle^q |\mathbf F(\boldsymbol z)| < \infty\).
When \(|\boldsymbol y| + \sqrt \tau \le 2\), we have \(\langle (\boldsymbol y, \sqrt \tau) \rangle \le \sqrt 5\) and \(\|\mathbf F\|_\infty \le M\), so \eqref{eq:kernel} gives
\[
  \bigl|(\mathcal K_\tau * \mathbf F)(\boldsymbol y)\bigr|
  \le \|\mathcal K_\tau\|_1 \|\mathbf F\|_\infty
  \le \frac{C M}{\sqrt \tau}
  \le \frac{C_q M}{\sqrt \tau \, \langle (\boldsymbol y, \sqrt \tau) \rangle^q}.
\]
When \(|\boldsymbol y| + \sqrt \tau > 2\), we have \(|\boldsymbol y| + \sqrt \tau \ge \frac{2}{3} \langle (\boldsymbol y, \sqrt \tau) \rangle\), and we split the convolution at \(|\boldsymbol z| = \frac{1}{2}(|\boldsymbol y| + \sqrt \tau)\).
On the outer region \(\{|\boldsymbol z| \ge \frac{1}{2}(|\boldsymbol y| + \sqrt \tau)\}\), we have \(\langle \boldsymbol z \rangle^{-q} \le 3^q \langle (\boldsymbol y, \sqrt \tau) \rangle^{-q}\).
On the inner region \(\{|\boldsymbol z| < \frac{1}{2}(|\boldsymbol y| + \sqrt \tau)\}\), the triangle inequality gives \(|\boldsymbol y - \boldsymbol z| + \sqrt \tau \ge \frac{1}{2}(|\boldsymbol y| + \sqrt \tau) \ge \frac{1}{3} \langle (\boldsymbol y, \sqrt \tau) \rangle\) and \(|\boldsymbol z| < \langle (\boldsymbol y, \sqrt \tau) \rangle\), so \eqref{eq:kernel} yields \(|\mathcal K_\tau(\boldsymbol y - \boldsymbol z)| \le C \langle (\boldsymbol y, \sqrt \tau) \rangle^{-n-1}\).
Consequently,
\[
  \bigl|(\mathcal K_\tau * \mathbf F)(\boldsymbol y)\bigr|
  \le \frac{3^q M \|\mathcal K_\tau\|_1}{\langle (\boldsymbol y, \sqrt \tau) \rangle^q}
      + \frac{C M}{\langle (\boldsymbol y, \sqrt \tau) \rangle^{n+1}}
        \int_{|\boldsymbol z| < \langle (\boldsymbol y, \sqrt \tau) \rangle}
          \frac{\mathrm d \boldsymbol z}{\langle \boldsymbol z \rangle^q},
\]
where evaluating the radial integral in the three cases gives
\[
  \begin{aligned}
  &\frac{1}{\langle (\boldsymbol y, \sqrt \tau) \rangle^{n+1}}
   \int_{|\boldsymbol z| < \langle (\boldsymbol y, \sqrt \tau) \rangle}
     \frac{\mathrm d \boldsymbol z}{\langle \boldsymbol z \rangle^q} \\
  \le{}&
  \begin{cases}
    \dfrac{C_q}{\langle (\boldsymbol y, \sqrt \tau) \rangle^{q+1}}
      \le \dfrac{C_q}{\sqrt \tau \, \langle (\boldsymbol y, \sqrt \tau) \rangle^q},
      & 0 < q < n, \\[10pt]
    \dfrac{C \log\bigl(2 + \langle (\boldsymbol y, \sqrt \tau) \rangle\bigr)}
          {\langle (\boldsymbol y, \sqrt \tau) \rangle^{n+1}}
      \le \dfrac{C}{\langle (\boldsymbol y, \sqrt \tau) \rangle^{n + 1/2}},
      & q = n, \\[10pt]
    \dfrac{C_q}{\langle (\boldsymbol y, \sqrt \tau) \rangle^{n+1}}
      \le \dfrac{C_q}{\langle (\boldsymbol y, \sqrt \tau) \rangle^{\max\{q, \, n + 1/2\}}},
      & n < q \le n + 1.
  \end{cases}
  \end{aligned}
\]
Here the first case uses \(\langle (\boldsymbol y, \sqrt \tau) \rangle \ge \sqrt \tau\) and the case \(q = n\) uses \(\langle (\boldsymbol y, \sqrt \tau) \rangle^{-1/2} \log\bigl(2 + \langle (\boldsymbol y, \sqrt \tau) \rangle\bigr) \le C\).
Combining all three cases with \(\|\mathcal K_\tau\|_1 \le C \tau^{-1/2}\) proves the first bound in \eqref{eq:kernelalltime}.

For the heat semigroup \(e^{\tau\Delta} h = G_\tau * h\), since \(\bigl(1 + |\boldsymbol y| / \sqrt \tau\bigr)^{n+1} e^{-|\boldsymbol y|^2 / (4\tau)} \le C\), the rescaled Gaussian kernel \(\tau^{-1/2} G_\tau\) obeys the same pointwise and \(L^1\) bounds as \(\mathcal K_\tau\) in \eqref{eq:kernel}:
\[
  0 \le \frac{1}{\sqrt \tau} G_\tau(\boldsymbol y)
  = \frac{(4\pi)^{-n/2}}{\tau^{(n+1)/2}} e^{-|\boldsymbol y|^2 / (4\tau)}
  \le \frac{C}{(|\boldsymbol y| + \sqrt \tau)^{n+1}},
  \qquad
  \left\|\frac{1}{\sqrt \tau} G_\tau\right\|_1 = \frac{1}{\sqrt \tau}.
\]
Applying the first bound in \eqref{eq:kernelalltime} to \(\tau^{-1/2} G_\tau\) in place of \(\mathcal K_\tau\) and multiplying through by \(\sqrt \tau\) proves the second bound in \eqref{eq:kernelalltime}.
Finally, for \(0 < \tau \le 1\), combining \(\langle (\boldsymbol y, \sqrt \tau) \rangle \ge \langle \boldsymbol y \rangle \ge 1\) with \(1 \le \tau^{-1/2}\) yields \eqref{eq:kernelweighted}.
\end{proof}

We next construct a local mild solution at unit scale \(R = 1\) and rescale to large \(R\).

\begin{lemma}\label{lem:evolution}
Let \(n \ge 5\), and let \(\mathbf u_0\) and \(\mathbf u_{0,R}\) (\(R \ge 1\)) be as in Lemma~\ref{lem:cutoff}.
Then there exists \(T_0 \in (0, 1/4]\) such that \(\mathbf u_{0,1}\) generates a unique bounded mild solution \(\mathbf u_1\) to \eqref{eq:NS} on \([0, 4T_0] \times \mathbb R^n\), which is smooth and divergence-free on \((0, 4T_0] \times \mathbb R^n\) and satisfies
\begin{equation}\label{eq:unitweighted}
  \sup_{0 \le t \le 4T_0} \sup_{\boldsymbol y \in \mathbb R^n}
    \sum_{k=0}^2 \langle \boldsymbol y \rangle^{k+1} |\nabla^k \mathbf u_1(t, \boldsymbol y)|
  \le C.
\end{equation}
Consequently, whenever \(R \ge 1\) satisfies \(R^{-2} \le T_0\), the rescaled velocity field
\begin{equation}\label{eq:rescale}
  \mathbf u_R(t, \boldsymbol y)
  := \frac{1}{R} \mathbf u_1\left(\frac{t}{R^2}, \, \frac{\boldsymbol y}{R}\right)
\end{equation}
is the unique bounded mild solution to \eqref{eq:NS} on \([0, 1] \times \mathbb R^n\) with initial datum \(\mathbf u_{0,R}\), and satisfies uniformly for \(0 \le t \le 1\)
\begin{equation}\label{eq:bgbounds}
  |\nabla^k \mathbf u_R(t, \boldsymbol y)|
  \le \frac{C}{(R + |\boldsymbol y|)^{k+1}}
  \quad (k = 0, 1, 2),
  \qquad
  \|\nabla \mathbf u_R(t)\|_\infty \le \frac{C}{R^2},
\end{equation}
with a constant \(C\) independent of \(R\).
\end{lemma}

\begin{proof}
For \(T \in (0, 1]\), let \(X_T\) be the Banach space of fields \(\mathbf v \in C([0, T]; C_0(\mathbb R^n; \mathbb R^n))\) such that \(\nabla^k \mathbf v \in C([0, T]; C_0)\) for \(k \in \{1, 2\}\) and
\[
  \|\mathbf v\|_{X_T}
  := \sup_{0 \le t \le T} \sup_{\boldsymbol y \in \mathbb R^n}
     \sum_{k=0}^2 \langle \boldsymbol y \rangle^{k+1} |\nabla^k \mathbf v(t, \boldsymbol y)| < \infty.
\]
We seek \(\mathbf u_1 \in X_{4T_0}\) solving the mild integral equation
\begin{equation}\label{eq:mild}
  \mathbf u_1(t)
  = e^{t\Delta} \mathbf u_{0,1}
    - \int_0^t \mathcal K_{t-s} * (\mathbf u_1 \otimes \mathbf u_1)(s) \, \mathrm d s.
\end{equation}

Since \(\nabla^k \mathbf u_{0,1} \in C_0\) for \(k \in \{0, 1, 2\}\), strong continuity of \(e^{t\Delta}\) on \(C_0(\mathbb R^n)\) together with \eqref{eq:cutoffbounds} and \eqref{eq:kernelweighted} implies \(e^{\cdot\Delta} \mathbf u_{0,1} \in X_T\) with \(\|e^{\cdot\Delta} \mathbf u_{0,1}\|_{X_T} \le c_0\).
For any \(\mathbf v, \mathbf w \in X_T\) and \(t \in [0, T]\), we use \(\nabla^k (\mathcal K_\tau * \mathbf F) = \mathcal K_\tau * \nabla^k \mathbf F\) and apply \eqref{eq:kernelweighted} at \(q = k + 1\) to \(\nabla^k (\mathbf v \otimes \mathbf w)\) for \(k \in \{0, 1, 2\}\), which by the product rule gives
\[
  \begin{aligned}
  &\sup_{\boldsymbol y \in \mathbb R^n}
    \sum_{k=0}^2 \langle \boldsymbol y \rangle^{k+1}
          \left|\nabla^k \int_0^t \mathcal K_{t-s} * (\mathbf v \otimes \mathbf w)(s) \, \mathrm d s\right| \\
  \le{}& C \int_0^t \frac{1}{\sqrt{t - s}}
    \sup_{\boldsymbol y \in \mathbb R^n}
    \sum_{k=0}^2 \langle \boldsymbol y \rangle^{k+1} |\nabla^k (\mathbf v \otimes \mathbf w)(s, \boldsymbol y)| \, \mathrm d s \\
  \le{}& \frac{c_1}{2} \|\mathbf v\|_{X_T} \|\mathbf w\|_{X_T}
    \int_0^t \frac{\mathrm d s}{\sqrt{t - s}}
    = c_1 \sqrt t \, \|\mathbf v\|_{X_T} \|\mathbf w\|_{X_T}.
  \end{aligned}
\]
Substituting \(\tau = t - s\) and applying dominated convergence with \(\|\mathcal K_\tau\|_1 \le C \tau^{-1/2} \in L^1(0, T)\) shows that the Duhamel integral and its spatial derivatives up to order two are continuous in time with values in unweighted \(C_0(\mathbb R^n)\); combined with the weighted supremum bound above, this shows that the Duhamel integral belongs to \(X_T\).
Writing
\[
  \mathbf v \otimes \mathbf v - \mathbf w \otimes \mathbf w
  = (\mathbf v - \mathbf w) \otimes \mathbf v + \mathbf w \otimes (\mathbf v - \mathbf w)
\]
and fixing \(T_0 \in (0, 1/4]\) with \(8 c_0 c_1 \sqrt{T_0} \le 1/2\), we see that the right-hand side of \eqref{eq:mild} is a contraction on \(\{\mathbf v \in X_{4T_0} \mid \|\mathbf v\|_{X_{4T_0}} \le 2 c_0\}\), and its fixed point \(\mathbf u_1\) satisfies \eqref{eq:unitweighted}.

Uniqueness among bounded mild solutions of \eqref{eq:mild} follows by iterating the resulting integral inequality with kernel \((t - s)^{-1/2}\) once (as in Proposition~\ref{prop:defect}) and applying Gr\"onwall's inequality.
Since \(\langle \boldsymbol y \rangle^{-1} \in L^{n+1}(\mathbb R^n)\), we have \(\mathbf u_1 \in C([0, 4T_0]; L^{n+1}(\mathbb R^n)) \cap L^\infty((0, 4T_0) \times \mathbb R^n)\); by the subcritical \(L^{n+1}\) mild-solution regularity theory \cite[Theorem~4]{Giga1986} (together with \(\operatorname{div} \mathbf u_{0,1} = 0\) and \eqref{eq:kerneldef}), \(\mathbf u_1\) is smooth and divergence-free on \((0, 4T_0] \times \mathbb R^n\).
Finally, when \(R^{-2} \le T_0\) and \(0 \le t \le 1\), we have \(t / R^2 \le T_0 \le 4T_0\), so parabolic scaling of \eqref{eq:cutoffbounds} and \eqref{eq:mild} shows that \(\mathbf u_R\) in \eqref{eq:rescale} is the unique bounded mild solution with initial datum \(\mathbf u_{0,R}\).
Moreover, for \(k \in \{0, 1, 2\}\), \eqref{eq:unitweighted} gives
\[
  |\nabla^k \mathbf u_R(t, \boldsymbol y)|
  = \frac{1}{R^{k+1}} \left|(\nabla^k \mathbf u_1)\left(\frac{t}{R^2}, \, \frac{\boldsymbol y}{R}\right)\right|
  \le \frac{C}{R^{k+1} \langle \boldsymbol y / R \rangle^{k+1}}
  \le \frac{C}{(R + |\boldsymbol y|)^{k+1}},
\]
proving \eqref{eq:bgbounds}.
\end{proof}

\section{Stationary head-pressure interior regularity}\label{app:stationary}

Interior regularity for the stationary Navier--Stokes equations in dimensions \(m \ge 5\) under a one-sided head-pressure bound originates in Frehse and R\r u\v zi\v cka \cite{FR94maximum,FR98criterion}, who also recovered the local energy inequality from weak regularity \cite[Theorems~1.8 and~4.2]{FR94bounded}.
Quantitative Morrey-space embeddings and Stokes estimates in this setting were developed by Li and Yang \cite{LY22}, and a related quantitative estimate assuming a local energy inequality appears in \cite[Proposition~4.1]{GLLWX26transverse}.
In Section~\ref{subsec:regularity}, Proposition~\ref{thm:stationaryreg} is applied to the lifted \(\mathbb R^{n+1}\) stationary solution \((\mathbf V, \widetilde \Pi)\), and its local energy equality is derived directly within the proof rather than assumed a priori.
Throughout this appendix, unadorned balls \(B_r, B_\rho(\boldsymbol z)\), Riesz transforms, and tensor conventions are understood in the ambient dimension \(m\).

Following \cite{FR94maximum,FR98criterion}, we derive scale-critical \(L^{2,m-2}\) Morrey bounds from a one-sided upper bound on the steady total head pressure \(\varpi := \pi + \frac{1}{2} |\mathbf v|^2\) by testing the weak pressure equation against regularized versions of \(|\boldsymbol \xi - \boldsymbol z|^{4-m}\), writing \(\langle \boldsymbol \xi - \boldsymbol z \rangle_\delta := (|\boldsymbol \xi - \boldsymbol z|^2 + \delta^2)^{1/2}\) for \(\delta > 0\).

\begin{lemma}\label{lem:prestest}
Let \(m \ge 5\), \(0 < R_1 < R_2\), and \(\boldsymbol z_0 \in \mathbb R^m\).
Suppose \(\mathbf v \in L^2(B_{R_2}(\boldsymbol z_0); \mathbb R^m)\) and \(\pi \in L^1(B_{R_2}(\boldsymbol z_0))\) satisfy the weak pressure equation
\begin{equation}\label{eq:weakpreseq}
  \int_{B_{R_2}(\boldsymbol z_0)} \pi \Delta \psi \, \mathrm d \boldsymbol \xi
  + \sum_{i=1}^m \sum_{j=1}^m \int_{B_{R_2}(\boldsymbol z_0)} v_i v_j \partial_i \partial_j \psi \, \mathrm d \boldsymbol \xi = 0,
  \qquad \forall \, \psi \in C_c^\infty(B_{R_2}(\boldsymbol z_0)).
\end{equation}
Set \(\varpi := \pi + \frac{1}{2} |\mathbf v|^2\).
Then for every \(\boldsymbol z \in B_{R_1}(\boldsymbol z_0)\),
\begin{equation}\label{eq:prestestbound}
  \begin{aligned}
  &\int_{B_{R_1}(\boldsymbol z_0)}
    \left(\frac{|\varpi(\boldsymbol \xi)|}{|\boldsymbol \xi - \boldsymbol z|^{m-2}}
          + \frac{[(\boldsymbol \xi - \boldsymbol z) \cdot \mathbf v(\boldsymbol \xi)]^2}
                 {|\boldsymbol \xi - \boldsymbol z|^m}\right) \, \mathrm d \boldsymbol \xi \\
  \le{}& C_m \int_{B_{R_2}(\boldsymbol z_0)}
        \frac{\varpi_+(\boldsymbol \xi)}{|\boldsymbol \xi - \boldsymbol z|^{m-2}} \, \mathrm d \boldsymbol \xi \\
  &+ \frac{C_m}{(R_2 - R_1)^{m-2}}
        \int_{B_{R_2}(\boldsymbol z_0) \setminus B_{R_1}(\boldsymbol z_0)}
          \bigl(|\varpi(\boldsymbol \xi)| + |\mathbf v(\boldsymbol \xi)|^2\bigr) \, \mathrm d \boldsymbol \xi.
  \end{aligned}
\end{equation}
Moreover,
\begin{equation}\label{eq:prestestmorrey}
  \begin{aligned}
  \|\mathbf v\|_{L^{2,m-2}(B_{R_1}(\boldsymbol z_0))}^2
  \le{}& C_m \sup_{\boldsymbol z \in B_{(R_1+R_2)/2}(\boldsymbol z_0)}
        \int_{B_{R_2}(\boldsymbol z_0)}
          \frac{\varpi_+(\boldsymbol \xi)}{|\boldsymbol \xi - \boldsymbol z|^{m-2}} \, \mathrm d \boldsymbol \xi \\
  &+ \frac{C_m}{(R_2 - R_1)^{m-2}}
        \int_{B_{R_2}(\boldsymbol z_0)}
          \bigl(|\varpi(\boldsymbol \xi)| + |\mathbf v(\boldsymbol \xi)|^2\bigr) \, \mathrm d \boldsymbol \xi.
  \end{aligned}
\end{equation}
\end{lemma}

\begin{proof}
Set \(r_1 := (R_1 + R_2)/2\) and \(r_2 := (R_1 + 3R_2)/4\), so that \(R_1 < r_1 < r_2 < R_2\), and choose \(\eta \in C_c^\infty(B_{R_2}(\boldsymbol z_0))\) with \(0 \le \eta \le 1\), \(\eta|_{B_{r_2}(\boldsymbol z_0)} = 1\), and \(|\nabla^k \eta| \le C (R_2 - R_1)^{-k}\) for \(k = 1, 2\).
For any \(\boldsymbol z \in B_{r_1}(\boldsymbol z_0)\) and \(\delta > 0\), direct differentiation gives
\[
  \begin{aligned}
  \partial_i \partial_j \left (\frac{1}{\langle \boldsymbol \xi - \boldsymbol z \rangle_\delta^{m-4}} \right )
  &= - \frac{(m - 4) \delta_{ij}}{\langle \boldsymbol \xi - \boldsymbol z \rangle_\delta^{m-2}}
     + \frac{(m - 4)(m - 2) (\xi_i - z_i)(\xi_j - z_j)}{\langle \boldsymbol \xi - \boldsymbol z \rangle_\delta^{m}}, \\
  \Delta \left(\frac{1}{\langle \boldsymbol \xi - \boldsymbol z \rangle_\delta^{m-4}}\right)
  &= - \frac{2(m - 4)}{\langle \boldsymbol \xi - \boldsymbol z \rangle_\delta^{m-2}}
     - \frac{(m - 4)(m - 2) \delta^2}{\langle \boldsymbol \xi - \boldsymbol z \rangle_\delta^{m}}.
  \end{aligned}
\]
Here the \(1 : 2\) ratio between the coefficient \(-(m - 4)\) of \(\delta_{ij}\) in the Hessian and \(-2(m - 4)\) in the Laplacian combines \(|\mathbf v|^2\) and \(2\pi\) into \(2\varpi = 2\pi + |\mathbf v|^2\), while reducing the singularity to the critical Morrey kernel \(|\boldsymbol \xi - \boldsymbol z|^{2-m}\).
Testing \eqref{eq:weakpreseq} with \(\psi(\boldsymbol \xi) = \eta(\boldsymbol \xi)^2 \langle \boldsymbol \xi - \boldsymbol z \rangle_\delta^{4-m}\), we separate the main terms in which both derivatives fall on \(\langle \boldsymbol \xi - \boldsymbol z \rangle_\delta^{4-m}\) from the cutoff terms supported in \(\operatorname{supp} \nabla \eta \subset B_{R_2}(\boldsymbol z_0) \setminus B_{r_2}(\boldsymbol z_0)\).
For the main terms, substituting \(\pi = \varpi - \frac{1}{2} |\mathbf v|^2\) and \(-\varpi = |\varpi| - 2\varpi_+\) gives
\[
  \begin{aligned}
  &\int_{B_{R_2}(\boldsymbol z_0)} \eta^2 \Bigl[ 
     -\varpi \frac{2(m - 4)}{\langle \boldsymbol \xi - \boldsymbol z \rangle_\delta^{m-2}}
     - \frac{(m - 4)(m - 2) \pi \delta^2}{\langle \boldsymbol \xi - \boldsymbol z \rangle_\delta^{m}} 
  + \frac{(m - 4)(m - 2) [(\boldsymbol \xi - \boldsymbol z) \cdot \mathbf v]^2}
       {\langle \boldsymbol \xi - \boldsymbol z \rangle_\delta^{m}}
   \Bigr] \, \mathrm d \boldsymbol \xi \\
  ={}& \int_{B_{R_2}(\boldsymbol z_0)} \eta^2 \Bigl[
     \frac{2(m - 4)\bigl(|\varpi| - 2\varpi_+\bigr)}{\langle \boldsymbol \xi - \boldsymbol z \rangle_\delta^{m-2}} 
  + \frac{(m - 4)(m - 2)\bigl(|\varpi| - 2\varpi_+ + \tfrac{1}{2}|\mathbf v|^2\bigr)
       \delta^2}{\langle \boldsymbol \xi - \boldsymbol z \rangle_\delta^{m}} \\
  &\quad + \frac{(m - 4)(m - 2)\bigl[(\boldsymbol \xi - \boldsymbol z) \cdot \mathbf v\bigr]^2}
       {\langle \boldsymbol \xi - \boldsymbol z \rangle_\delta^{m}}
   \Bigr] \, \mathrm d \boldsymbol \xi.
  \end{aligned}
\]
On \(\operatorname{supp} \nabla \eta \subset B_{R_2}(\boldsymbol z_0) \setminus B_{r_2}(\boldsymbol z_0)\), every \(\boldsymbol z \in B_{r_1}(\boldsymbol z_0)\) satisfies \(|\boldsymbol \xi - \boldsymbol z| \ge r_2 - r_1 = (R_2 - R_1)/4\).
Every term in \(\partial_i \partial_j \psi\) and \(\Delta \psi\) where at least one derivative hits \(\eta^2\) is therefore bounded by \(C_m (R_2 - R_1)^{2-m}\).
Using \(|\pi| \le |\varpi| + \frac{1}{2} |\mathbf v|^2\), their contribution to \eqref{eq:weakpreseq} is bounded by
\[
  \frac{C_m}{(R_2 - R_1)^{m-2}}
  \int_{B_{R_2}(\boldsymbol z_0) \setminus B_{R_1}(\boldsymbol z_0)}
    \bigl(|\varpi(\boldsymbol \xi)| + |\mathbf v(\boldsymbol \xi)|^2\bigr) \, \mathrm d \boldsymbol \xi.
\]
Dropping the nonnegative terms \((|\varpi| + \frac{1}{2} |\mathbf v|^2) \delta^2 \ge 0\) in the main integral, bounding \(\delta^2 \langle \boldsymbol \xi - \boldsymbol z \rangle_\delta^{-m} \le \langle \boldsymbol \xi - \boldsymbol z \rangle_\delta^{2-m}\) in the \(-2\varpi_+\) term, and sending \(\delta \downarrow 0\) via Fatou's lemma on \(B_{r_2}(\boldsymbol z_0)\) (where \(\eta = 1\)), we obtain for every \(\boldsymbol z \in B_{r_1}(\boldsymbol z_0)\):
\begin{equation}\label{eq:midballpot}
  \begin{aligned}
  &\int_{B_{r_2}(\boldsymbol z_0)}
    \left(\frac{|\varpi(\boldsymbol \xi)|}{|\boldsymbol \xi - \boldsymbol z|^{m-2}}
          + \frac{[(\boldsymbol \xi - \boldsymbol z) \cdot \mathbf v(\boldsymbol \xi)]^2}
                 {|\boldsymbol \xi - \boldsymbol z|^m}\right) \, \mathrm d \boldsymbol \xi \\
  \le{}& C_m \int_{B_{R_2}(\boldsymbol z_0)}
        \frac{\varpi_+(\boldsymbol \xi)}{|\boldsymbol \xi - \boldsymbol z|^{m-2}} \, \mathrm d \boldsymbol \xi \\
  &+ \frac{C_m}{(R_2 - R_1)^{m-2}}
        \int_{B_{R_2}(\boldsymbol z_0) \setminus B_{R_1}(\boldsymbol z_0)}
          \bigl(|\varpi(\boldsymbol \xi)| + |\mathbf v(\boldsymbol \xi)|^2\bigr) \, \mathrm d \boldsymbol \xi.
  \end{aligned}
\end{equation}
Since \(B_{R_1}(\boldsymbol z_0) \subset B_{r_1}(\boldsymbol z_0) \subset B_{r_2}(\boldsymbol z_0)\), \eqref{eq:prestestbound} follows immediately.

To deduce \eqref{eq:prestestmorrey}, fix any \(\boldsymbol z_1 \in B_{R_1}(\boldsymbol z_0)\) and \(0 < \rho \le 2 R_1\).
If \(\rho \ge (R_2 - R_1)/4\), we have directly
\[
  \frac{1}{\rho^{m-2}} \int_{B_{R_1}(\boldsymbol z_0) \cap B_\rho(\boldsymbol z_1)} |\mathbf v|^2 \, \mathrm d \boldsymbol \xi
  \le \frac{4^{m-2}}{(R_2 - R_1)^{m-2}}
      \int_{B_{R_2}(\boldsymbol z_0)} |\mathbf v|^2 \, \mathrm d \boldsymbol \xi.
\]
If \(0 < \rho < (R_2 - R_1)/4\), then \(B_\rho(\boldsymbol z_1) \subset B_{2\rho}(\boldsymbol z_1) \subset B_{r_1}(\boldsymbol z_0) \subset B_{r_2}(\boldsymbol z_0)\), and \(B_\rho(\boldsymbol \xi) \subset B_{2\rho}(\boldsymbol z_1)\) for every \(\boldsymbol \xi \in B_\rho(\boldsymbol z_1)\).
Passing to polar coordinates \(\boldsymbol z = \boldsymbol \xi - r \boldsymbol \omega\) (\(0 < r < \rho\), \(\boldsymbol \omega \in \mathbb S^{m-1}\)) and using rotational symmetry
\[
  \int_{\mathbb S^{m-1}} (\boldsymbol \omega \cdot \mathbf v(\boldsymbol \xi))^2 \, \mathrm d \sigma(\boldsymbol \omega)
  = \frac{|\mathbb S^{m-1}|}{m} |\mathbf v(\boldsymbol \xi)|^2
\]
together with \(\int_0^\rho r \, \mathrm d r = \frac{1}{2} \rho^2\), Tonelli's theorem gives
\[
  \begin{aligned}
  \frac{|\mathbb S^{m-1}|}{2m} \rho^2 \int_{B_\rho(\boldsymbol z_1)} |\mathbf v(\boldsymbol \xi)|^2 \, \mathrm d \boldsymbol \xi
  &= \int_{B_\rho(\boldsymbol z_1)} \! \int_{B_\rho(\boldsymbol \xi)}
       \frac{[(\boldsymbol \xi - \boldsymbol z) \cdot \mathbf v(\boldsymbol \xi)]^2}{|\boldsymbol \xi - \boldsymbol z|^m}
       \, \mathrm d \boldsymbol z \, \mathrm d \boldsymbol \xi \\
  &\le \int_{B_{2\rho}(\boldsymbol z_1)} \! \int_{B_{r_2}(\boldsymbol z_0)}
       \frac{[(\boldsymbol \xi - \boldsymbol z) \cdot \mathbf v(\boldsymbol \xi)]^2}{|\boldsymbol \xi - \boldsymbol z|^m}
       \, \mathrm d \boldsymbol \xi \, \mathrm d \boldsymbol z.
  \end{aligned}
\]
Bounding the inner integral via \eqref{eq:midballpot}, dividing by \(\rho^m\) using \(|B_{2\rho}(\boldsymbol z_1)| = C_m \rho^m\), and taking the supremum over \(\boldsymbol z_1 \in B_{R_1}(\boldsymbol z_0)\) and \(0 < \rho \le 2 R_1\) yields \eqref{eq:prestestmorrey}.
\end{proof}

We next record the Morrey-space consequence of Ledoux's improved Sobolev inequality \cite{Ledoux03}.

\begin{lemma}[Ledoux \cite{Ledoux03}]\label{lem:morreyl4}
Let \(m \ge 3\) and \(d \ge 1\).
For every \(\mathbf f \in H^1(\mathbb R^m; \mathbb R^d)\) with \(\|\mathbf f\|_{L^{2,m-2}(\mathbb R^m)} < \infty\),
\begin{equation}\label{eq:l4bound}
  \|\mathbf f\|_{L^4(\mathbb R^m)}^4
  \le C_m \|\mathbf f\|_{L^{2,m-2}(\mathbb R^m)}^2 \|\nabla \mathbf f\|_{L^2(\mathbb R^m)}^2.
\end{equation}
In particular, for any concentric balls \(B_{R_1}(\boldsymbol z_0) \subset B_{R_2}(\boldsymbol z_0)\) (\(0 < R_1 < R_2\)) and any \(\mathbf f \in H^1(B_{R_2}(\boldsymbol z_0); \mathbb R^d)\),
\[
  \begin{aligned}
  &\|\mathbf f\|_{L^4(B_{R_1}(\boldsymbol z_0))}^4 \\
  \le{}& C_m \left(\frac{R_2}{R_2 - R_1}\right)^{m-2}
       \|\mathbf f\|_{L^{2,m-2}(B_{R_2}(\boldsymbol z_0))}^2
       \left(\|\nabla \mathbf f\|_{L^2(B_{R_2}(\boldsymbol z_0))}^2
             + \frac{\|\mathbf f\|_{L^2(B_{R_2}(\boldsymbol z_0))}^2}{(R_2 - R_1)^2}\right).
  \end{aligned}
\]
\end{lemma}

\begin{proof}
Assume \(\|\mathbf f\|_{L^{2,m-2}(\mathbb R^m)} > 0\), and fix a nonnegative radial mollifier \(\zeta_1 \in C_c^\infty(B_1)\) with \(\int_{\mathbb R^m} \zeta_1(\boldsymbol \xi) \, \mathrm d \boldsymbol \xi = 1\), setting \(\zeta_\delta(\boldsymbol \xi) := \delta^{-m} \zeta_1(\boldsymbol \xi / \delta)\) for \(\delta > 0\).
By the Cauchy--Schwarz inequality and the standard Riesz-potential bound in terms of the Hardy--Littlewood maximal operator \(\mathcal M\), for almost every \(\boldsymbol \xi \in \mathbb R^m\) and every \(\delta > 0\),
\[
  \begin{aligned}
  |(\mathbf f * \zeta_\delta)(\boldsymbol \xi)|
  &\le \frac{C_m}{\delta^{m/2}} \|\mathbf f\|_{L^2(B_\delta(\boldsymbol \xi))}
   \le \frac{C_m}{\delta} \|\mathbf f\|_{L^{2,m-2}(\mathbb R^m)}, \\
  |\mathbf f(\boldsymbol \xi) - (\mathbf f * \zeta_\delta)(\boldsymbol \xi)|
  &\le C_m \int_{B_\delta(\boldsymbol \xi)}
       \frac{|\nabla \mathbf f(\boldsymbol z)|}{|\boldsymbol \xi - \boldsymbol z|^{m-1}} \, \mathrm d \boldsymbol z
   \le C_m \delta \, \mathcal M(|\nabla \mathbf f|)(\boldsymbol \xi).
  \end{aligned}
\]
Let \(C_m \ge 1\) be a common dimensional constant valid in both estimates.
Almost every \(\boldsymbol \xi \in \mathbb R^m\) with \(\mathbf f(\boldsymbol \xi) \ne \boldsymbol 0\) is a Lebesgue point of \(\mathbf f\) satisfying \(\mathcal M(|\nabla \mathbf f|)(\boldsymbol \xi) < \infty\); taking the scale \(\delta = 2 C_m \|\mathbf f\|_{L^{2,m-2}(\mathbb R^m)} / |\mathbf f(\boldsymbol \xi)|\) in the triangle inequality \(|\mathbf f(\boldsymbol \xi)| \le |(\mathbf f * \zeta_\delta)(\boldsymbol \xi)| + |\mathbf f(\boldsymbol \xi) - (\mathbf f * \zeta_\delta)(\boldsymbol \xi)|\) yields
\[
  |\mathbf f(\boldsymbol \xi)|^2
  \le 4 C_m^2 \|\mathbf f\|_{L^{2,m-2}(\mathbb R^m)} \mathcal M(|\nabla \mathbf f|)(\boldsymbol \xi).
\]
Squaring, integrating over \(\mathbb R^m\), and using the \(L^2(\mathbb R^m)\) boundedness of \(\mathcal M\) gives \eqref{eq:l4bound}.

For the local bound on \(B_{R_1}(\boldsymbol z_0) \subset B_{R_2}(\boldsymbol z_0)\), pick \(\eta \in C_c^\infty(B_{(R_1+R_2)/2}(\boldsymbol z_0))\) with \(0 \le \eta \le 1\), \(\eta|_{B_{R_1}(\boldsymbol z_0)} = 1\), and \(|\nabla \eta| \le C / (R_2 - R_1)\).
For any ball \(B_\rho(\boldsymbol z_1)\) meeting \(\operatorname{supp} \eta \subset B_{(R_1+R_2)/2}(\boldsymbol z_0)\), when \(0 < \rho \le (R_2 - R_1)/2\) we have \(\boldsymbol z_1 \in B_{R_2}(\boldsymbol z_0)\) and \(\operatorname{supp}(\eta \mathbf f) \subset B_{R_2}(\boldsymbol z_0)\), so
\[
  \frac{1}{\rho^{m-2}} \int_{B_\rho(\boldsymbol z_1)} |\eta \mathbf f|^2 \, \mathrm d \boldsymbol \xi
  = \frac{1}{\rho^{m-2}} \int_{B_{R_2}(\boldsymbol z_0) \cap B_\rho(\boldsymbol z_1)} |\eta \mathbf f|^2 \, \mathrm d \boldsymbol \xi
  \le \|\mathbf f\|_{L^{2,m-2}(B_{R_2}(\boldsymbol z_0))}^2;
\]
whereas for \(\rho > (R_2 - R_1)/2\) we use \(R_2^{2-m} \|\mathbf f\|_{L^2(B_{R_2}(\boldsymbol z_0))}^2 \le \|\mathbf f\|_{L^{2,m-2}(B_{R_2}(\boldsymbol z_0))}^2\) to bound
\[
  \frac{1}{\rho^{m-2}} \int_{B_\rho(\boldsymbol z_1)} |\eta \mathbf f|^2 \, \mathrm d \boldsymbol \xi
  \le \left(\frac{2 R_2}{R_2 - R_1}\right)^{m-2} \|\mathbf f\|_{L^{2,m-2}(B_{R_2}(\boldsymbol z_0))}^2.
\]
Applying \eqref{eq:l4bound} to \(\eta \mathbf f\) completes the proof.
\end{proof}

\begin{proposition}\label{thm:stationaryreg}
Let \(m \ge 5\), let \(D \subset \mathbb R^m\) be open, and let \(\mathbf v \in H^1_{\mathrm{loc}}(D; \mathbb R^m)\) and \(\pi \in L^1_{\mathrm{loc}}(D)\) solve the unforced stationary Navier--Stokes system in \(\mathscr{D}'(D)\):
\begin{equation}\label{eq:steadyNS}
  -\Delta \mathbf v + (\mathbf v \cdot \nabla) \mathbf v + \nabla \pi = \boldsymbol 0,
  \qquad \operatorname{div} \mathbf v = 0.
\end{equation}
If \(\varpi := \pi + \frac{1}{2} |\mathbf v|^2\) is essentially bounded above on every compact subset of \(D\), then \(\mathbf v\) and \(\pi\) admit \(C^\infty(D)\) representatives.
More quantitatively, if \(B_2(\boldsymbol z_0) \subset D\) and there exists \(M \ge 1\) such that
\[
  \pi + \frac{1}{2} |\mathbf v|^2 \le M \quad \text{a.e. in } B_2(\boldsymbol z_0),
  \qquad
  \|\mathbf v\|_{H^1(B_2(\boldsymbol z_0))} + \|\pi\|_{L^1(B_2(\boldsymbol z_0))} \le M,
\]
then for every integer \(k \ge 0\),
\begin{equation}\label{eq:quantitativeCk}
  \|\mathbf v\|_{C^k(B_1(\boldsymbol z_0))} \le C_{m,k,M}.
\end{equation}
\end{proposition}

\begin{proof}
By translation invariance we may take \(\boldsymbol z_0 = \boldsymbol 0\) and establish \(\|\mathbf v\|_{C^k(B_1)} \le C_{m,k,M}\) from the hypotheses on \(B_2\); qualitative smoothness on \(D\) then follows by rescaling any ball \(B_{2\rho}(\boldsymbol z) \Subset D\) to \(B_2\) and recovering \(\pi \in C^\infty(D)\) via \(\nabla \pi = \Delta \mathbf v - (\mathbf v \cdot \nabla) \mathbf v\).
Throughout, \(C\) denotes constants depending only on \(m\), \(M\), and the chosen concentric radii.

\smallskip
\noindent\emph{Step 1.}
Let \(0 < r < 2\) and choose intermediate radii \(r < \rho_1 < \rho_2 < 2\).
Testing \eqref{eq:steadyNS} on \(B_2\) against \(\nabla \psi\) (\(\psi \in C_c^\infty(B_2)\)) verifies \eqref{eq:weakpreseq} on \(B_2\).
Since \(\varpi_+ \le M\) a.e.\ on \(B_2\), applying Lemma~\ref{lem:prestest} on \(B_{\rho_2} \subset B_2\) gives
\[
  \sup_{\boldsymbol z \in B_{\rho_2}} \int_{B_{\rho_2}}
    \frac{|\varpi(\boldsymbol \xi)|}{|\boldsymbol \xi - \boldsymbol z|^{m-2}} \, \mathrm d \boldsymbol \xi
  + \|\mathbf v\|_{L^{2,m-2}(B_{\rho_2})}^2
  \le C.
\]
The local estimate in Lemma~\ref{lem:morreyl4} on \(B_{\rho_1} \subset B_{\rho_2}\) then yields \(\|\mathbf v\|_{L^4(B_{\rho_1})} \le C\).
On \(B_{\rho_1}\), the difference
\[
  \pi - \sum_{i=1}^m \sum_{j=1}^m \mathcal R_i \mathcal R_j (\mathbf 1_{B_{\rho_1}} v_i v_j)
\]
is harmonic with bounded \(L^1(B_{\rho_1})\) norm, so Riesz-transform boundedness on \(L^2(\mathbb R^m)\) and interior harmonic estimates on \(B_r \Subset B_{\rho_1}\) give \(\|\pi\|_{L^2(B_r)} + \|\varpi\|_{L^2(B_r)} \le C\).
Hence, for every \(0 < r < 2\),
\begin{equation}\label{eq:B35pot}
  \sup_{\boldsymbol z \in B_r} \int_{B_r}
    \frac{|\varpi(\boldsymbol \xi)|}{|\boldsymbol \xi - \boldsymbol z|^{m-2}} \, \mathrm d \boldsymbol \xi
  + \|\mathbf v\|_{L^{2,m-2}(B_r)}^2
  + \|\mathbf v\|_{L^4(B_r)}
  + \|\varpi\|_{L^2(B_r)}
  \le C.
\end{equation}
Consequently, \eqref{eq:steadyNS} takes the rotational form \(-\Delta \mathbf v + \Omega_{\mathbf v} \mathbf v + \nabla \varpi = \boldsymbol 0\) on \(B_r\), where \(\Omega_{\mathbf v} \mathbf v \in L^{4/3}(B_r; \mathbb R^m)\) is pointwise orthogonal to \(\mathbf v\).
Testing against compactly supported mollifications of \(\mathbf v \varphi \in H^1_0(B_r; \mathbb R^m) \cap L^4(B_r; \mathbb R^m)\) yields the local energy equality
\begin{equation}\label{eq:localenergy}
  \int_{B_r} \nabla \mathbf v : \nabla (\mathbf v \varphi) \, \mathrm d \boldsymbol \xi
  = \int_{B_r} \varpi \mathbf v \cdot \nabla \varphi \, \mathrm d \boldsymbol \xi,
  \qquad \forall \, \varphi \in C_c^\infty(B_r).
\end{equation}

\smallskip
\noindent\emph{Step 2.}
Choose concentric radii \(1 < r_0 < r < 2\) once and for all (for instance, trisecting \((1, 2)\)), so that \eqref{eq:B35pot}--\eqref{eq:localenergy} hold on \(B_r\).
Because \(\varpi \in L^2(B_r)\) has no gradient bound in \eqref{eq:localenergy}, we introduce the Dirichlet potential \(\Gamma \in W^{2,2}(B_r) \cap H^1_0(B_r)\) solving \(-\Delta \Gamma = \varpi\) in \(B_r\) with \(\Gamma|_{\partial B_r} = 0\).
Combining the Green-function bound \(0 \le G_{B_r}(\boldsymbol z, \boldsymbol \xi) \le C_m |\boldsymbol \xi - \boldsymbol z|^{2-m}\) with \eqref{eq:B35pot} gives
\begin{equation}\label{eq:GammaLinfty}
  \|\Gamma\|_{L^\infty(B_r)}
  \le C_m \sup_{\boldsymbol z \in B_r}
      \int_{B_r} \frac{|\varpi(\boldsymbol \xi)|}{|\boldsymbol \xi - \boldsymbol z|^{m-2}} \, \mathrm d \boldsymbol \xi
  \le C.
\end{equation}
Set \(\rho_0 := \frac{1}{4}(r - r_0)\), so that \(B_{2\rho}(\boldsymbol z) \subset B_{(r_0+r)/2} \Subset B_r\) whenever \(\boldsymbol z \in B_{r_0}\) and \(0 < \rho \le \rho_0\).
For any \(\boldsymbol z \in B_{r_0}\) and \(0 < \rho \le \rho_0\), define
\[
  I(\rho)
  := \int_{B_\rho(\boldsymbol z)}
     \left(\frac{|\nabla \mathbf v(\boldsymbol \xi)|^2}{|\boldsymbol \xi - \boldsymbol z|^{m-4}}
           + \frac{|\mathbf v(\boldsymbol \xi)|^2 + |\nabla \Gamma(\boldsymbol \xi)|^2 + |\varpi(\boldsymbol \xi)|}
                  {|\boldsymbol \xi - \boldsymbol z|^{m-2}}\right) \, \mathrm d \boldsymbol \xi,
\]
and pick \(\eta_\rho \in C_c^\infty(B_{2\rho}(\boldsymbol z))\) with \(0 \le \eta_\rho \le 1\), \(\eta_\rho|_{B_\rho(\boldsymbol z)} = 1\), and \(|\nabla^k \eta_\rho| \le C \rho^{-k}\) (\(k = 1, 2\)).

First, applying \eqref{eq:prestestbound} on \(B_\rho(\boldsymbol z) \subset B_{2\rho}(\boldsymbol z)\) with \(\varpi_+ \le M\), using \(\rho^{2-m} \le 2^{m-2} |\boldsymbol \xi - \boldsymbol z|^{2-m}\) on \(B_{2\rho}(\boldsymbol z) \setminus B_\rho(\boldsymbol z)\), and adding the annular integral of \(|\varpi| |\boldsymbol \xi - \boldsymbol z|^{2-m}\) to both sides gives
\begin{equation}\label{eq:holevarpi}
  \int_{B_{2\rho}(\boldsymbol z)} \frac{|\varpi(\boldsymbol \xi)|}{|\boldsymbol \xi - \boldsymbol z|^{m-2}} \, \mathrm d \boldsymbol \xi
  \le C \rho^2
      + C \int_{B_{2\rho}(\boldsymbol z) \setminus B_\rho(\boldsymbol z)}
          \frac{|\varpi(\boldsymbol \xi)| + |\mathbf v(\boldsymbol \xi)|^2}{|\boldsymbol \xi - \boldsymbol z|^{m-2}} \, \mathrm d \boldsymbol \xi.
\end{equation}
Second, letting \(\bar \Gamma_\rho\) denote the mean of \(\Gamma\) over \(B_{2\rho}(\boldsymbol z) \setminus B_\rho(\boldsymbol z)\), we test \(-\Delta \Gamma = \varpi\) on \(B_r\) against \((\Gamma - \bar \Gamma_\rho) \eta_\rho^2 \langle \boldsymbol \xi - \boldsymbol z \rangle_\delta^{2-m} \in H^1_0(B_r)\).
Integrating by parts twice gives
\[
  \begin{aligned}
  &\int_{B_r} \varpi (\Gamma - \bar \Gamma_\rho)
     \frac{\eta_\rho^2}{\langle \boldsymbol \xi - \boldsymbol z \rangle_\delta^{m-2}} \, \mathrm d \boldsymbol \xi \\
  ={}& \int_{B_r} \frac{\eta_\rho^2 |\nabla \Gamma|^2}{\langle \boldsymbol \xi - \boldsymbol z \rangle_\delta^{m-2}} \, \mathrm d \boldsymbol \xi
     + \frac{1}{2} \int_{B_r} \nabla \bigl[(\Gamma - \bar \Gamma_\rho)^2\bigr]
       \cdot \nabla \left(\frac{\eta_\rho^2}{\langle \boldsymbol \xi - \boldsymbol z \rangle_\delta^{m-2}}\right) \mathrm d \boldsymbol \xi \\
  ={}& \int_{B_r} \frac{\eta_\rho^2 |\nabla \Gamma|^2}{\langle \boldsymbol \xi - \boldsymbol z \rangle_\delta^{m-2}} \, \mathrm d \boldsymbol \xi
     + \frac{1}{2} \int_{B_r} (\Gamma - \bar \Gamma_\rho)^2
       \left[-\Delta \left(\frac{\eta_\rho^2}{\langle \boldsymbol \xi - \boldsymbol z \rangle_\delta^{m-2}}\right)\right] \mathrm d \boldsymbol \xi.
  \end{aligned}
\]
Direct differentiation gives
\[
  -\Delta \bigl(\langle \boldsymbol \xi - \boldsymbol z \rangle_\delta^{2-m}\bigr)
  = \frac{m(m - 2) \delta^2}{\langle \boldsymbol \xi - \boldsymbol z \rangle_\delta^{m+2}} \ge 0
  \qquad \text{on } \mathbb R^m,
\]
while all terms in \(-\Delta \bigl(\eta_\rho^2 \langle \boldsymbol \xi - \boldsymbol z \rangle_\delta^{2-m}\bigr)\) in which at least one derivative hits \(\eta_\rho^2\) are supported in \(B_{2\rho}(\boldsymbol z) \setminus B_\rho(\boldsymbol z)\) and bounded in modulus by \(C \rho^{-m}\).
Dropping the nonnegative Laplacian term on \(B_r\) and using \(\|\Gamma - \bar \Gamma_\rho\|_{L^\infty(B_r)} \le 2 \|\Gamma\|_{L^\infty(B_r)} \le C\) from \eqref{eq:GammaLinfty}, we apply Poincar\'e's inequality on the rescaled annulus \(B_{2\rho}(\boldsymbol z) \setminus B_\rho(\boldsymbol z)\) (where \(\rho \ge \frac{1}{2} |\boldsymbol \xi - \boldsymbol z|\)):
\[
  \begin{aligned}
  \int_{B_r} \frac{\eta_\rho^2 |\nabla \Gamma|^2}{\langle \boldsymbol \xi - \boldsymbol z \rangle_\delta^{m-2}} \, \mathrm d \boldsymbol \xi
  &\le C \int_{B_{2\rho}(\boldsymbol z)} \frac{|\varpi(\boldsymbol \xi)|}{|\boldsymbol \xi - \boldsymbol z|^{m-2}} \, \mathrm d \boldsymbol \xi
       + \frac{C}{\rho^m} \int_{B_{2\rho}(\boldsymbol z) \setminus B_\rho(\boldsymbol z)} |\Gamma - \bar \Gamma_\rho|^2 \, \mathrm d \boldsymbol \xi \\
  &\le C \int_{B_{2\rho}(\boldsymbol z)} \frac{|\varpi(\boldsymbol \xi)|}{|\boldsymbol \xi - \boldsymbol z|^{m-2}} \, \mathrm d \boldsymbol \xi
       + C \int_{B_{2\rho}(\boldsymbol z) \setminus B_\rho(\boldsymbol z)} \frac{|\nabla \Gamma(\boldsymbol \xi)|^2}{|\boldsymbol \xi - \boldsymbol z|^{m-2}} \, \mathrm d \boldsymbol \xi.
  \end{aligned}
\]
Sending \(\delta \downarrow 0\) via Fatou's lemma on \(B_\rho(\boldsymbol z)\), adding the annular integral of \(|\nabla \Gamma|^2 |\boldsymbol \xi - \boldsymbol z|^{2-m}\) to both sides, and substituting \eqref{eq:holevarpi} yields
\begin{equation}\label{eq:holeGamma}
  \int_{B_{2\rho}(\boldsymbol z)} \frac{|\nabla \Gamma(\boldsymbol \xi)|^2}{|\boldsymbol \xi - \boldsymbol z|^{m-2}} \, \mathrm d \boldsymbol \xi
  \le C \rho^2
      + C \int_{B_{2\rho}(\boldsymbol z) \setminus B_\rho(\boldsymbol z)}
          \frac{|\nabla \Gamma(\boldsymbol \xi)|^2 + |\varpi(\boldsymbol \xi)| + |\mathbf v(\boldsymbol \xi)|^2}
               {|\boldsymbol \xi - \boldsymbol z|^{m-2}} \, \mathrm d \boldsymbol \xi.
\end{equation}

Third, we test \eqref{eq:localenergy} with \(\varphi(\boldsymbol \xi) = \eta_\rho(\boldsymbol \xi)^2 \langle \boldsymbol \xi - \boldsymbol z \rangle_\delta^{4-m}\) (\(\delta > 0\)).
Since \(-\Delta \langle \boldsymbol \xi - \boldsymbol z \rangle_\delta^{4-m} \ge 2(m - 4) \langle \boldsymbol \xi - \boldsymbol z \rangle_\delta^{2-m}\) and \(\operatorname{supp} \nabla \eta_\rho \subset B_{2\rho}(\boldsymbol z) \setminus B_\rho(\boldsymbol z)\), integrating by parts on the left-hand side of \eqref{eq:localenergy} yields
\[
  \begin{aligned}
  &\int_{B_r} \nabla \mathbf v : \nabla(\mathbf v \varphi) \, \mathrm d \boldsymbol \xi \\
  ={}& \int_{B_r} \left(|\nabla \mathbf v|^2 \varphi + \frac{1}{2} |\mathbf v|^2 (-\Delta \varphi)\right) \, \mathrm d \boldsymbol \xi \\
  \ge{}& \int_{B_r} \eta_\rho^2 \Bigl(
        \frac{|\nabla \mathbf v|^2}{\langle \boldsymbol \xi - \boldsymbol z \rangle_\delta^{m-4}}
        + \frac{(m - 4) |\mathbf v|^2}{\langle \boldsymbol \xi - \boldsymbol z \rangle_\delta^{m-2}}
      \Bigr) \, \mathrm d \boldsymbol \xi
  - C \int_{B_{2\rho}(\boldsymbol z) \setminus B_\rho(\boldsymbol z)}
        \frac{|\mathbf v|^2}{|\boldsymbol \xi - \boldsymbol z|^{m-2}} \, \mathrm d \boldsymbol \xi.
  \end{aligned}
\]
On the right-hand side of \eqref{eq:localenergy}, substituting \(\varpi = -\Delta \Gamma\) and integrating by parts gives
\[
  \int_{B_r} \varpi \mathbf v \cdot \nabla \varphi \, \mathrm d \boldsymbol \xi
  = \sum_{i=1}^m \sum_{j=1}^m \int_{B_r}
      \partial_j \Gamma \bigl((\partial_j v_i) \partial_i \varphi + v_i \partial_i \partial_j \varphi\bigr) \, \mathrm d \boldsymbol \xi.
\]
Using \(0 \le \eta_\rho \le 1\), we have
\[
  \begin{aligned}
  |\nabla \varphi|
  &\le C \frac{\eta_\rho}{\langle \boldsymbol \xi - \boldsymbol z \rangle_\delta^{m-3}}
     + C \frac{\mathbf 1_{B_{2\rho}(\boldsymbol z) \setminus B_\rho(\boldsymbol z)}}{\rho^{m-3}}, \\
  |\nabla^2 \varphi|
  &\le C \frac{\eta_\rho}{\langle \boldsymbol \xi - \boldsymbol z \rangle_\delta^{m-2}}
     + C \frac{\mathbf 1_{B_{2\rho}(\boldsymbol z) \setminus B_\rho(\boldsymbol z)}}{\rho^{m-2}}.
  \end{aligned}
\]
Splitting the singular exponents via \(\frac{m - 4}{2} + \frac{m - 2}{2} = m - 3\) and \(\frac{m - 2}{2} + \frac{m - 2}{2} = m - 2\) in the main terms as
\[
  \begin{aligned}
  \frac{\eta_\rho |\nabla \mathbf v| |\nabla \Gamma|}{\langle \boldsymbol \xi - \boldsymbol z \rangle_\delta^{m-3}}
  &= \left(\frac{\eta_\rho |\nabla \mathbf v|}{\langle \boldsymbol \xi - \boldsymbol z \rangle_\delta^{(m-4)/2}}\right)
     \left(\frac{|\nabla \Gamma|}{\langle \boldsymbol \xi - \boldsymbol z \rangle_\delta^{(m-2)/2}}\right), \\
  \frac{\eta_\rho |\mathbf v| |\nabla \Gamma|}{\langle \boldsymbol \xi - \boldsymbol z \rangle_\delta^{m-2}}
  &= \left(\frac{\eta_\rho |\mathbf v|}{\langle \boldsymbol \xi - \boldsymbol z \rangle_\delta^{(m-2)/2}}\right)
     \left(\frac{|\nabla \Gamma|}{\langle \boldsymbol \xi - \boldsymbol z \rangle_\delta^{(m-2)/2}}\right),
  \end{aligned}
\]
applying Young's inequality, and using \(\langle \boldsymbol \xi - \boldsymbol z \rangle_\delta^{2-m} \le |\boldsymbol \xi - \boldsymbol z|^{2-m}\) on the \(\nabla \Gamma\) terms yields
\[
  \begin{aligned}
  &\int_{B_r} \varpi \mathbf v \cdot \nabla \varphi \, \mathrm d \boldsymbol \xi \\
  \le{}& \frac{1}{2} \int_{B_r} \eta_\rho^2 \Bigl(
         \frac{|\nabla \mathbf v|^2}{\langle \boldsymbol \xi - \boldsymbol z \rangle_\delta^{m-4}}
         + \frac{(m - 4) |\mathbf v|^2}{\langle \boldsymbol \xi - \boldsymbol z \rangle_\delta^{m-2}}
       \Bigr) \, \mathrm d \boldsymbol \xi + C \int_{B_{2\rho}(\boldsymbol z)}
             \frac{|\nabla \Gamma|^2}{|\boldsymbol \xi - \boldsymbol z|^{m-2}} \, \mathrm d \boldsymbol \xi \\
  &+ C \int_{B_{2\rho}(\boldsymbol z) \setminus B_\rho(\boldsymbol z)}
         \left(\frac{|\nabla \mathbf v|^2}{|\boldsymbol \xi - \boldsymbol z|^{m-4}}
               + \frac{|\mathbf v|^2}{|\boldsymbol \xi - \boldsymbol z|^{m-2}}\right) \, \mathrm d \boldsymbol \xi.
  \end{aligned}
\]
We absorb the first integral into the left-hand side at fixed \(\delta > 0\) and send \(\delta \downarrow 0\) via Fatou's lemma.
Bounding the \(\nabla \Gamma\) integral via \eqref{eq:holeGamma} and adding \eqref{eq:holevarpi}--\eqref{eq:holeGamma} restricted to \(B_\rho(\boldsymbol z)\) then yields a constant \(c_0 > 0\) depending only on \(m\) and \(M\) such that, for all \(\boldsymbol z \in B_{r_0}\) and \(0 < \rho \le \rho_0\),
\[
  I(\rho)
  \le c_0 \int_{B_{2\rho}(\boldsymbol z) \setminus B_\rho(\boldsymbol z)}
      \left(\frac{|\nabla \mathbf v(\boldsymbol \xi)|^2}{|\boldsymbol \xi - \boldsymbol z|^{m-4}}
            + \frac{|\mathbf v(\boldsymbol \xi)|^2 + |\nabla \Gamma(\boldsymbol \xi)|^2 + |\varpi(\boldsymbol \xi)|}
                   {|\boldsymbol \xi - \boldsymbol z|^{m-2}}\right) \, \mathrm d \boldsymbol \xi
      + c_0 \rho^2.
\]
At \(\rho = \rho_0\), since \(|\boldsymbol \xi - \boldsymbol z| \ge \rho_0\) on \(B_{2\rho_0}(\boldsymbol z) \setminus B_{\rho_0}(\boldsymbol z)\), this estimate immediately yields \(\sup_{\boldsymbol z \in B_{r_0}} I(\rho_0) \le C < \infty\).
For \(0 < \rho \le \rho_0 / 2\), the annular integral equals \(I(2\rho) - I(\rho)\), giving Widman's hole-filling inequality \cite{Giaquinta1983}
\[
  I(\rho) \le c_0 \bigl[I(2\rho) - I(\rho)\bigr] + c_0 \rho^2,
  \qquad \text{whence} \qquad
  I(\rho) \le \frac{c_0}{1 + c_0} I(2\rho) + c_0 \rho^2.
\]
Iterating dyadically from \(\rho_0 / 2\) yields an exponent \(\sigma \in (0, 1)\) depending only on \(m\) and \(M\) such that \(I(\rho) \le C \rho^\sigma\) for all \(\boldsymbol z \in B_{r_0}\) and \(0 < \rho \le \rho_0 / 2\).
Bounding the range \(\rho_0 / 2 < \rho \le 2 r_0\) directly by \((\rho_0 / 2)^{-(m - 4 + \sigma)} \|\nabla \mathbf v\|_{L^2(B_2)}^2\), we conclude that
\begin{equation}\label{eq:strictmorrey}
  \|\nabla \mathbf v\|_{L^{2, m - 4 + \sigma}(B_{r_0})} \le C.
\end{equation}

\smallskip
\noindent\emph{Step 3.}
Select intermediate radii \(1 < r_1 < r_2 < r_0\) and set \(\mu := m - 4 + \sigma\), so that \(m - \mu = 4 - \sigma \in (3, 4)\).
For any concentric balls \(B_{\rho_1} \Subset B_{\rho_2} \subset B_{r_0}\), Adams' Morrey embedding \cite{Adams75} (see \cite[Appendix~A, Theorem~A.2]{LY22}) gives
\[
  \|\mathbf f\|_{L^{q^*,\mu}(B_{\rho_1})}
  \le C_q \bigl(\|\nabla \mathbf f\|_{L^{q,\mu}(B_{\rho_2})} + \|\mathbf f\|_{L^1(B_{\rho_2})}\bigr),
  \qquad 1 < q < 4 - \sigma, \quad \frac{1}{q^*} = \frac{1}{q} - \frac{1}{4 - \sigma},
\]
whenever \(\mathbf f \in W^{1,q}(B_{\rho_2}) \cap L^{q,\mu}(B_{\rho_2})\) with \(\nabla \mathbf f \in L^{q,\mu}(B_{\rho_2})\).
Moreover, let \((\mathbf v, \pi)\) be a distributional solution of \(-\Delta \mathbf v + \nabla \pi = \mathbf g\) and \(\operatorname{div} \mathbf v = 0\) with \(\mathbf g \in L^{q,\mu}(B_{\rho_2})\) (\(1 < q < \infty\)).
Mollifying at scales \(0 < \delta < (\rho_2 - \rho_1)/2\), applying the linear Morrey--Stokes estimate \cite[Appendix~A, Theorem~A.3]{LY22} via a finite covering of \(B_{\rho_1}\) by rescaled balls inside \(B_{(\rho_1+\rho_2)/2}\), and letting \(\delta \downarrow 0\) via weak \(L^q(B_{\rho_1})\) lower semicontinuity gives
\[
  \|\nabla^2 \mathbf v\|_{L^{q,\mu}(B_{\rho_1})} + \|\nabla \pi\|_{L^{q,\mu}(B_{\rho_1})}
  \le C_q \bigl(\|\mathbf g\|_{L^{q,\mu}(B_{\rho_2})} + \|\mathbf v\|_{L^1(B_{\rho_2})}\bigr).
\]

By \eqref{eq:B35pot} (since \(\mu < m - 2\)) and \eqref{eq:strictmorrey}, we have \(\mathbf v, \nabla \mathbf v \in L^{2,\mu}(B_{r_0})\).
Choosing an intermediate radius \(\tilde r_0 \in (r_2, r_0)\), Adams' embedding at \(q = 2\) yields \(\mathbf v \in L^{s_0,\mu}(B_{\tilde r_0})\), with
\[
  \frac{1}{s_0} := \frac{1}{2} - \frac{1}{4 - \sigma},
  \qquad\text{so that}\qquad
  s_0 = \frac{2(4 - \sigma)}{2 - \sigma} > 4 - \sigma.
\]
Starting from \(p_0 = 2\), whenever \(p_j < 4 - \sigma\), H\"older's inequality shows that \(\mathbf g = -(\mathbf v \cdot \nabla) \mathbf v \in L^{q_j,\mu}\), and the Morrey--Stokes estimate followed by Adams' embedding yields \(\nabla \mathbf v \in L^{p_{j+1},\mu}\) on a slightly smaller concentric ball, where
\[
  \frac{1}{q_j} = \frac{1}{p_j} + \frac{1}{s_0} \in \left(\frac{1}{4 - \sigma}, \, 1\right),
  \qquad
  \frac{1}{p_{j+1}} = \frac{1}{q_j} - \frac{1}{4 - \sigma}
  = \frac{1}{p_j} - \left(\frac{2}{4 - \sigma} - \frac{1}{2}\right).
\]
Since \(\frac{2}{4 - \sigma} - \frac{1}{2} > 0\), a finite number \(J = J(m, M)\) of equally spaced concentric steps from \(B_{\tilde r_0}\) to \(B_{r_2}\) reaches \(p_J \ge 4 - \sigma\), giving \(\nabla \mathbf v \in L^{p_J,\mu}(B_{r_2})\) together with \(\mathbf v \in L^{s_0,\mu}(B_{r_2})\).

Because \(\min\{p_J, s_0\} \ge 4 - \sigma\), H\"older's inequality gives \(\mathbf v, \nabla \mathbf v \in L^{q,\mu}(B_{r_2})\) for every \(q \in (2, 4 - \sigma)\), and Adams' embedding on \(B_{r_1} \Subset B_{r_2}\) yields \(\mathbf v \in L^{q^*,\mu}(B_{r_1}) \subset L^{q^*}(B_{r_1})\).
For any \(s \in [1, \infty)\), choosing \(q \in (2, 4 - \sigma)\) with \(q^* \ge s\) yields \(\|\mathbf v\|_{L^s(B_{r_1})} \le C_{m,M,s}\).
In particular, taking \(s > 2m\) gives \(\mathbf v \otimes \mathbf v \in L^{s/2}(B_{r_1})\) with \(s / 2 > m\).
Linear \(W^{1,s/2}\) interior Stokes regularity \cite{Galdi2011} applied to \(-\Delta \mathbf v + \nabla \pi = -\operatorname{div}(\mathbf v \otimes \mathbf v)\) then gives \(\mathbf v \in W^{1,s/2}_{\mathrm{loc}}(B_{r_1}) \hookrightarrow L^\infty_{\mathrm{loc}}(B_{r_1})\), and standard linear Stokes bootstrapping across concentric balls between \(B_{r_1}\) and \(B_1\) establishes \eqref{eq:quantitativeCk}.
\end{proof}

\addtocontents{toc}{\SkipTocEntry}\section*{Acknowledgments}
Baishun Lai was partially supported by the National Natural Science Foundation of China under Grant No.~12371202.
Ruiyi Zhang and Xiaoxin Zheng were partially supported by the National Natural Science Foundation of China under Grant No.~12371231.

\medskip
\noindent\textbf{Conflict of interest.}
The authors declare that they have no conflict of interest.

\medskip
\noindent\textbf{Data availability.}
Data sharing is not applicable to this article as no datasets were generated or analyzed during the current study.

\medskip
\noindent\textbf{Declaration of generative AI use.}
During the preparation of this work, generative AI tools were used to assist with calculation checks, literature organization, and language editing.
The authors verified all mathematical arguments and take full responsibility for the content of the paper.

\end{document}